\documentclass[11pt, A4paper]{article}

\usepackage{amsmath}
\usepackage{amssymb}
\usepackage{graphicx}
\usepackage{epstopdf}
\usepackage[compatibility=false]{caption}
\usepackage{subcaption}
\usepackage{float}
\usepackage{placeins}

\newtheorem{prop}{Proposition}

\newtheorem{lem}{Lemma}

\newtheorem{remark}{Remark}

\title{A fluid approach to total-progeny-dependent birth-and-death processes}
\author{Sophie Hautphenne and Minyuan Li\footnote{Corresponding author: minyuanl@student.unimelb.edu.au}\\The University of Melbourne}
\date{}

\begin{document}
\maketitle

\begin{abstract}
We introduce a class of branching processes in which the reproduction or lifetime distribution at a given time depends on the total cumulative number of individuals who have been born in the population until that time. 
We focus on a
continuous-time version of these processes, called \emph{total-progeny-dependent birth-and-death processes}, and study some of their properties through the analysis of their fluid (deterministic) approximation. These properties include the maximum population size, the total progeny size at extinction, the time to reach the maximum population size, and the time until extinction. As the fluid approach does not allow us to approximate the time until extinction directly, we propose several methods to complement this approach. We also use the fluid approach to study the behaviour of the processes as we increase the magnitude of the individual birth rate. 
\medskip

\end{abstract}

\section{Introduction}

Many resources are not renewable. Examples affecting human populations include fossil fuels (coal, petroleum, natural gas) and nuclear energy. Animal populations also often modify their environment by building burrows or nests, eating other living beings, etc. Any population that uses resources at unsustainable rates risks a severe decline or extinction. 
In this paper, we introduce a class of  branching processes  equipped to model the negative impact of non-renewable resource consumption. 
More specifically, we assume that
each individual's lifetime or reproductive success depends on the current \textit{total progeny size,} that is, the total cumulative number of individuals who have been born in the population. It is reasonable, for example, to define the individual fertility rate as a decreasing function of the total progeny size since the more individuals have been born in the population, the more resources have been consumed, and the less propitious the conditions for reproduction have become. Assuming a one-to-one correspondence between the total progeny size and the amount of consumed resources, the level of resources is then not modelled directly, but rather indirectly through an endogenous feature of the population. The resulting process, which we name \emph{total-progeny-dependent branching process}, is a two-dimensional Markov chain which exhibits interesting features, as we will demonstrate.

\begin{figure}[h]
\centering \includegraphics[width=0.9\textwidth]{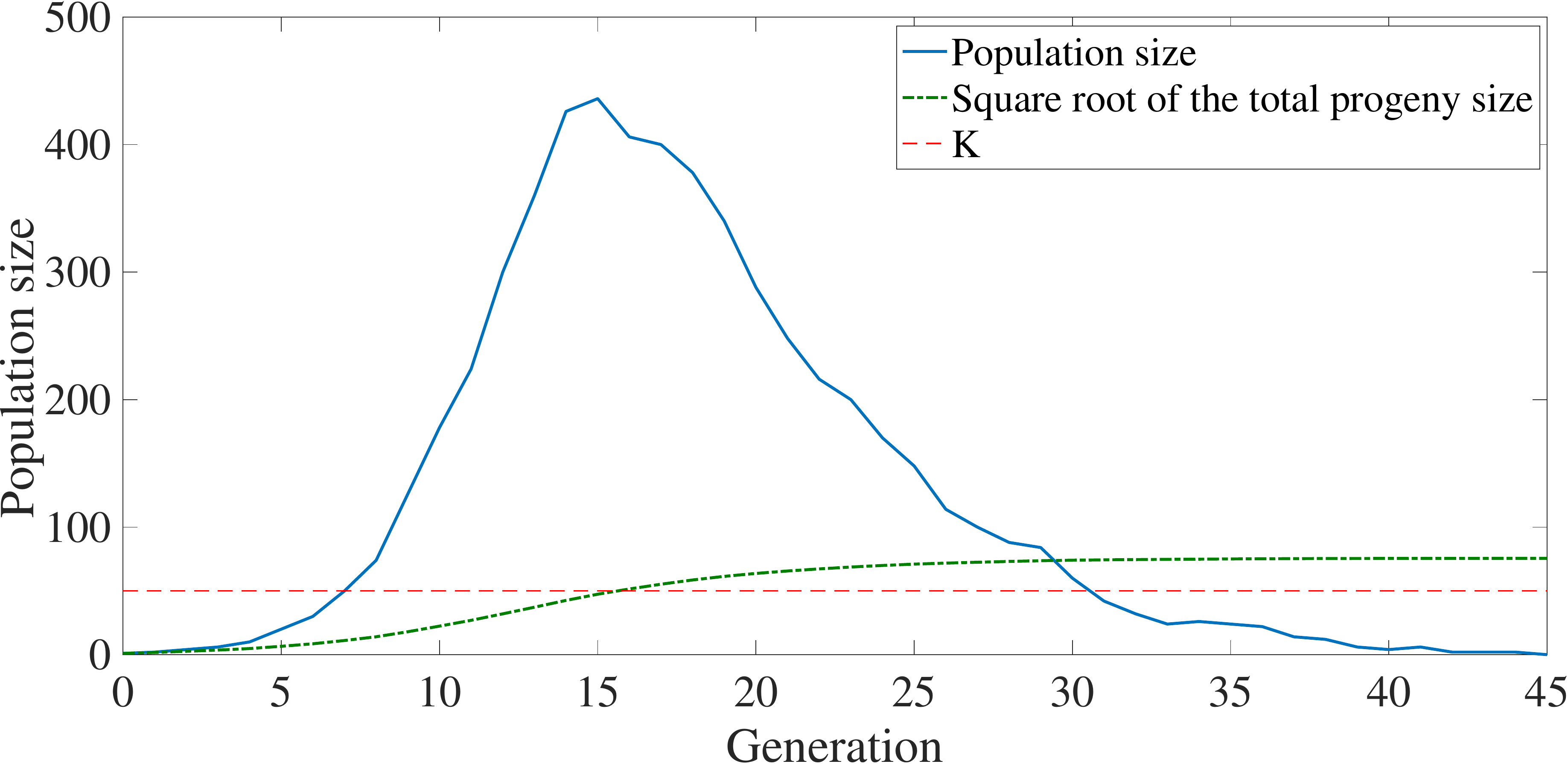}
 \caption{\label{f4}A trajectory of a total-progeny-dependent branching process with binary splitting where $K=50$.}
\end{figure}

An example of a typical trajectory of a discrete-time total-progeny-dependent branching process is shown in Figure \ref{f4}. In this binary splitting model, at each generation, individuals produce two children with probability $p_2(x)=K/(K+\sqrt{x})$ (where $K$ is a positive constant), or no children with probability $p_0(x)=1-p_2(x)$, when the current total progeny size is $x$ ($\sqrt{x}$ is represented by the green curve). We see that the population starts declining as soon as $\sqrt{x}$ crosses the constant $K$, that is, when the process moves from being supercritical to being subcritical. 

\begin{figure}[h]
\centering\includegraphics[width=0.9\textwidth]{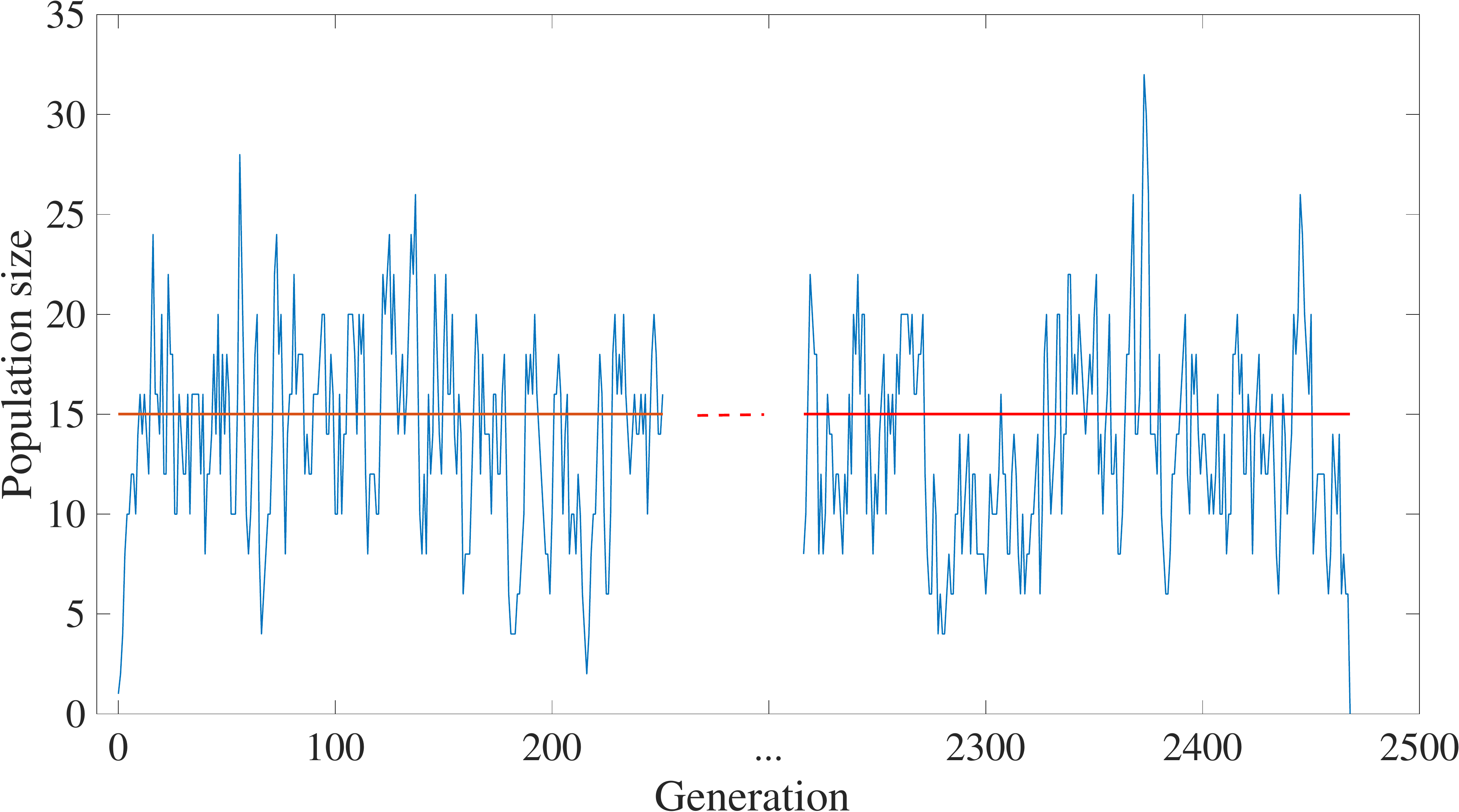}
 \caption{\label{f5}A portion of trajectory of a population-size-dependent branching processes with binary splitting and {carrying capacity} $K=15$.}
\end{figure}

In a setting where resources do not renew, populations are never able to sustain themselves, and ultimate extinction happens not long after the population size has reached a peak. A contrasting class of well-studied models with dependences are the so-called \emph{population-size-dependent branching processes}, in which the individual's lifetime or reproductive success depends on the current \emph{population size}. In these processes, the population can grow until it reaches the \textit{maximum population size} a particular habitat  can support, called the \textit{carrying capacity}. If it is able to reach the vicinity of the carrying capacity, then the population lingers around that value for a very long time ---more precisely, an expected time of the order  $e^{cK}$ for some constant $c>0$, where $K$ is the carrying capacity \cite{hamza2016establishment}. Due to the presence of the carrying capacity and stochastic events, extinction eventually occurs with probability one \cite{jagers1992stabilities}.
 These well-established asymptotic properties are common for a large class of processes, see \cite{hamza2016establishment,jk11,jagers2020populations}. Figure \ref{f5} shows an example of a typical trajectory of a binary splitting discrete-time population-size-dependent branching processes in which $p_2(z)=K/(K+z)$, where $z$ is the current population size. We observe that, despite the similar form of the offspring laws, the trajectories in Figures \ref{f4} and \ref{f5} are completely different. In particular, there is no concept of carrying capacity in a total-progeny-dependent branching process.

To the best of our knowledge, total-progeny-dependent branching processes have not previously appeared in the literature. It is worth mentioning two notable related models. First, the stochastic Susceptible-Infectious-Removed (SIR) epidemic process is a special case of total-progeny-dependent birth-and-death process where the population size has an upper bound, and the individual birth rate is a linearly decreasing function of the total progeny size (see Remark \ref{remSIR} for more detail). The class of models considered here is much richer since the birth and death rates can be any (positive) function of the total progeny, and there is no bound on the population size. Second, the \emph{resource-dependent branching process} introduced in \cite{bruss2015resource} is a discrete-time branching process in which individuals produce and consume resources; a society rule determines how resources are distributed among individuals, who have a means of interaction through claims: they stay in the society only if their claim is met, otherwise they leave. In contrast with the models considered here, the reproduction law of individuals in \cite{bruss2015resource} is constant over time. The resource-dependent branching processes may thus be seen as particular types of \emph{controlled branching processes} \cite{sevast1974controlled,velasco2017controlled}.

As a convenient starting point for the analysis of total-progeny-dependent branching processes, here we consider the simplest continuous-time process, namely, the \emph{total-progeny-dependent birth-and-death process}. In this process, the birth and death rates per individuals are some functions $b(x)$ and $d(x)$, respectively, of the current total progeny size $x$. 
It is generally known that the equivalent deterministic (fluid) models, which are easier to analyse, are for practical purposes as good as the stochastic models when the populations are sufficiently large \cite{kurtz1971limit}. We compare several properties of  total-progeny-dependent birth-and-death processes with their fluid equivalent, including the mean maximum population size, the mean total progeny at extinction, the mean time to reach the maximum population size, and the mean time until extinction. To this end, we focus on two toy models in which only the birth rate depends on the total progeny $x$: in Model 1, $b(x)=\lambda/x,$ and in Model~2, $b(x)= \lambda \exp(-{x}/{\alpha})$, where $\lambda$ and $\alpha$ are positive parameters which control the magnitude of the birth rate. We stress that these particular models are not chosen for realistic reasons, but rather because they are amenable to some explicit fluid analysis and they exhibit interesting properties. For example,  we show that, as $\lambda\to\infty$, the (fluid) time to reach the maximum population size converges to a positive constant in Model 1 and to zero in Model 2. Because the mean time until extinction cannot be directly approximated using the fluid approach, we propose several ways to complement that approach. We observe the existence of a value of $\lambda$ minimising the (approximated) mean time until extinction in Model 2. 

The models presented in this paper offer many potential extensions. Other more realistic models of total-progeny-dependent birth-and-death processes, in which the birth rate decreases more slowly with the total progeny size than in Models 1 and 2, include the cases where $b(x)=\lambda/x^p$ for $0\leq p\leq 1$, and  $b(x)=1-e^{-e^{-\gamma(x-\beta)}}$ (the negative Gompertz function); those models, however, are less tractable. It is also possible to relax the assumption that the resources are not renewable by supposing that the individual's lifetime and reproductive success at time $t$ depend on the total cumulative number of individuals who have been present in the population since time $(t-\delta)^+$ for some constant $\delta\geq0$ (where the particular case $
\delta=0$ corresponds to the population-size-dependent case).

The rest of the paper is organised as follows. In Section \ref{sec_prel} we define the total-progeny-dependent birth-and-death process and its deterministic equivalent, and we introduce some preliminary results. Sections \ref{sec_mod1} and  \ref{sec_mod2} are devoted to the analysis of properties of Models 1 and 2, respectively. Finally, Section \ref{proofs} contains the proofs of some results.

\section{Preliminaries}\label{sec_prel}
A continuous-time \emph{total-progeny-dependent birth-and-death process}   is a two-dimensional Markov process $\{(Z_t, X_t)\}_{t\geq 0}$, where $Z_t$ denotes the population size at time $t$ and $X_t$ denotes the total progeny until time $t$, in which the \emph{per individual} birth and death rates, $b(x)$ and $d(x)$, depend on the current value $x$ of the total progeny. For any fixed value of $x$, individuals reproduce and die independently of each other.
The two possible transitions from state $(z,x)\in \mathbb N_0^2$ are therefore
\begin{itemize}
\item $(z, x) \to (z+1, x+1)$, at rate $z b(x)$,
\item $(z, x) \to (z-1, x)$, at rate $z d(x)$.
\end{itemize}
Any state $(0,x)$ is absorbing.

The deterministic equivalent of the Markov process $\{(Z_t, X_t)\}$ with initial value $(Z_0,X_0)=(1,1)$, which we call its \emph{fluid approximation}, is the real-valued vector function $(y_1(t),y_2(t))$, $t\in \mathbb R_0^+$, which satisfies the system of ordinary differential equations
\begin{align}
\frac{dy_1(t)}{dt} &= y_1(t)\{b(y_2(t)) - d(y_2(t))\} \label{eq1}\\ 
\frac{dy_2(t)}{dt} &= y_1(t)b(y_2(t)), \label{eq2}
\end{align}
with initial conditions $(y_1(0),y_2(0))=(1,1)$. If the function $F(z,x):=(z\{b(x)-d(x)\},z b(x))$ is Lipschitz continuous, then a unique solution is guaranteed on $[0,T]$ for some value $T>0$ (by the Picard–Lindel\"{o}f theorem). 
The function $y_1(t)$  approximates the number of individuals in the population at time $t$, and the function $y_2(t)$ approximates the total progeny until time~$t$. 
Eq. \eqref{eq1} has solution
\begin{equation}\label{sol_eq1}y_1(t)=\exp\left(\int_0^t\{b(y_2(u)) - d(y_2(u))\} du\right)>0,\end{equation} which implies by Eq. \eqref{eq2} that $y_2(t)$ is monotonically increasing. If \begin{equation}\label{cond_ext}\lim_{t\to\infty}\int_0^t\{b(y_2(u)) - d(y_2(u))\} du= -\infty,\end{equation}then $y_1(t)\to 0$ (eventual \emph{extinction}). 

Note that we do not aim at formally deriving $(y_1(t),y_2(t))$ as the \emph{fluid limit} (as $N\to\infty$) of a sequence of Markov processes $\{(Z_t^N, X_t^N)\}_{t\geq 0, N\geq 1}$ under proper parametrisation and scaling so that the results in \cite{kurtz1970solutions,kurtz1971limit,darling2002fluid} can be applied. Instead, here we study the quality of the fluid approximation of some total-progeny-dependent birth-and-death processes and several of its properties to get insights in the processes for practical purposes.

\begin{remark}\label{remSIR}The stochastic SIR epidemic process is a particular total-progeny-dependent birth-and-death process in which the resources are the susceptible individuals in a closed population of size $N$. In this process, $Z_t$ corresponds to the number of infectious individuals at time $t$, $X_t$ corresponds to the total cumulative number of cases until time $t$ (that is, the sum of the number of infectious and recovered individuals),  $b(x):=\beta (1-x/N)$, and $d(x):=\gamma$, $0\leq x\leq N$, where  $\beta$ is the transmission rate and $\gamma$ is the recovery rate. The equivalence between the deterministic SIR model and the fluid limit of SIR stochastic processes is well known, see for example \cite{kurtz1971limit}.
\end{remark}

To simplify our analysis we consider two models  in which the death rate is constant, $d(x) = \mu$, and in which the birth rate decreases with the current value $x$ of the total progeny size,
\begin{itemize}
\item \textbf{Model 1: }$b_1(x) = \dfrac{\lambda}{x},$ $x\geq 1$,
\item \textbf{Model 2: }$b_2(x) = \lambda e^{-\frac{x}{\alpha}},$ $x\geq 1$,
\end{itemize}
where
$\lambda, \alpha, \mu$ are strictly positives real constants. In particular, $\lambda$ and $\alpha$ control the magnitude of the birth rate.

\begin{figure}[H]
\begin{subfigure}[h]{0.5\linewidth}
         \centering
         \includegraphics[width=\linewidth]{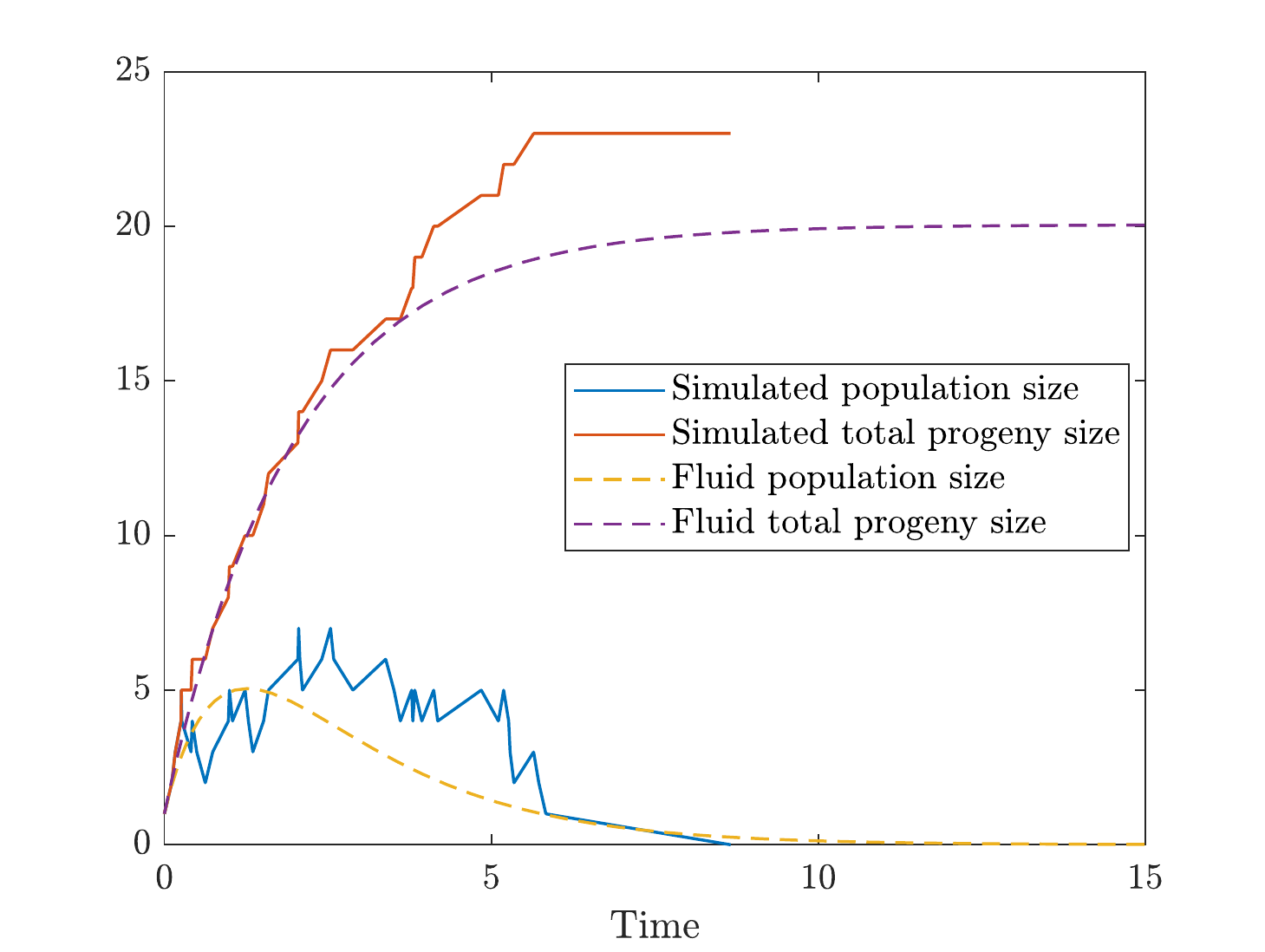}
         \caption{$\lambda=10$}
         \label{fig: simulation 100 vs ode}
     \end{subfigure}
     \hfill
     \begin{subfigure}[h]{0.5\linewidth}
         \centering
         \includegraphics[width=\textwidth]{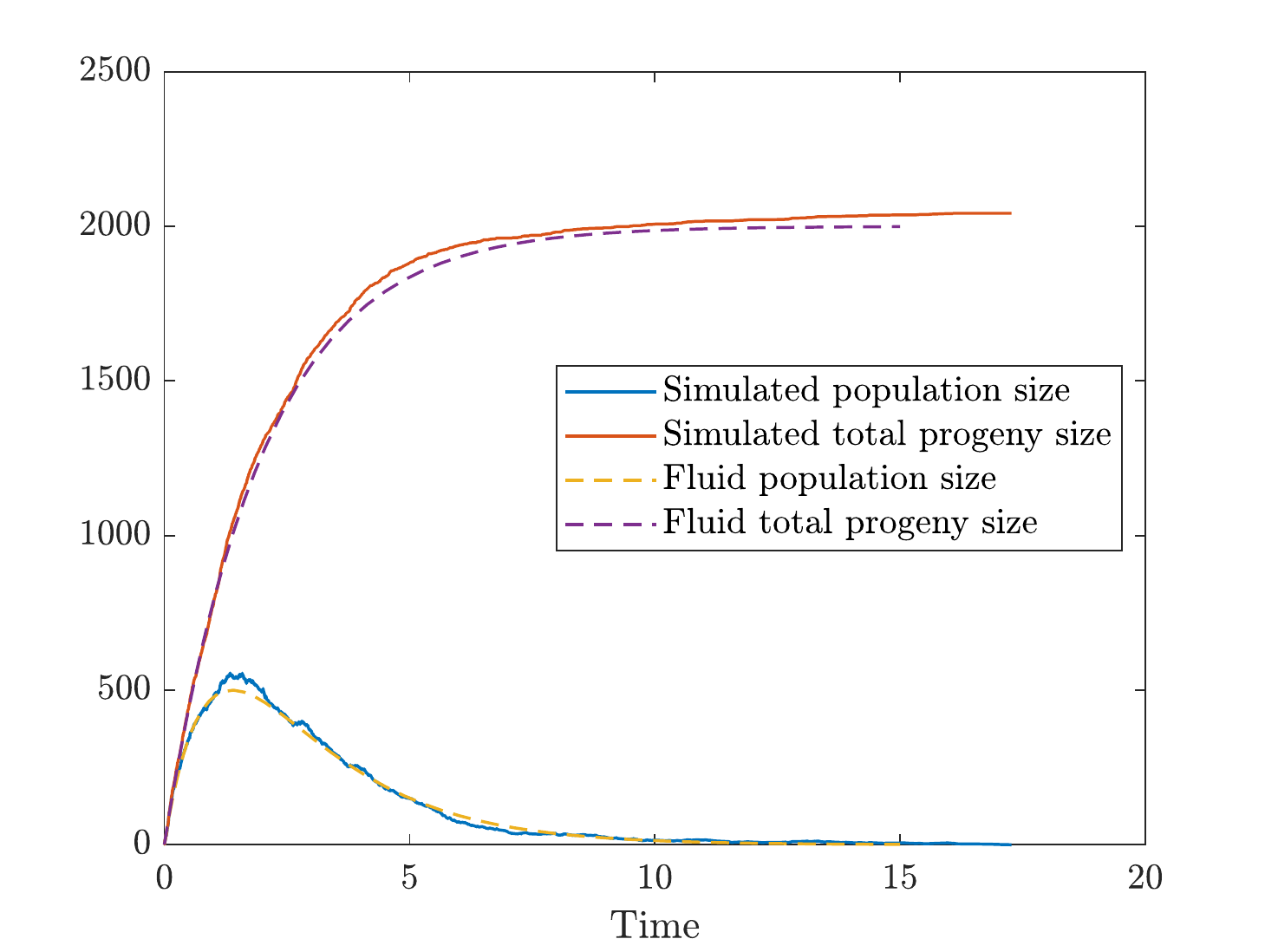}
         \caption{$\lambda=1000$}
         \label{fig: simulation 10000 vs ode}
     \end{subfigure}
	\caption{\textbf{Model 1 }--- Single trajectories of the stochastic process $\{(Z_t, X_t)\}$ and its fluid approximation $(y_1(t),y_2(t))$.}
	\label{fig: simulation vs ode}
\end{figure}

\begin{figure}[H]
\begin{subfigure}[h]{0.5\linewidth}
         \centering
          \includegraphics[width=\textwidth]{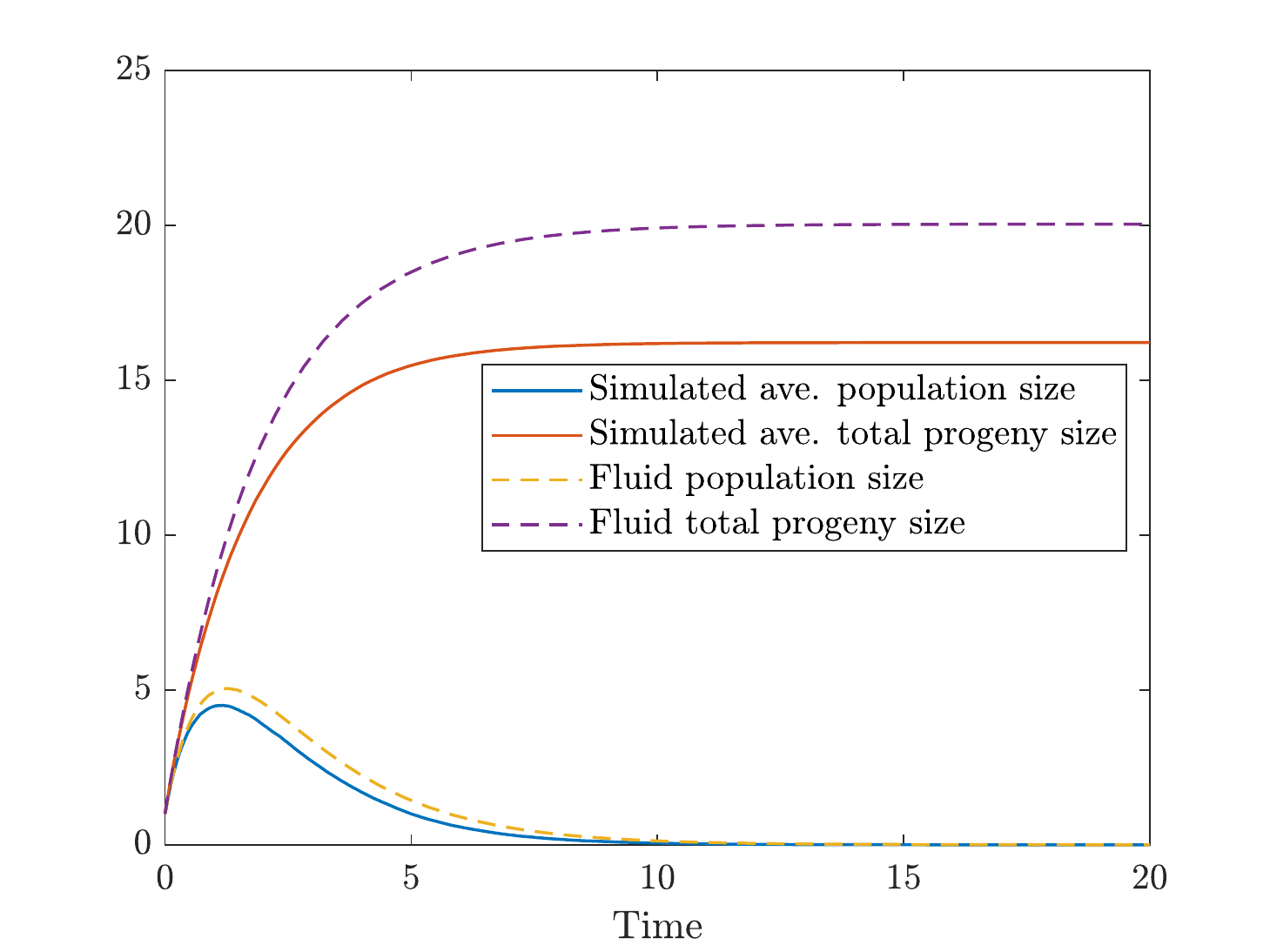}
         \caption{$\lambda=10$}
         \label{fig: ave simulation 100 vs ode}
       
     \end{subfigure}
     \hfill
     \begin{subfigure}[h]{0.5\linewidth}
         \centering
          \includegraphics[width=\linewidth]{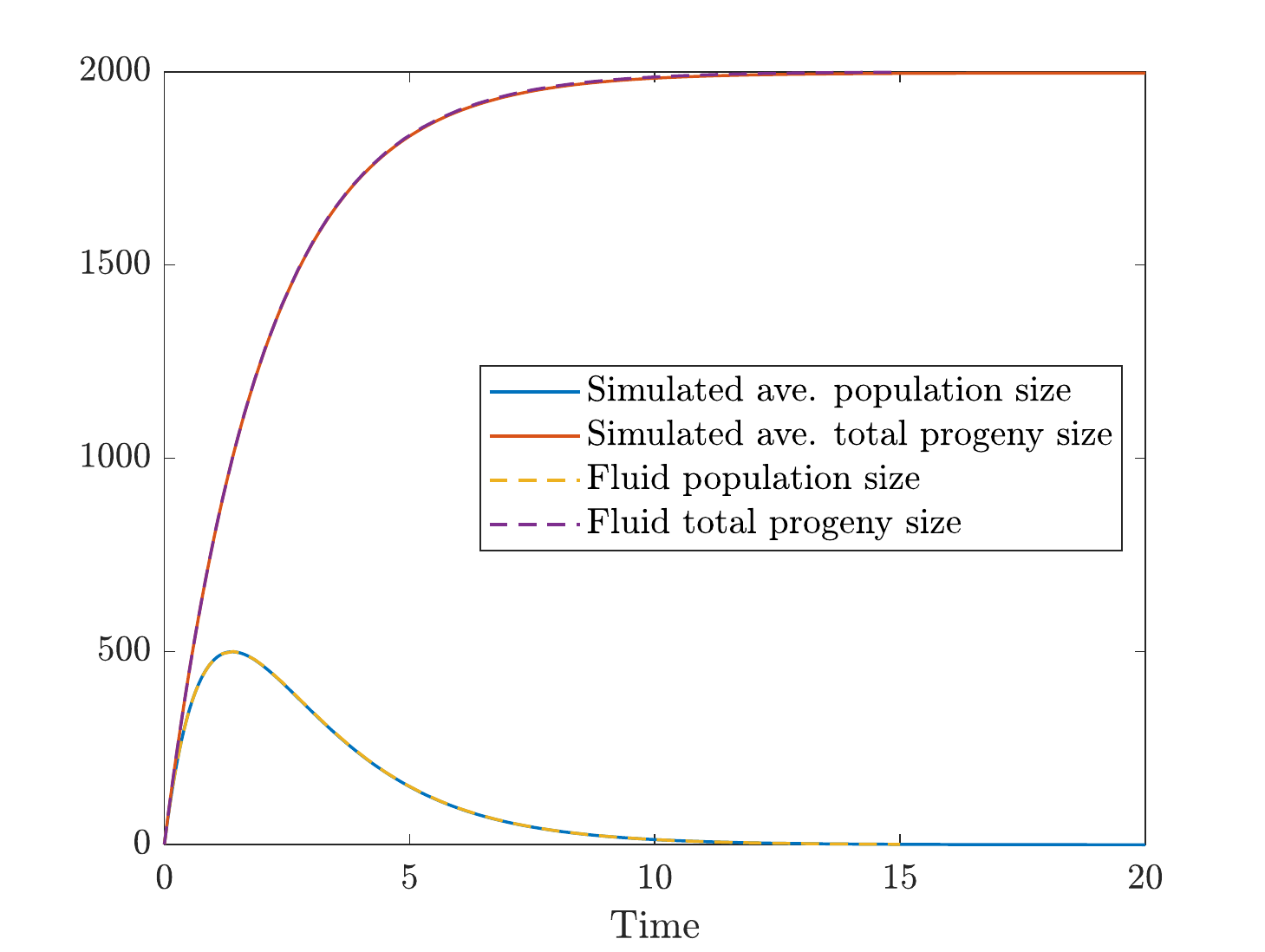}
         \caption{$\lambda=1000$}
         \label{fig: ave simulation 1000 vs ode}
     \end{subfigure}
	\caption{\textbf{Model 1 }--- Average of 10000 trajectories of the stochastic process $\{(Z_t, X_t)\}$ and its fluid approximation $(y_1(t),y_2(t))$.}
	\label{fig: ave simulation vs ode}
\end{figure}

\begin{figure}[H]
\begin{subfigure}[h]{0.5\linewidth}
         \centering
         \includegraphics[width=\linewidth]{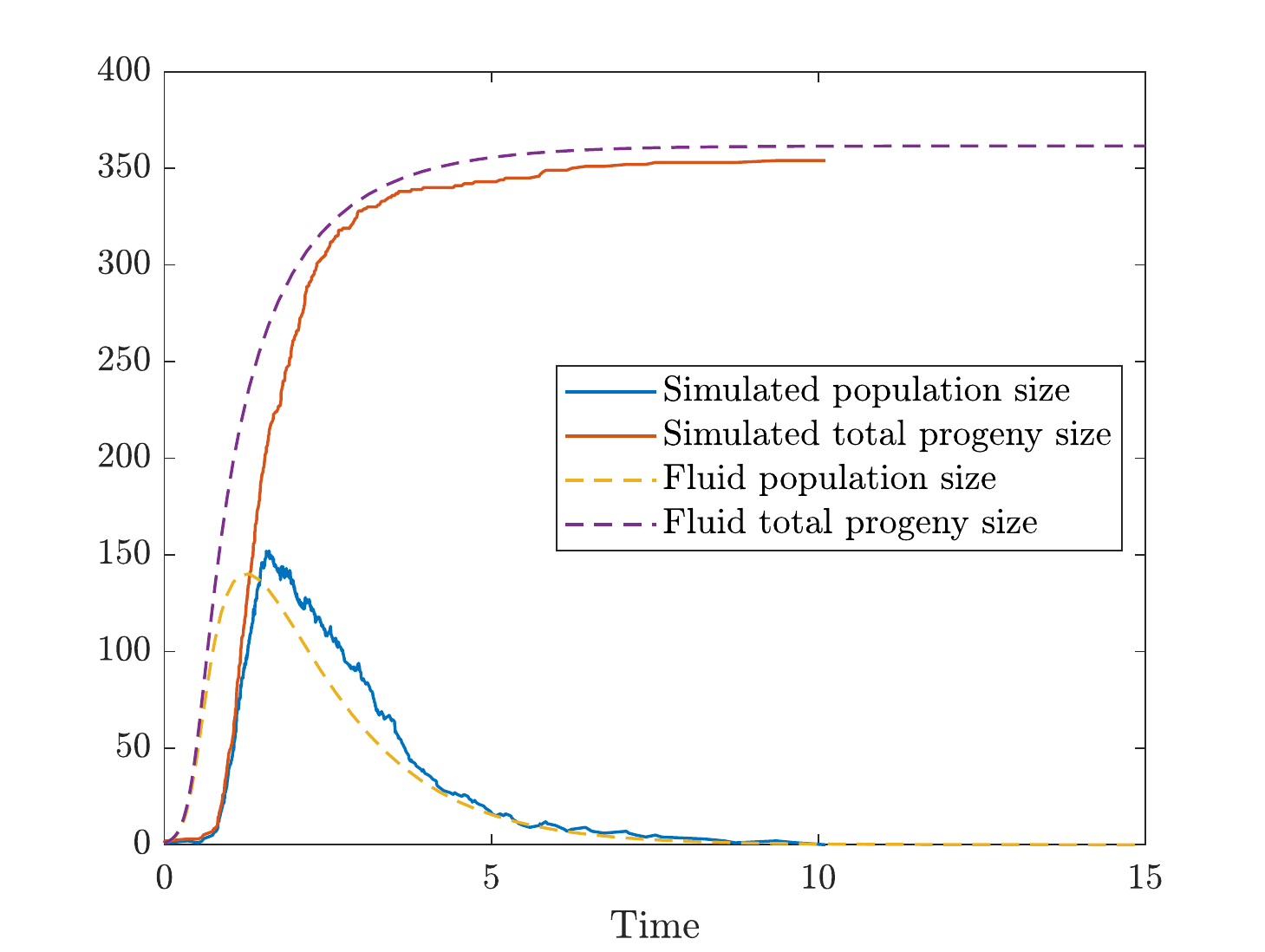}
         \caption{$\lambda=10$}
         \label{fig: simulation 100 vs ode}
     \end{subfigure}
     \hfill
     \begin{subfigure}[h]{0.5\linewidth}
         \centering
         \includegraphics[width=\textwidth]{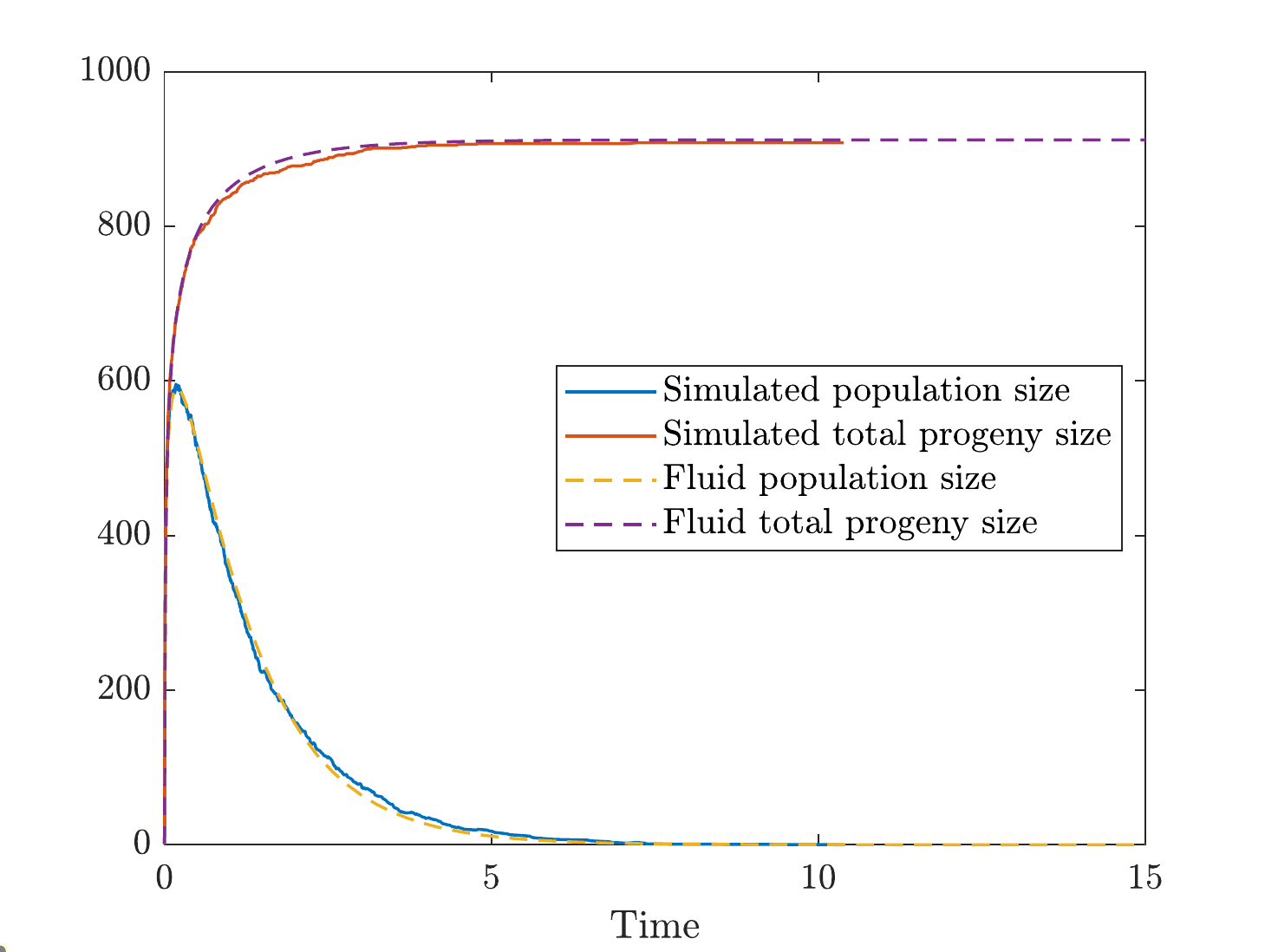}
         \caption{$\lambda=1000$}
         \label{fig: simulation 10000 vs ode}
     \end{subfigure}
	\caption{\textbf{Model 2 }--- Single trajectories of the stochastic process $\{(Z_t, X_t)\}$ and its fluid approximation $(y_1(t),y_2(t))$.}
	\label{fig: simulation vs ode m2}
\end{figure}

\begin{figure}[H]
\begin{subfigure}[h]{0.5\linewidth}
         \centering
          \includegraphics[width=\textwidth]{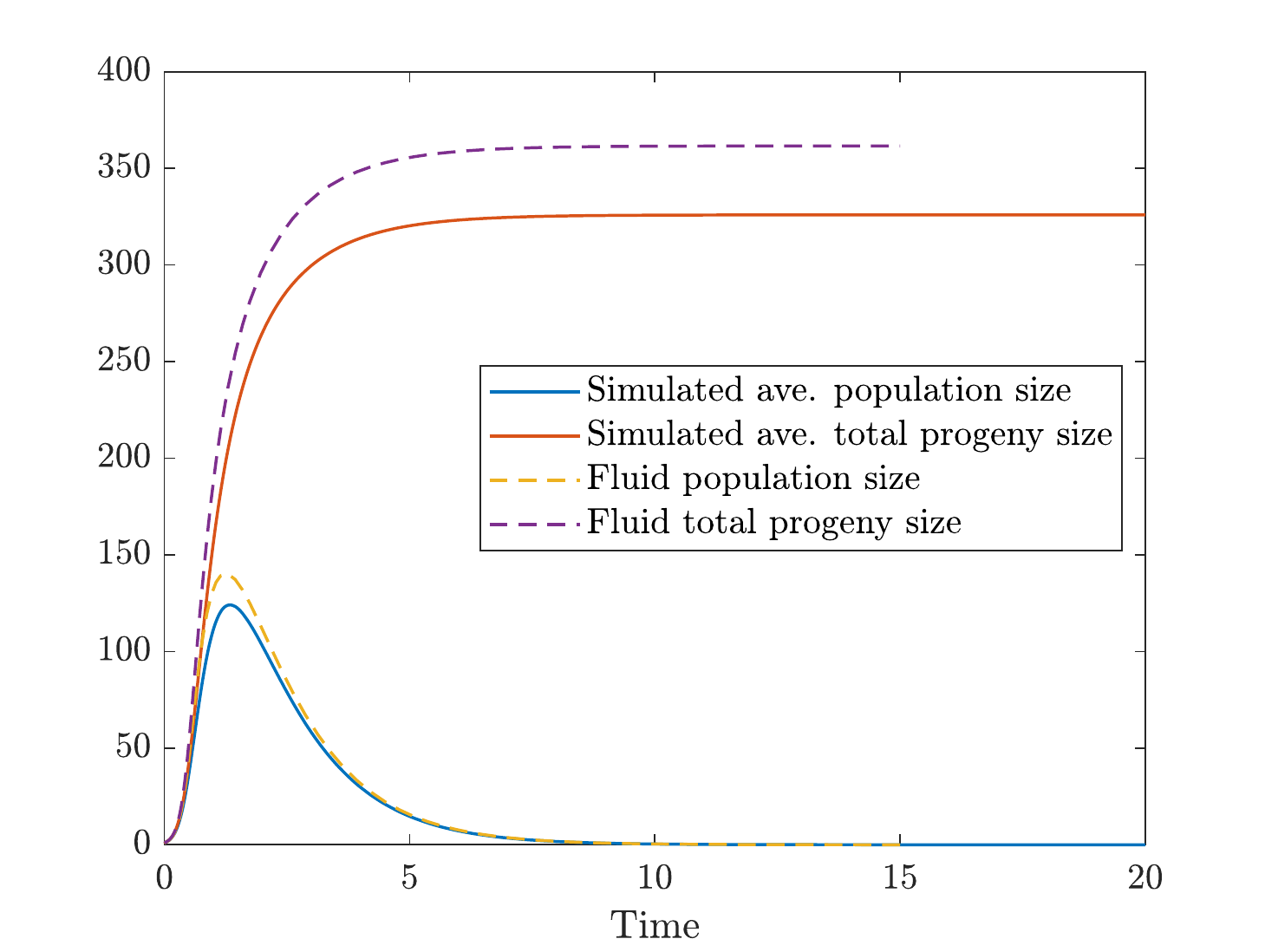}
         \caption{$\lambda=10$}
         \label{fig: ave simulation 100 vs ode}
       
     \end{subfigure}
     \hfill
     \begin{subfigure}[h]{0.5\linewidth}
         \centering
          \includegraphics[width=\linewidth]{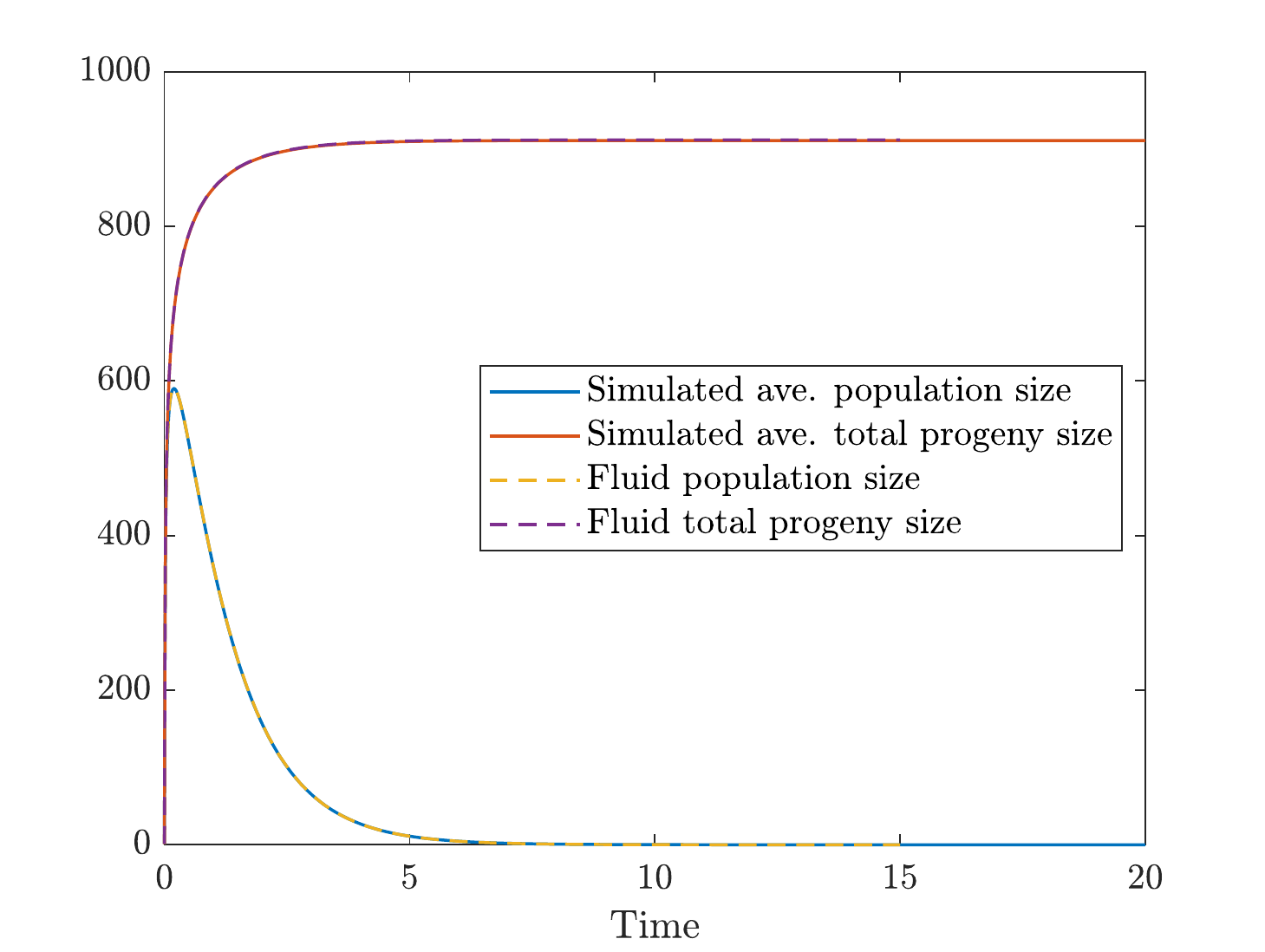}
         \caption{$\lambda=1000$}
         \label{fig: ave simulation 1000 vs ode}
     \end{subfigure}
	\caption{\textbf{Model 2 }--- Average of 10000 trajectories of the stochastic process $\{(Z_t, X_t)\}$ and its fluid approximation $(y_1(t),y_2(t))$.}
	\label{fig: ave simulation vs ode m2}
\end{figure}

Figures \ref{fig: simulation vs ode} and \ref{fig: ave simulation vs ode} show typical single trajectories and averaged trajectories, respectively, of the stochastic process $\{(Z_t, X_t)\}$ and its fluid approximation $(y_1(t),y_2(t))$, for Model 1 with $\mu=1$, and with $\lambda=10$ (left panels) and $\lambda=1000$ (right panels). We see that the population sizes grow to a maximum value before slowly declining to extinction due to a shortage of non-renewable resources. We also note that the fluid curves accurately approximate the mean behaviour of the stochastic process (as well as single trajectories) when $\lambda=1000$.  Figures \ref{fig: simulation vs ode m2} and \ref{fig: ave simulation vs ode m2} show equivalent trajectories for Model 2 with $\mu=1$, $\alpha=100$, and with $\lambda=10$ (left panels) and $\lambda=1000$ (right panels). When comparing the figures, we observe that the peak in the population size is sharper in Model 2, that is, extinction happens faster after the maximum population size has been reached, and we also see that the total progeny curve has an inflection point in Model 2 but not in Model~1.

Motivated by the fact that the averaged trajectories of the stochastic processes  in Figures \ref{fig: ave simulation vs ode} and \ref{fig: ave simulation vs ode m2} are well approximated by their deterministic counterparts when $\lambda$ is sufficiently large, we will use the fluid approximations to study properties of total-progeny-dependent birth-and-death processes such as the mean maximum population size, the mean total progeny at extinction, the mean time at which maximum is reached, and  the mean time until extinction. In addition, the fluid approach will allow us to investigate the behaviour of these quantities as the magnitude of the individual birth rate  increases.

Before pursuing, observe that
by dividing Eq. (\ref{eq1}) by Eq. (\ref{eq2}), we obtain
\begin{eqnarray} \label{div}
\frac{dy_1(t)}{dy_2(t)} &=& \frac{b(y_2(t)) - d(y_2(t))}{b(y_2(t))} = 1 - \frac{d(y_2(t))}{b(y_2(t))},
\end{eqnarray} 
which provides us with an expression for the (fluid) population size at time $t$, $y_1(t)$, as a function of the (fluid) total progeny until time $t$, $y_2(t)$:
\begin{eqnarray}  \label{y_1 in y_2} y_1(t) &=& y_2(t) - \int \frac{d(y_2(t))}{b(y_2(t))} \,dy_2(t).
\end{eqnarray}
The (fluid) maximum population size is obtained by setting $dy_1(t)/dt=0$ in Eq.~(\ref{eq1}), which leads to an equation for the total progeny at the time this maximum is reached, 
\begin{equation} \label{at max population}
b(y_2(t)) = d(y_2(t)).
\end{equation}In other words, the process starts its descent to extinction when the birth rates becomes smaller than the death rate. 

In the next two sections we study the properties of  Models 1 and 2. We provide more details for Model 1, whose fluid approximation is amenable to a complete analytical analysis.

\section{Properties of Model 1}\label{sec_mod1}

Using $b_1(x) = \lambda/x$ and $d(x) = \mu$ in Eqs. \eqref{div} and (\ref{y_1 in y_2}), we obtain
\begin{equation} \label{model 1 y1 in y2 with constant}
\frac{dy_1(t)}{dy_2(t)} = 1 - \frac{\mu}{\lambda}y_2(t) \quad \Longrightarrow \quad y_1(t) = y_2(t) - \frac{\mu}{2\lambda}y_2(t)^2 + C_1,
\end{equation}for some constant $C_1$. Using the initial condition 
\begin{equation} \label{initial condition}
y_1(0) = y_2(0) = 1
\end{equation}in Eq. (\ref{model 1 y1 in y2 with constant}), we obtain $C_1 = \frac{\mu}{2\lambda}$, which leads to the final expression
\begin{equation} \label{model 1 y1 in y2}
 y_1(t) = y_2(t) - \frac{\mu}{2\lambda}y_2(t)^2 + \frac{\mu}{2\lambda}.
\end{equation}
Eq. (\ref{model 1 y1 in y2}) establishes a relationship between the current population size $y_1(t)$ and the total progeny size $y_2(t)$ at any time $t \geq 0$ which will be useful in order to derive the model properties.

\subsection{Maximum population size}
\begin{figure}[H]
\begin{subfigure}[h]{0.5\linewidth}
         \centering
         \includegraphics[width=\linewidth]{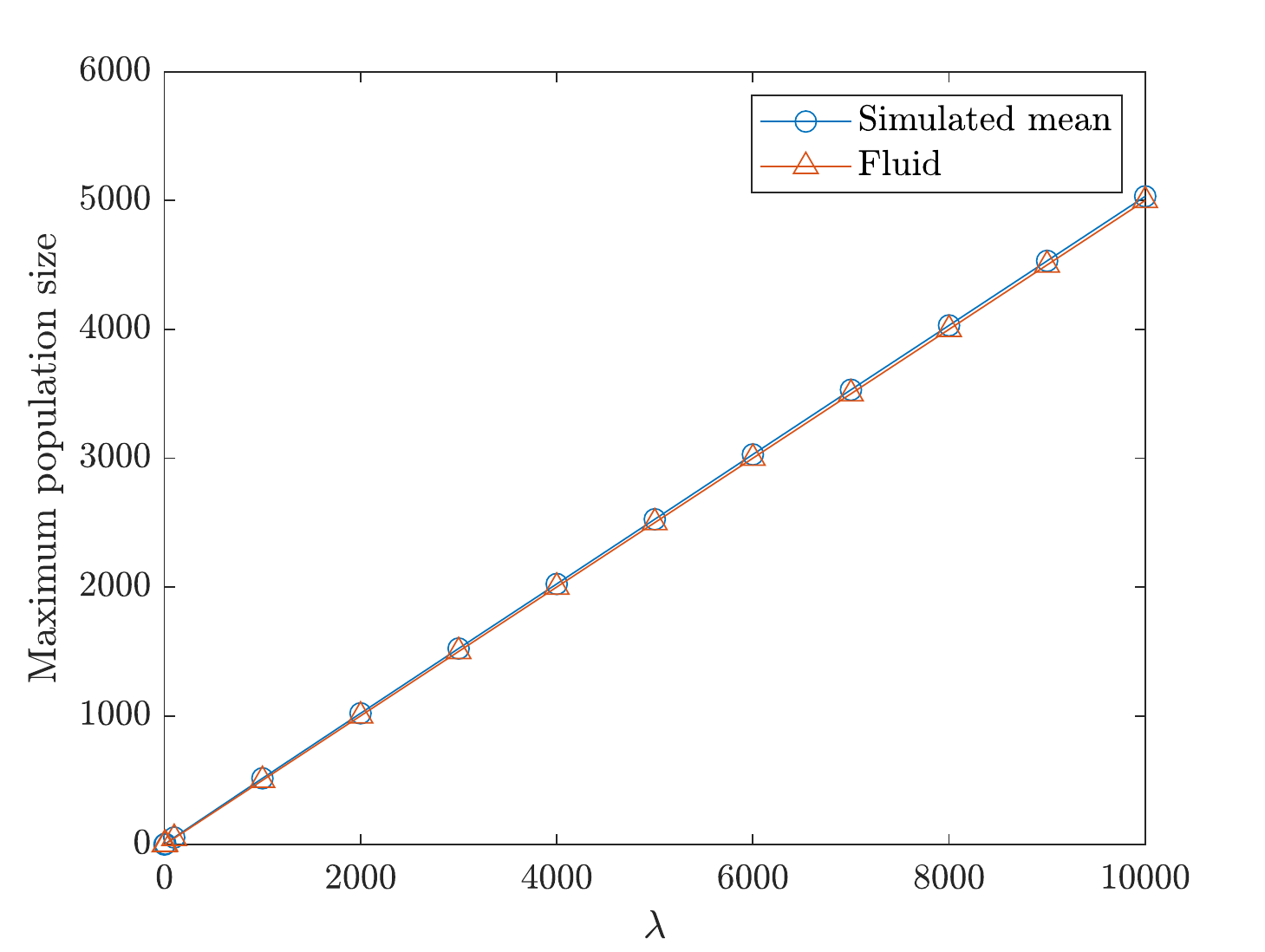}
     \end{subfigure}
     \hfill
     \begin{subfigure}[h]{0.5\linewidth}
         \centering
         \includegraphics[width=\textwidth]{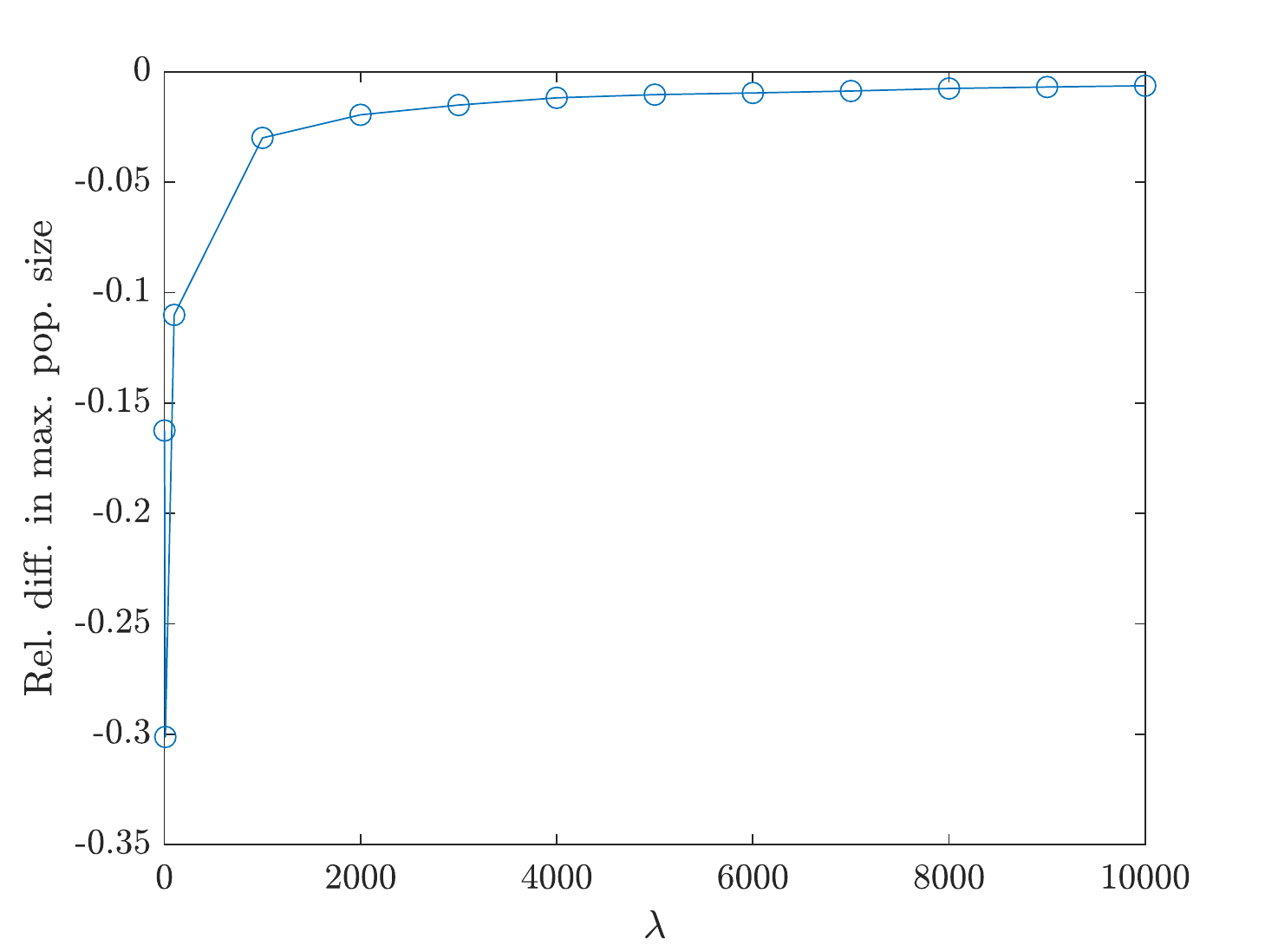}
     \end{subfigure}
	\caption{\textbf{Model 1} --- Left: Mean maximum population size as a function of $\lambda$, when $\mu = 1$. Right: Relative difference between the average of the simulations and the fluid approximation. }
	\label{fig: model 1 max pop}
\end{figure}

Let $t_{max}$ be the time when the population size reaches its maximum value.  Eqs. (\ref{at max population}) and (\ref{model 1 y1 in y2}) provide us with an approximation of the total progeny and population size at time $t_{max}$, respectively:
\begin{align}
y_2(t_{max}) &= \frac{\lambda}{\mu}, \label{Model 1: y2(t_max)}\\
y_1(t_{max}) &= \frac{\lambda}{2\mu} + \frac{\mu}{2\lambda}. \label{model 1 y_1 max pop}
\end{align}
These equations indicate that, for large values of $\lambda$, the population alive at time $t_{max}$ accounts for about half of the total progeny until that time.
Figure \ref{fig: model 1 max pop} (left) compares the value of $y_1(t_{max})$ given in Eq. \eqref{model 1 y_1 max pop} and the average maximum population size based on 5000 simulations of the stochastic process for increasing values of $\lambda$. The relative difference is highlighted in the right panel of the figure and indicates that the fluid approach provides an accurate approximation of the mean maximum population size.

\subsection{Total progeny at extinction}
Extinction in the stochastic process corresponds to the event $\{(Z(t),X(t))=(0,x)$ for some $t>0$ and $x>0\}$. By definition, a total-progeny-dependent birth-and-death process becomes subcritical as soon as the total progeny becomes larger than $y_2(t_{max})$; from that point, the process  $Z(t)$ is stochastically dominated by a subcritical linear birth-and-death process, which implies  that extinction occurs almost surely. In the fluid approximation, if $y_1(0)>0$, there is no time $t$ such that $y_1(t)=0$, but instead we have $y_1(\infty):=\lim_{t \to \infty} y_1(t) = 0$. We denote by $y_2(\infty):=\lim_{t \to \infty} y_2(t)$ the  total progeny at extinction.
Letting $t\to\infty$ in Eq. (\ref{model 1 y1 in y2}) we obtain
\begin{equation}\label{model 1: y2 ext}
y_2(\infty) = \frac{\lambda + \sqrt{\lambda^2 + \mu^2}}{\mu}.
\end{equation}

\begin{figure}[H]
\begin{subfigure}[h]{0.5\linewidth}
         \centering
         \includegraphics[width=\linewidth]{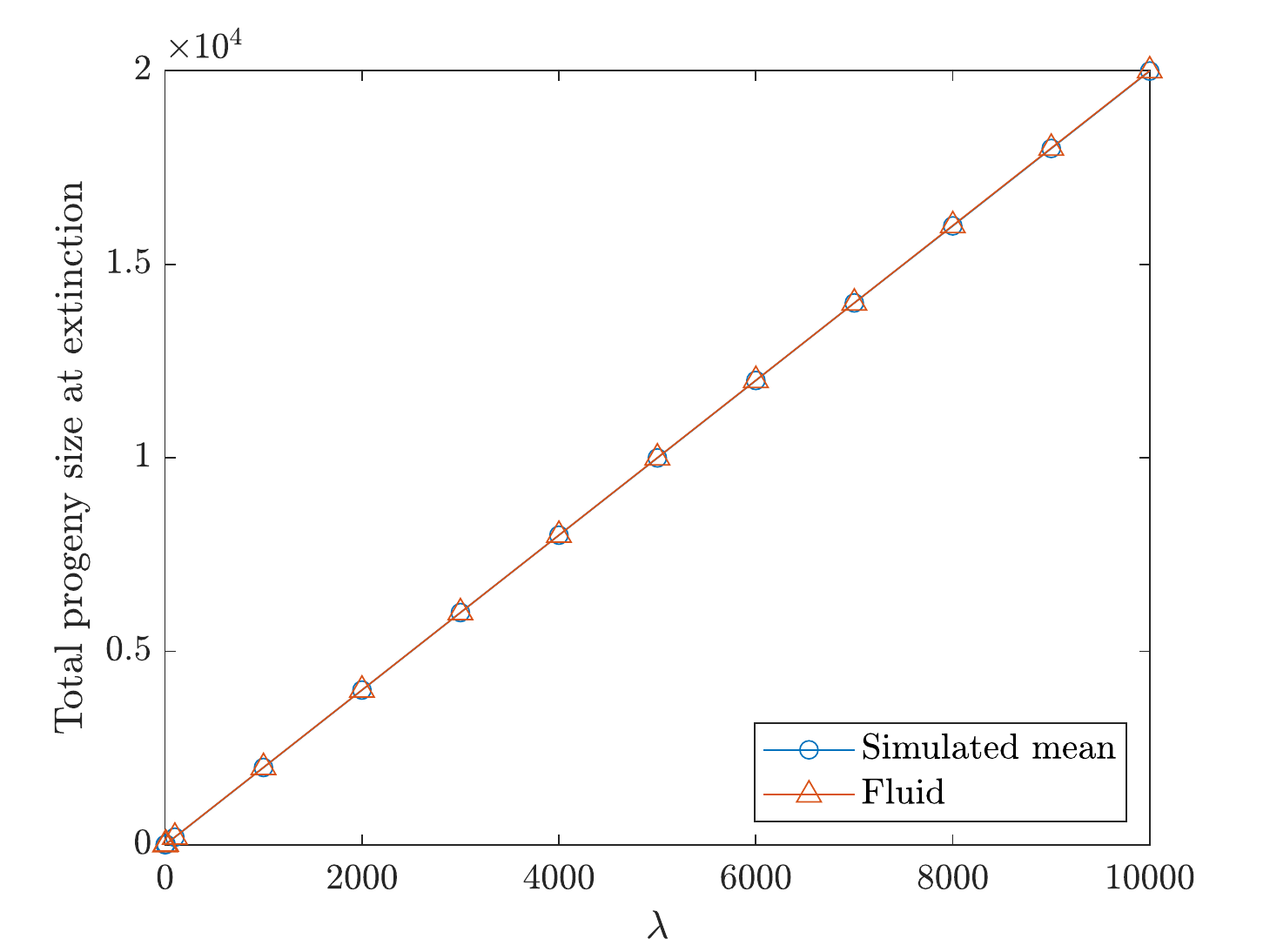}
     \end{subfigure}
     \hfill
     \begin{subfigure}[h]{0.5\linewidth}
         \centering
         \includegraphics[width=\textwidth]{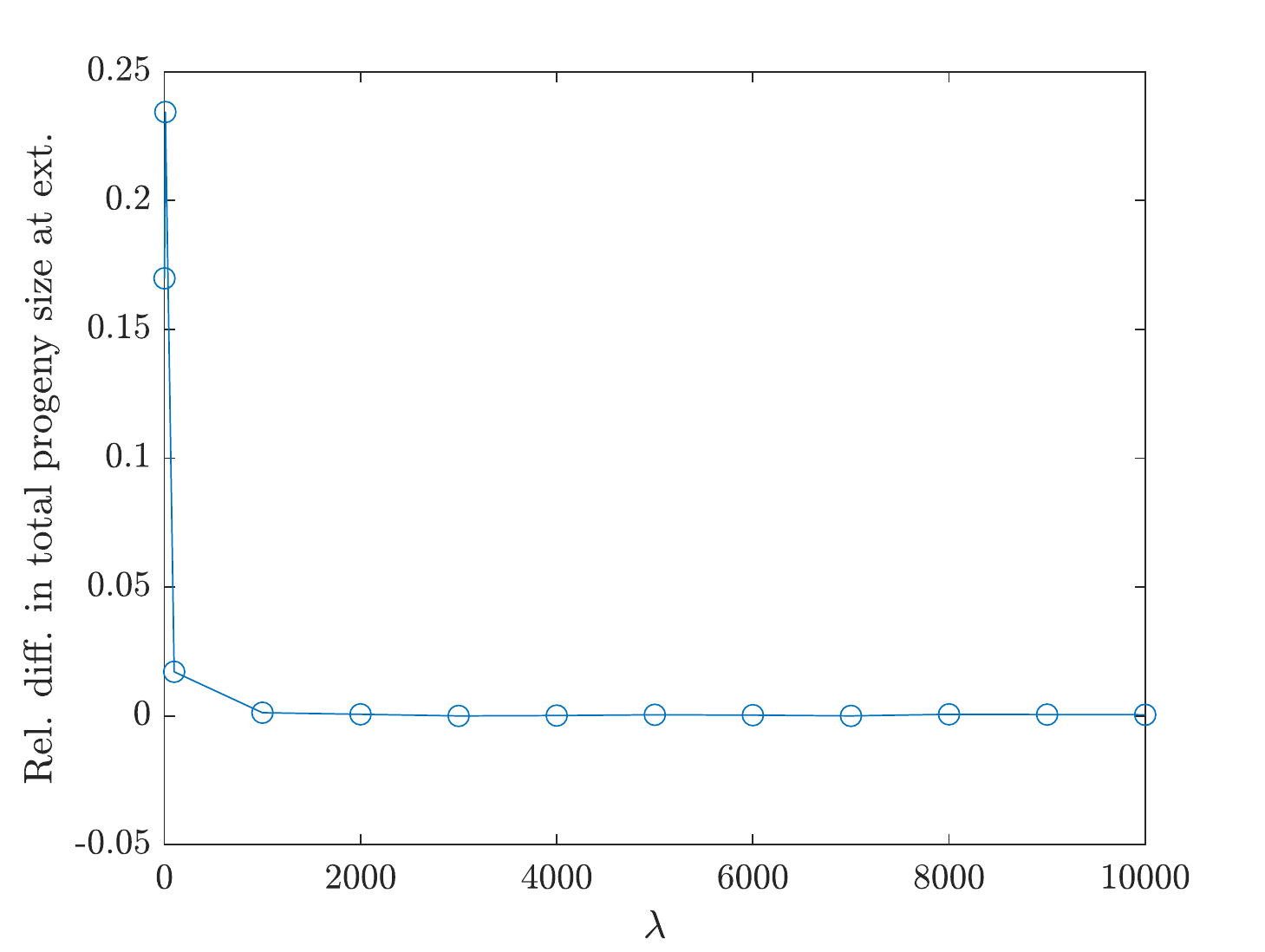}
     \end{subfigure}
	\caption{\textbf{Model 1} --- Left: Mean total progeny at extinction as a function of $\lambda$, with $\mu = 1$. Right: Relative difference between the average of the simulations and the fluid approximation. }
	\label{fig: model 1 total prog at ext}
\end{figure}

Figure \ref{fig: model 1 total prog at ext} indicates that, as $\lambda$ increases, the fluid approach provides an accurate approximation to the mean total progeny size.
However, this approach does not directly provide us with the mean time until extinction  (which, in the fluid model, is infinite due to the nature of the solution \eqref{sol_eq1}). Methods to approximate the mean time until extinction are presented in Section \ref{sec:text}.

\subsection{Time to reach a given total progeny size}

\begin{figure}[h!]
\centering
\includegraphics[width=0.9\linewidth]{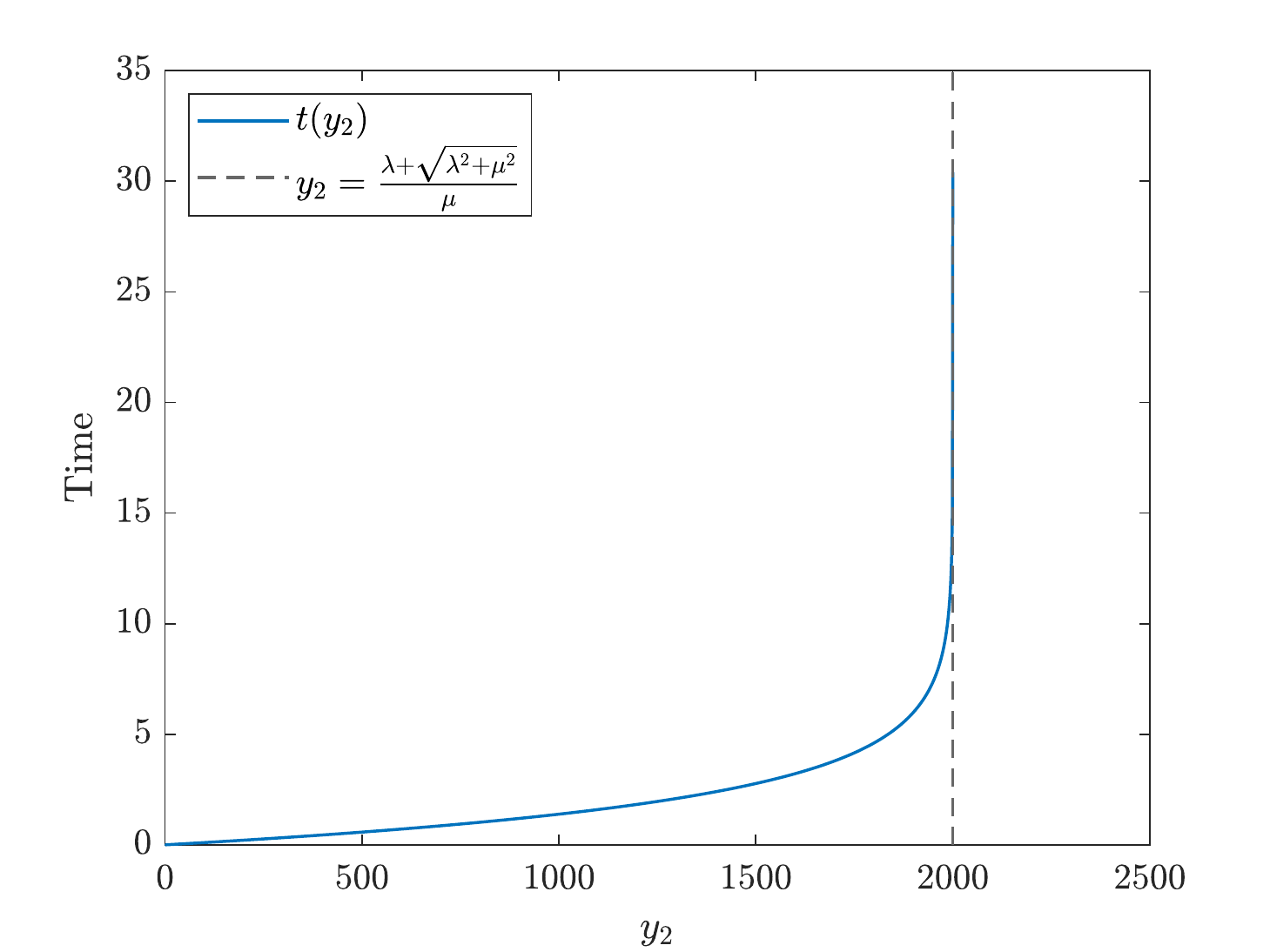}
\caption{\textbf{Model 1} ---  Time at which the total progeny reaches a given value $y_2\leq y_2(\infty)$ given by Eq. (\ref{model 1 t with y2}), with $\lambda=1000$.}
\label{fig: model 1 t vs y2}
\end{figure}

Since $y_2(t)$ increases monotonically, its inverse $t(y_2)$,  which records the time at which the total progeny reaches a given value $y_2\leq y_2(\infty)$, exists and increases monotonically. We now derive an expression for $t(y_2)$.
From Eqs.~(\ref{eq2}) and (\ref{model 1 y1 in y2}), we have
\begin{equation}
\frac{dy_2}{dt} = y_1b(y_2) = \left(y_2 - \frac{\mu}{2\lambda}y_2^2 + \frac{\mu}{2\lambda}\right)\frac{\lambda}{y_2} = \lambda - \frac{\mu}{2}y_2 + \frac{\mu}{2y_2},
\end{equation}
which gives
\begin{equation}
\label{model 1 ode}
 \frac{dt}{dy_2} = \frac{y_2}{\lambda y_2 - \frac{\mu}{2}y_2^2 + \frac{\mu}{2}} = -\frac{1}{\mu}\left(\frac{1 - \frac{\lambda}{\mu c}}{y_2 - \frac{\lambda}{\mu} + c} + \frac{1 + \frac{\lambda}{\mu c}}{y_2 - \frac{\lambda}{\mu} - c}\right),
\end{equation}
where \begin{equation}\label{c}c^2 = 1 + \frac{\lambda^2}{\mu^2}.\end{equation}
By integrating Eq. (\ref{model 1 ode}) with the initial condition $y_2(0) = 1$, we obtain 
\begin{equation} \label{model 1 t with y2}
t(y_2) = \frac{1}{\mu}\left[ (1 - \frac{\lambda}{\mu c}) \log\left(\dfrac{1 - \frac{\lambda}{\mu} + c}{y_2 - \frac{\lambda}{\mu} + c} \right) + (1 + \frac{\lambda}{\mu c}) \log\left(\dfrac{\frac{\lambda}{\mu} + c - 1}{\frac{\lambda}{\mu} + c - y_2}\right)\right]. 
\end{equation}Note that the constant $c$ in \eqref{c} has two possible values, however the solution $t(y_2)$ is independent of the choice of value. 
Figure \ref{fig: model 1 t vs y2} depicts the function $t(y_2)$ and highlights its asymptotic behaviour as $y_2$ approaches the total progeny at extinction, $y_2(\infty)$ (represented by the vertical dashed line).

\subsection{Time to  reach the maximum population size}

By substituting Eq. \eqref{Model 1: y2(t_max)} into Eq. \eqref{model 1 t with y2} we obtain $t_{max}:=t(y_2(t_{max}))$, the (fluid) time to  reach the maximum population size:
\begin{equation} \label{t_max expression}
t_{max} = \frac{1}{\mu}\left[ (1 - \frac{\lambda}{\mu c}) \log\left(1 - \frac{\lambda}{\mu} + c\right) + (1 + \frac{\lambda}{\mu c}) \log\left(\frac{\lambda}{\mu} + c - 1\right)-2\log c\right],
\end{equation}where $c$ is given in Eq. \eqref{c}.

\begin{figure}[H]
     \begin{subfigure}[h]{0.5\linewidth}
         \centering
         \includegraphics[width=\linewidth]{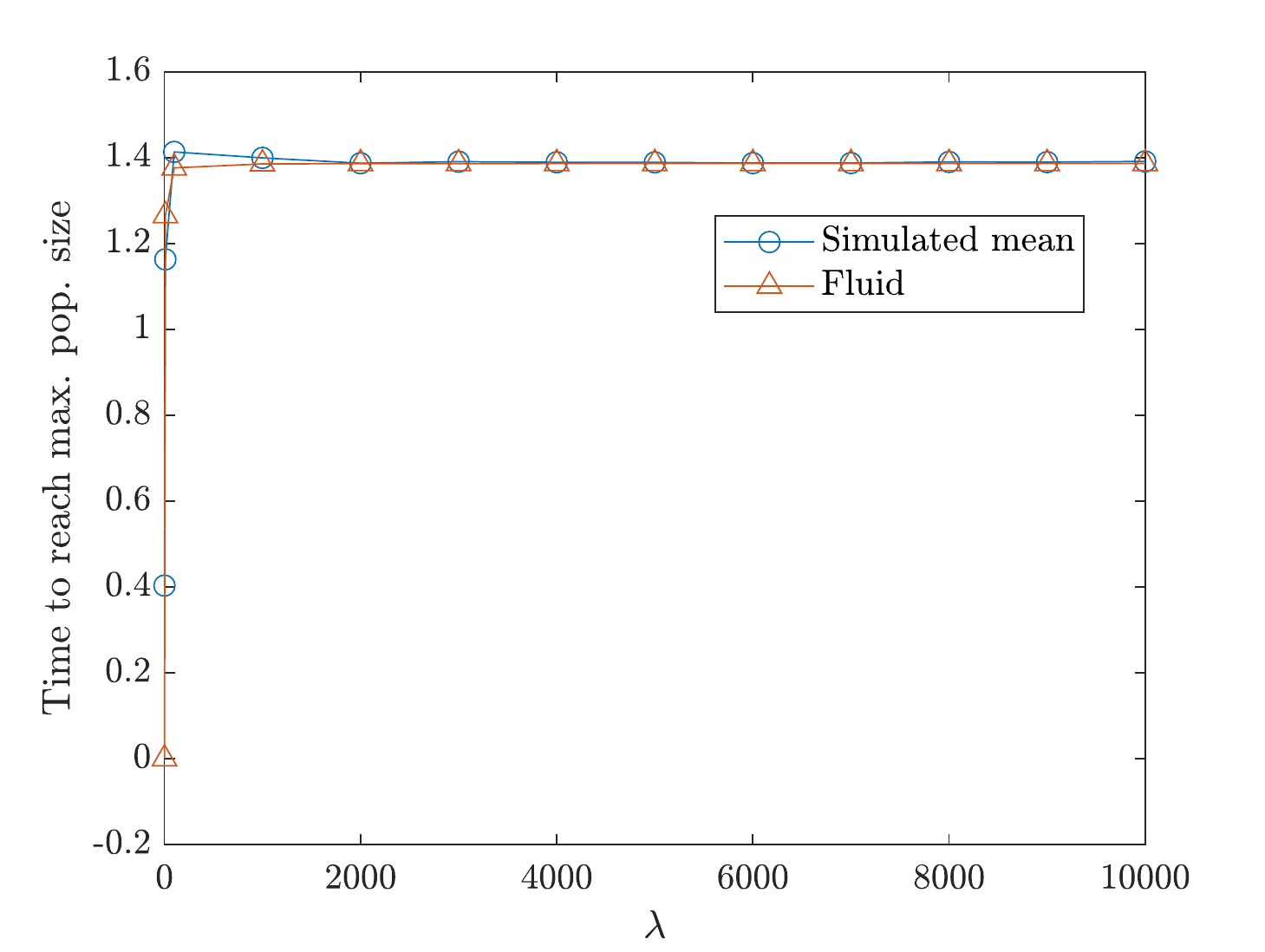}
     \end{subfigure}
     \hfill
     \begin{subfigure}[h]{0.5\linewidth}
         \centering
         \includegraphics[width=\linewidth]{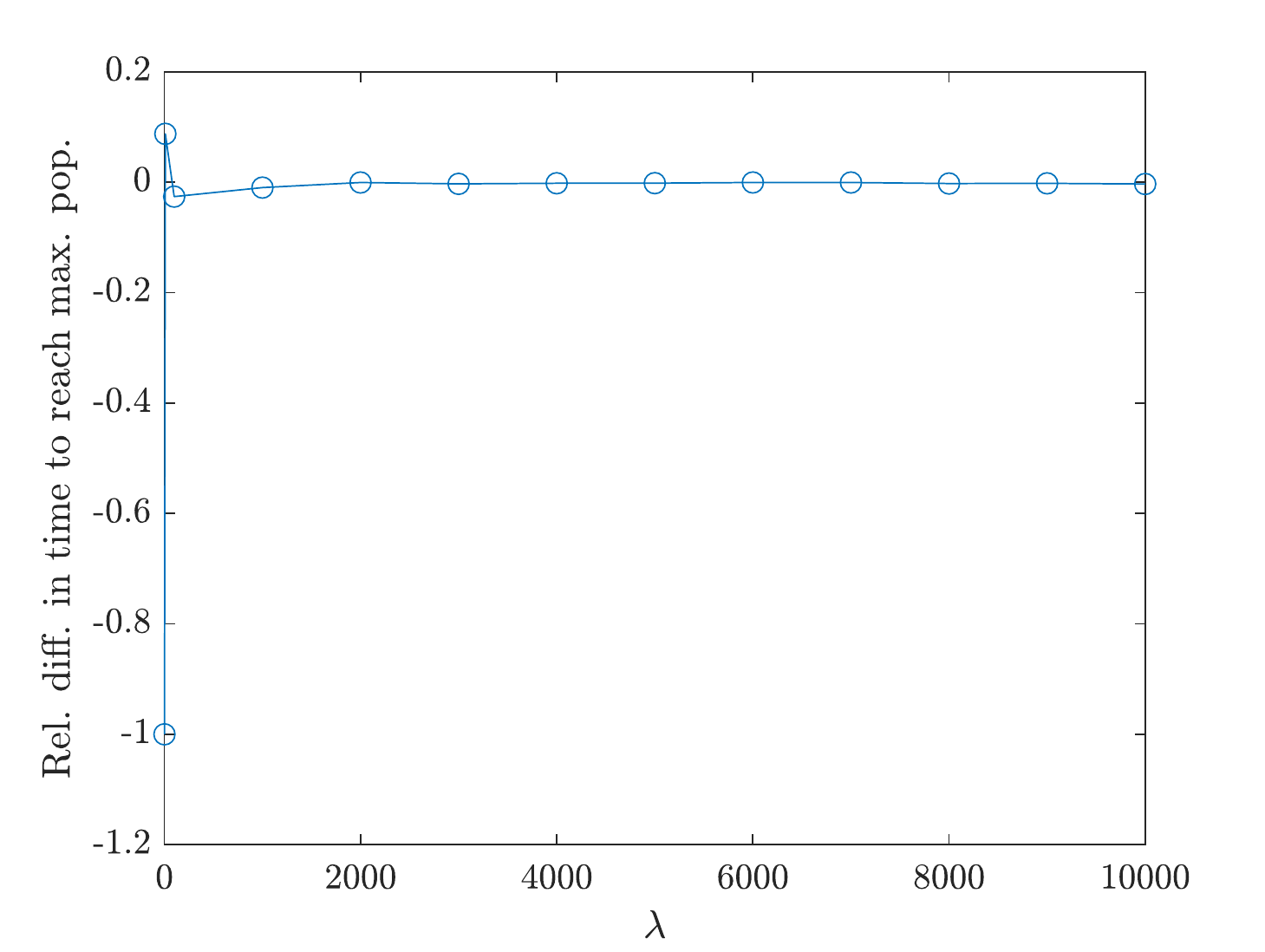}
     \end{subfigure}
        \caption{\textbf{Model 1} --- Left: Mean time to reach the maximal population size as a function of $\lambda$, with $\mu = 1$. Right: Relative difference between the average of the simulations and the fluid approximation.}
        \label{fig:model 1 time at max pop}
\end{figure}

Figure \ref{fig:model 1 time at max pop} indicates that for $\lambda$ suitably large, the fluid approach provides an accurate approximation to the mean time to reach the maximum population size. The left panel of the figure also highlights an interesting  property: as $\lambda$ increases, the mean time to reach the maximum population size converges to a constant. We show this result for the fluid function $t_{max} =t_{max} (\lambda)$ (the proof can be found in Section \ref{proofs}):

\begin{prop}\label{prop1}
$t_{max} \to \frac{2}{\mu}\log(2)$ as $\lambda\to\infty$.
\end{prop}
This means that, as $\lambda$ increases, the birth rate $b(y_2(t))$ increases in a way which perfectly balances the asymptotic linear increase in the maximum population size $y_1(t_{max})$.

\subsection{Time until extinction}\label{sec:text}

As illustrated in Figure \ref{fig: model 1 t vs y2},  the fluid approach does not allow us to compute the time until extinction directly since $t(y_2(\infty))=\infty$. We therefore need additional steps to approximate this time. 

Let $\varepsilon\geq 1$ denote a small population size, and let $t_\varepsilon$ denote the time at which the (fluid) population reaches size $\varepsilon$ in its descent, that is, such that $y_1(t_\varepsilon)=\varepsilon$, $t_\varepsilon>t_{max}$. The first approach is to approximate the mean time until extinction, $t_{ext}$, by 
\begin{equation}\label{text_approx1}t_{ext}^{(\varepsilon)}:= t_\varepsilon+\mathbb E[\max(X_1,\ldots,X_\varepsilon)],\end{equation}
where the random variables $X_i$ are i.i.d. exponential random variables with rate $d(y_2(t_\varepsilon))$, $1\leq i\leq \varepsilon$. In other words, we assume that after time $t_\varepsilon$, no more birth events occur and the process behaves like a pure-death process until extinction. This approximation is reasonable for values of $\lambda$ such that $y_2(t_\varepsilon)$ is large enough  for the birth rate $b(y_2(t_\varepsilon))$ to be negligible compared to the death rate $d(y_2(t_\varepsilon))$. Note that the expectation in \eqref{text_approx1} has a closed-form expression:
$$\mathbb E[\max(X_1,\ldots,X_\varepsilon)]=\dfrac{1}{d(y_2(t_\varepsilon))}\sum_{j=1}^\varepsilon \binom{\varepsilon}{j} \dfrac{(-1)^{j-1}}{j}.$$
To find $t_\varepsilon$, we let $y_1(t_\varepsilon)=\varepsilon$ in Eq. (\ref{model 1 y1 in y2}), which leads to an equation for $y_2(t_{\varepsilon})$,
\begin{equation}
\varepsilon = y_2(t_\varepsilon) - \frac{\mu}{2\lambda}y_2(t_\varepsilon)^2 + \frac{\mu}{2\lambda},
\end{equation}
whose larger solution is
\begin{equation}\label{extinction y2}
y_2(t_{\varepsilon}) = \frac{\lambda + \sqrt{\lambda^2 + \mu^2 - 2\varepsilon \mu \lambda}}{\mu}.
\end{equation}
Substituting Eq. (\ref{extinction y2}) into Eq. (\ref{model 1 t with y2}), we obtain
\begin{eqnarray}\nonumber
t_{\varepsilon} &=& \frac{1}{\mu}\left[ (1 - \frac{\lambda}{\mu c}) \log\left(\dfrac{1 - \frac{\lambda}{\mu} + c}{\frac{ \sqrt{\lambda^2 + \mu^2 - 2\varepsilon \mu \lambda}}{\mu}  + c} \right) \right.\\ \label{t in lambda and eplison}&&\left.\qquad+ (1 + \frac{\lambda}{\mu c}) \log\left(\dfrac{\frac{\lambda}{\mu} + c - 1}{ c - \frac{ \sqrt{\lambda^2 + \mu^2 - 2\varepsilon \mu \lambda}}{\mu}}\right)\right],
\end{eqnarray}where $c$ is given in Eq. \eqref{c}. 

Figure \ref{fig: model 1 t vs lambda} shows $t_{ext}^{(1)}=t_1+1/\mu$ (that is, $t_{ext}^{(\varepsilon)}$ with $\varepsilon=1$),
and the average time until extinction based on simulations. We see that there is a distinct gap between the two curves, which indicates that the choice of $\varepsilon=1$ may not be suitable. Note that when $\varepsilon=1$, Eq. \eqref{extinction y2} reduces to $y_2(t_{1})=2(\lambda/\mu) -1$.

\begin{figure}[H]
     \begin{subfigure}[h]{0.5\linewidth}
         \centering
         \includegraphics[width=\linewidth]{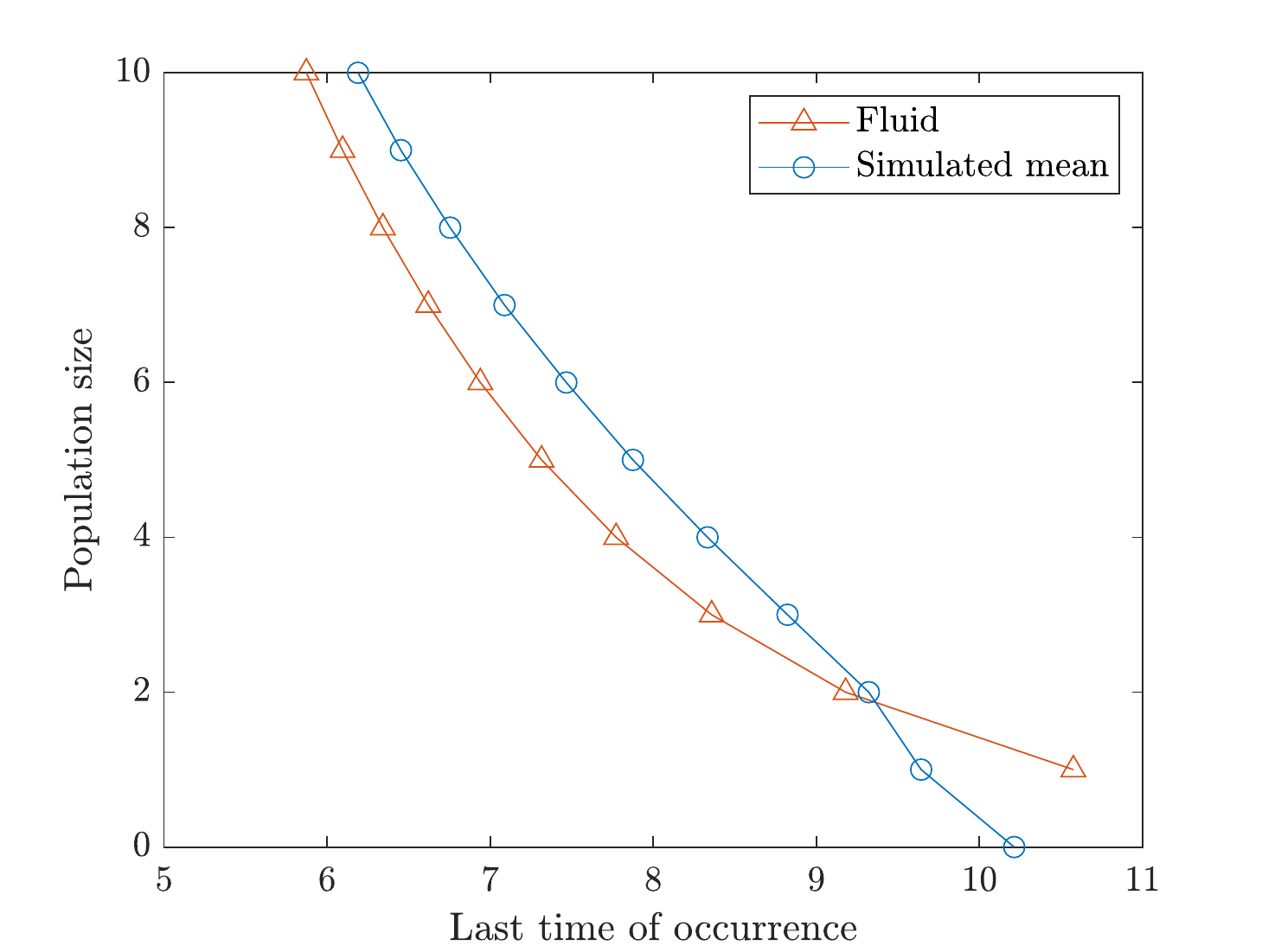}
         \caption{$\lambda = 100$}
         \label{fig:model 1 small pop size vs t 100}
     \end{subfigure}
     \hfill
     \begin{subfigure}[h]{0.5\linewidth}
         \centering
         \includegraphics[width=\linewidth]{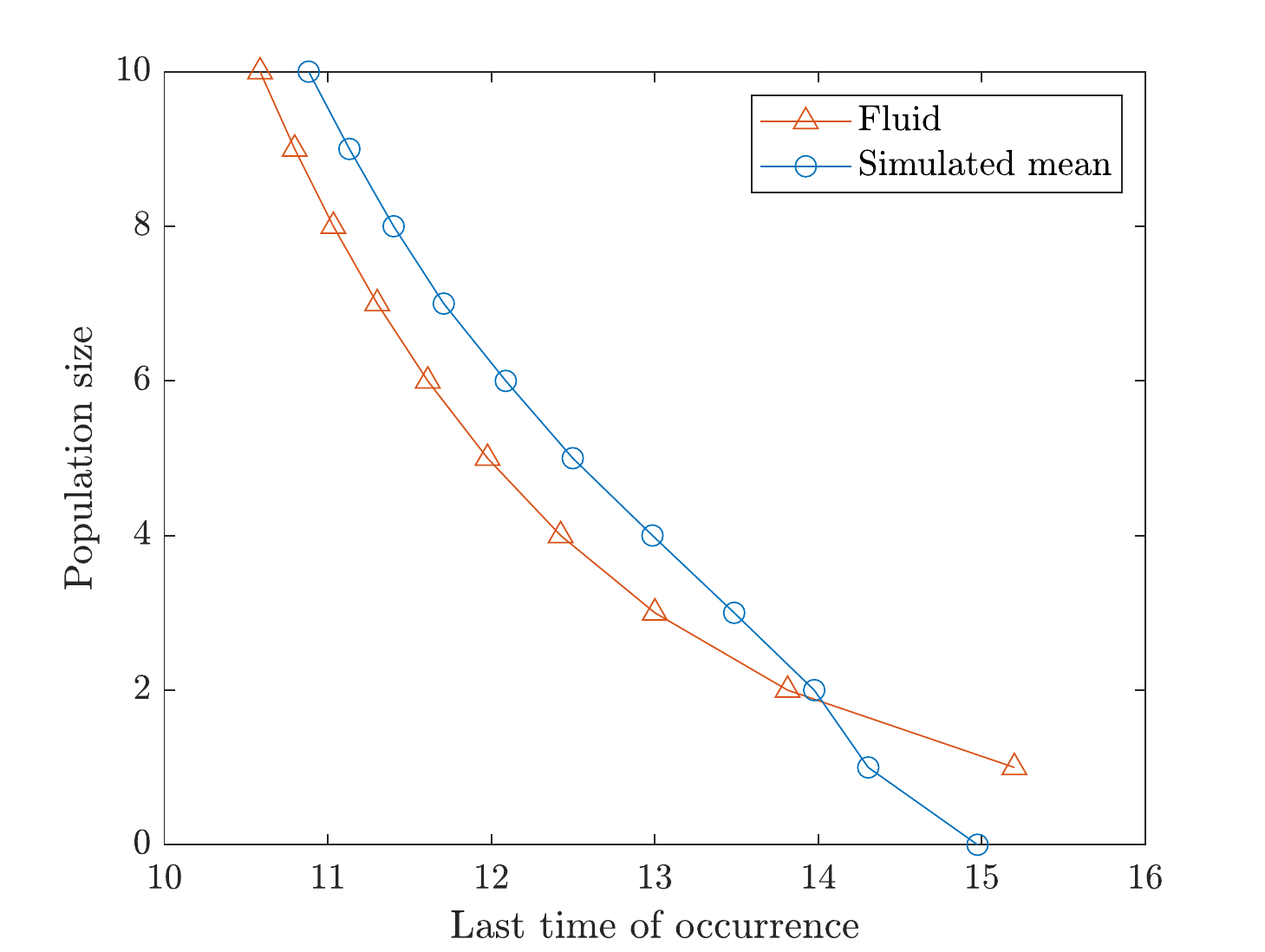}
         \caption{$\lambda = 1000$}
         \label{fig:model 1 small pop size vs t 1000}
     \end{subfigure}
        \caption{\textbf{Model 1} --- Small population sizes and their last occurrence time.}
        \label{fig: model 1 small pop size vs t}
\end{figure}

To further investigate the differences between the stochastic process and its fluid approximation, in Figure \ref{fig: model 1 small pop size vs t} we plot the average time at which the population reaches some small sizes $\varepsilon$ for the last time before extinction in the stochastic process against $t_{\varepsilon}$,  for  $\lambda = 100$ and $\lambda=1000$. The graphs indicate that choosing $\varepsilon=2$ may be more appropriate than $\varepsilon=1$; this is confirmed in Figure \ref{fig: model 1 t vs lambda}, which we comment after introducing our second method.

\begin{figure}[H]
\begin{subfigure}[h]{0.5\linewidth}
         \centering
         \includegraphics[width=\linewidth]{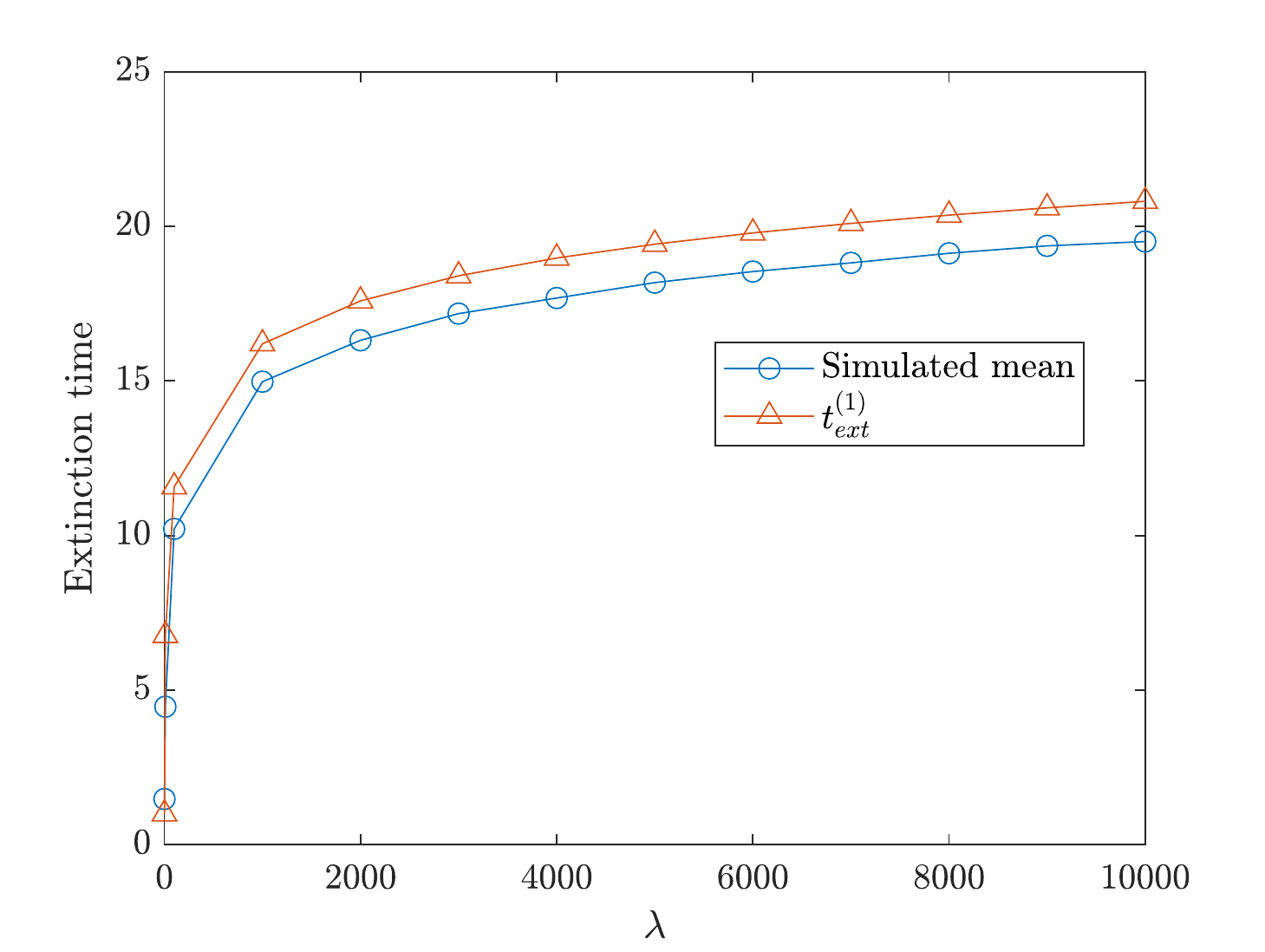}
     \end{subfigure}
     \hfill
     \begin{subfigure}[h]{0.5\linewidth}
         \centering
         \includegraphics[width=\textwidth]{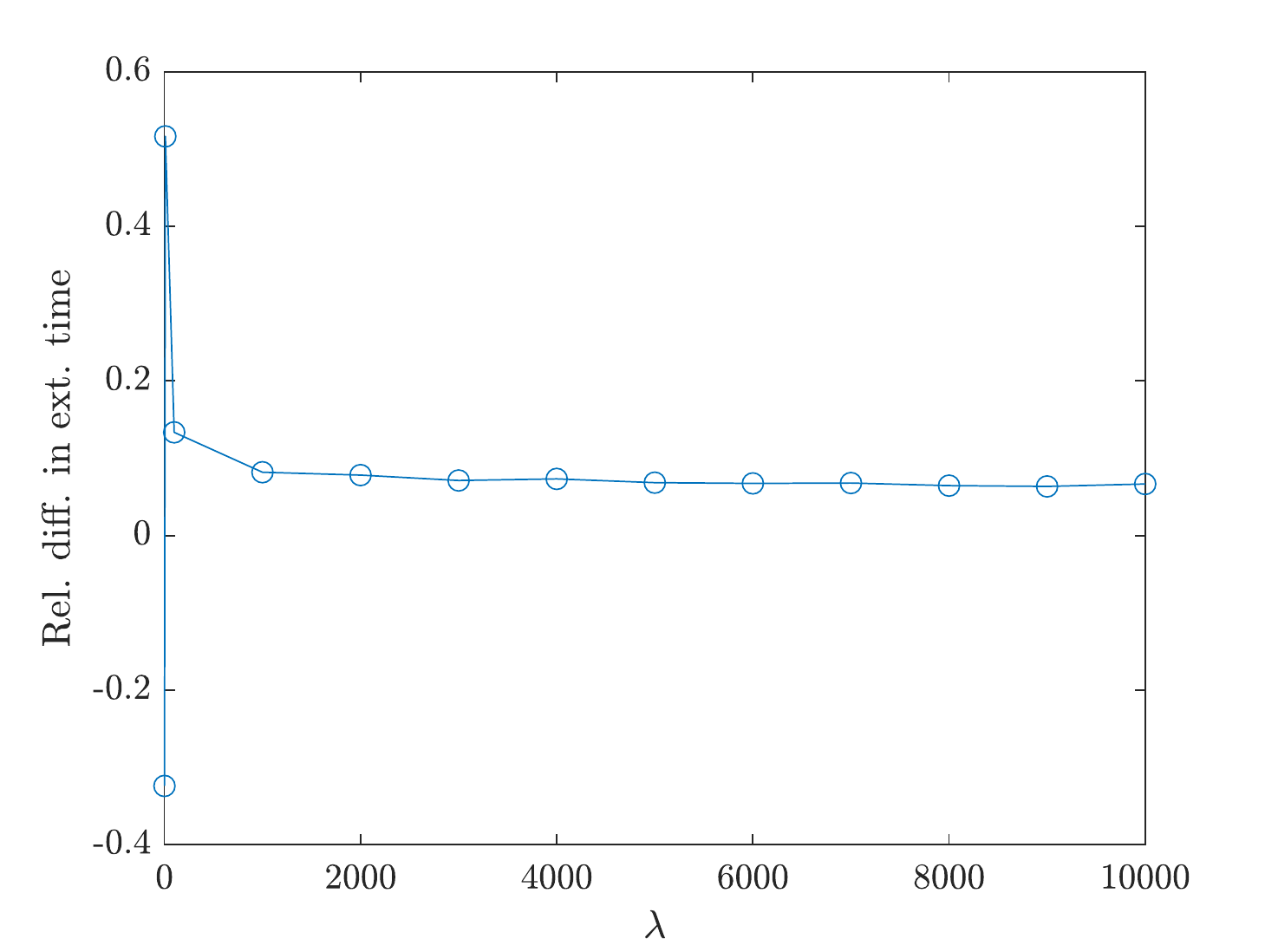}
     \end{subfigure}
	\caption{\textbf{Model 1} ---  Left: Approximated mean time until extinction $t_{ext}^{(\varepsilon)}$ as a function of $\lambda$, with $\mu = 1$ and $\varepsilon=1$. Right: Relative difference between the average of the simulations and the fluid approximation.
	}
	\label{fig: model 1 t vs lambda}
\end{figure}

A second approach to approximate the mean time until extinction is obtained as follows: let $z:=y_2(\infty)-1$, then $t(z)$ (obtained by replacing $y_2$ by $z$ in Eq. \eqref{model 1 t with y2}) approximates the time of the last birth (arguably it could also approximate the time of the second last birth). It follows that  $y_1(t(z))$ (obtained by replacing $y_2$ by $z$ in Eq. \eqref{model 1 y1 in y2}) is the population size at that time, and we assume that the process behaves like a pure death process from then on, so $t_{ext}$ can be approximated by
\begin{equation}\label{text_approx2}t_{ext}^{\star}:= t(z)+\mathbb E[\max(X_1,\ldots,X_{y_1(t(z))})],\end{equation}
where the random variables $X_i$ are i.i.d. exponential random variables with rate $d(z)$, $1\leq i\leq y_1(t(z)).$

\begin{figure}[H]
\centering
\includegraphics[width=0.8\linewidth]{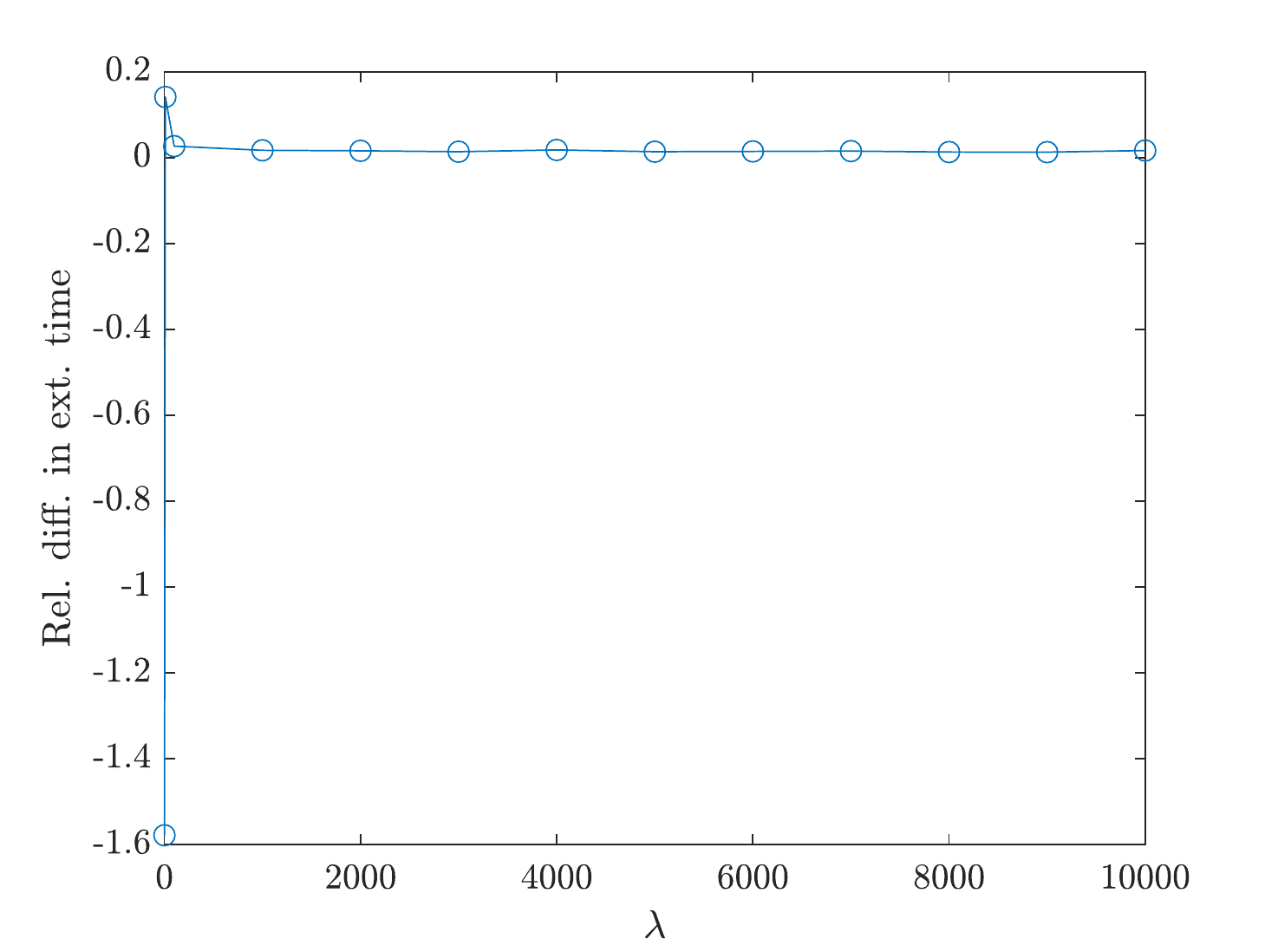}
\caption{\label{fig: model 1 t vs lambda_rel}\textbf{Model 1} --- Relative error between the average time to extinction based on simulations and the approximation  $t_{ext}^{(2)}$ as a function of $\lambda$.
}
\end{figure}

\begin{figure}[H]
\centering
\includegraphics[width=0.8\linewidth]{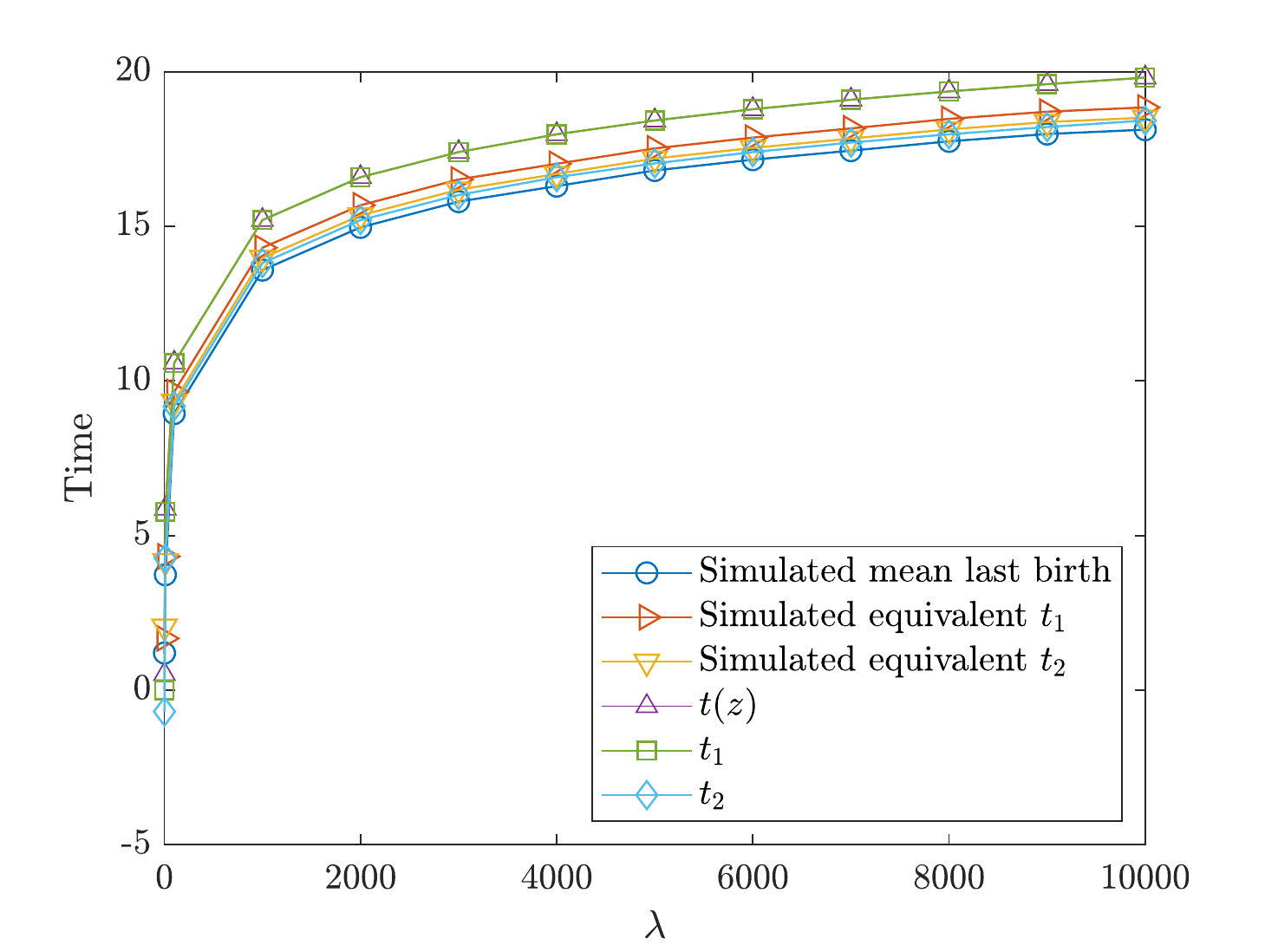}
\caption{\label{fig: model1_time_6_vs_lambda}\textbf{Model 1} --- Time of the last birth, and times to reach population sizes $\varepsilon=1,2$.
}
\end{figure}

Figure \ref{fig: model 1 t vs lambda} compares the methods of approximation $t_{ext}^{(\varepsilon)}$ for $\varepsilon=1,2,$ and $t_{ext}^{\star}$,  with the average time until  extinction based on simulations. We see that for Model 1, the most accurate method is $t_{ext}^{(2)}$. In this figure, we also observe that $t_{ext}^{(1)}$  coincides with $t_{ext}^{\star}$ because  $t_{\varepsilon=1}\approx t(z)$ due to the fact that $y_2(t_{1})=2(\lambda/\mu) -1\approx y_2(\infty)-1$. The relative error of $t_{ext}^{(2)}$ is shown in Figure \ref{fig: model 1 t vs lambda_rel}.
In Figure \ref{fig: model1_time_6_vs_lambda} we plot together $t(z)$, the simulated time of last birth, $t_{\varepsilon}$ for $\varepsilon=1,2$ and their simulated equivalent, as functions of $\lambda$. This figure shows that $t(z)$ over-estimates the  expected time until the last birth, which indicates that our assumption that the process behaves like a pure death process after time $t(z)$ is not appropriate here.

\begin{figure}[H]
\centering
\includegraphics[width=0.8\linewidth]{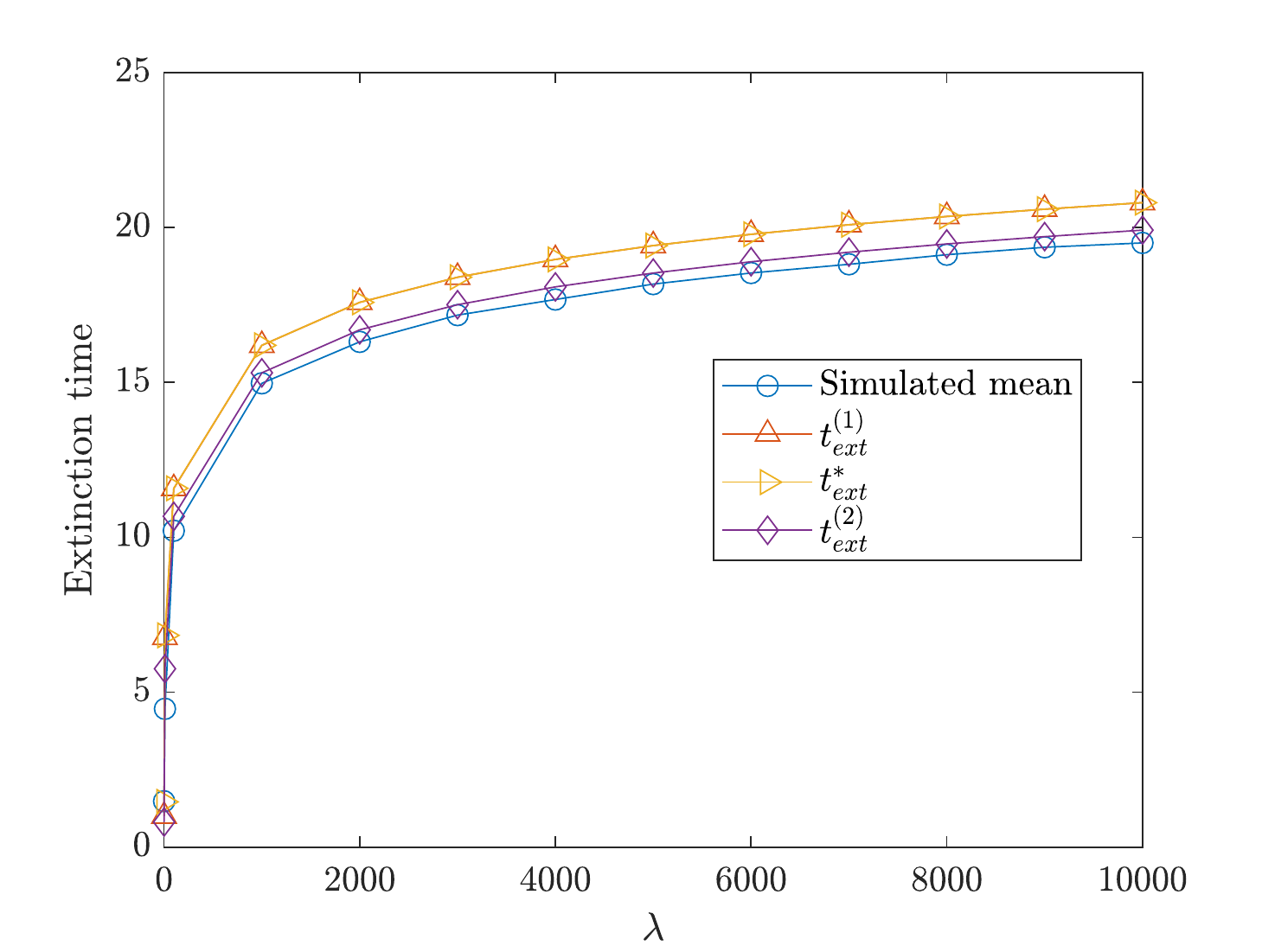}
\caption{\label{fig: model 1 t vs lambda}\textbf{Model 1} --- Mean time until extinction obtained with different approximation methods, as a function of $\lambda$.
}
\end{figure}

\section{Properties of Model 2}\label{sec_mod2}

Recall that in Model 2, $b_2(y_2) = \lambda e^{-\frac{y_2}{\alpha}}$ (both  $\lambda$ and $\alpha$ control the amplitude of the individual birth rate), and $d(y_2) = \mu$. 
Similar to Model 1, from Eqs.~\eqref{div} and (\ref{y_1 in y_2}), we have
\begin{equation} \label{model 2 y1 in y2 with constant}
\frac{dy_1(t)}{dy_2(t)} = 1 - \frac{\mu}{\lambda} e^{\frac{y_2(t)}{\alpha}} \quad \Longrightarrow \quad y_1(t) = y_2(t) - \frac{\alpha\mu}{\lambda}e^{\frac{y_2(t)}{\alpha}} + C_2,
\end{equation}for some constant $C_2$. Using the initial condition \eqref{initial condition}, we obtain $C_2 = \frac{\alpha \mu}{\lambda}e^{\alpha^{-1}}$, hence
\begin{equation} \label{model 2 y1 in y2}
\quad y_1(t) = y_2(t) - \frac{\alpha \mu}{\lambda}e^{\frac{y_2(t)}{\alpha}} +  \frac{\alpha\mu}{\lambda}e^{\alpha^{-1}},\qquad t\geq 0.
\end{equation}

\subsection{Maximum population size}
\begin{figure}[H]
     \begin{subfigure}[h]{0.5\linewidth}
         \centering
         \includegraphics[width=\linewidth]{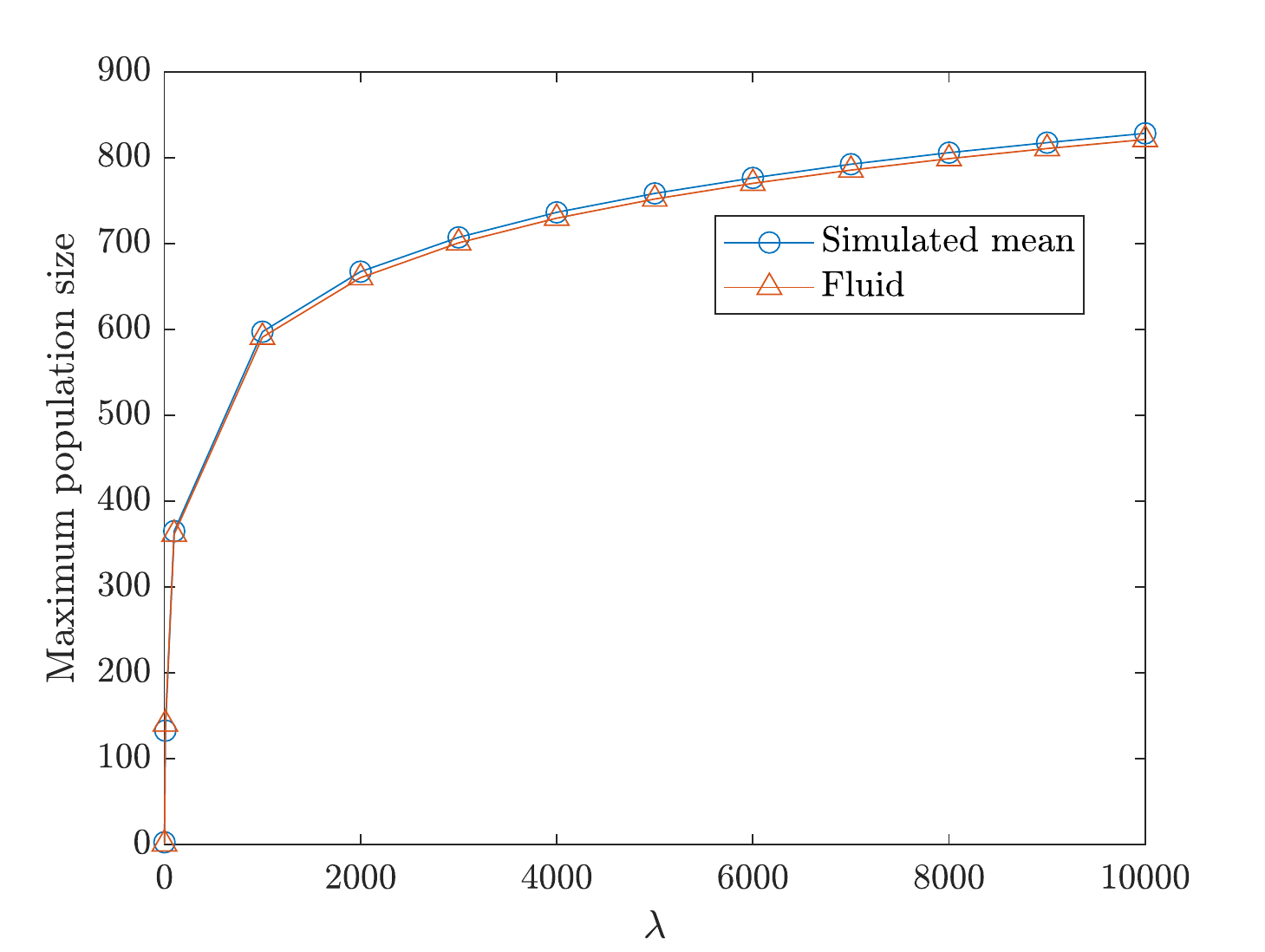}
     \end{subfigure}
     \hfill
     \begin{subfigure}[h]{0.5\linewidth}
         \centering
         \includegraphics[width=\linewidth]{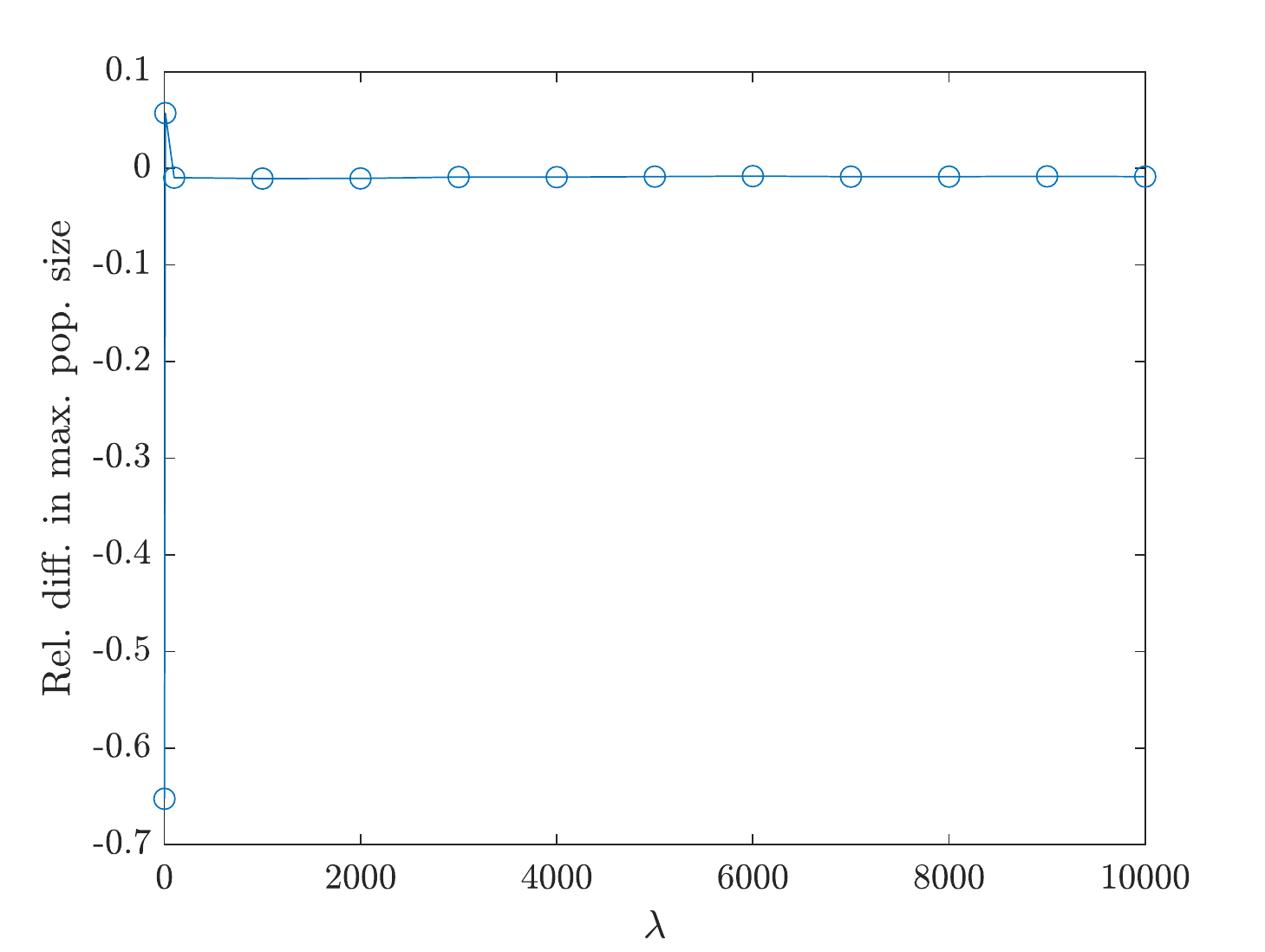}
     \end{subfigure}
        \caption{
        \textbf{Model 2} --- Left: Mean maximum population size as a function of $\lambda$, when $\alpha = 100,\mu = 1$. Right: Relative difference between the average of the simulations and the fluid approximation. }
        \label{fig: model 2: max pop vs lambda}
\end{figure}

\begin{figure}[H]
     \begin{subfigure}[h]{0.5\linewidth}
         \centering
         \includegraphics[width=\linewidth]{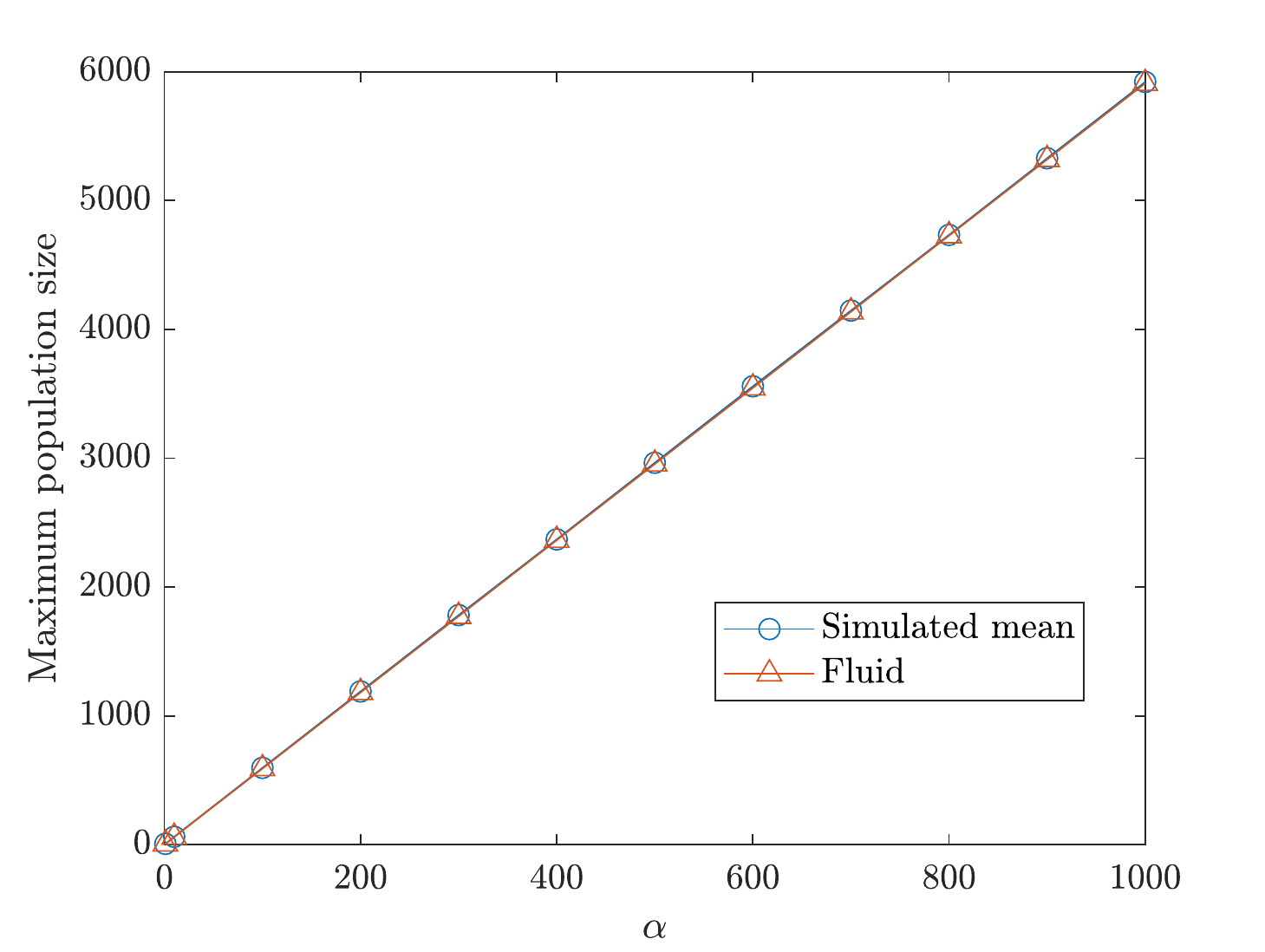}
     \end{subfigure}
     \hfill
     \begin{subfigure}[h]{0.5\linewidth}
         \centering
         \includegraphics[width=\linewidth]{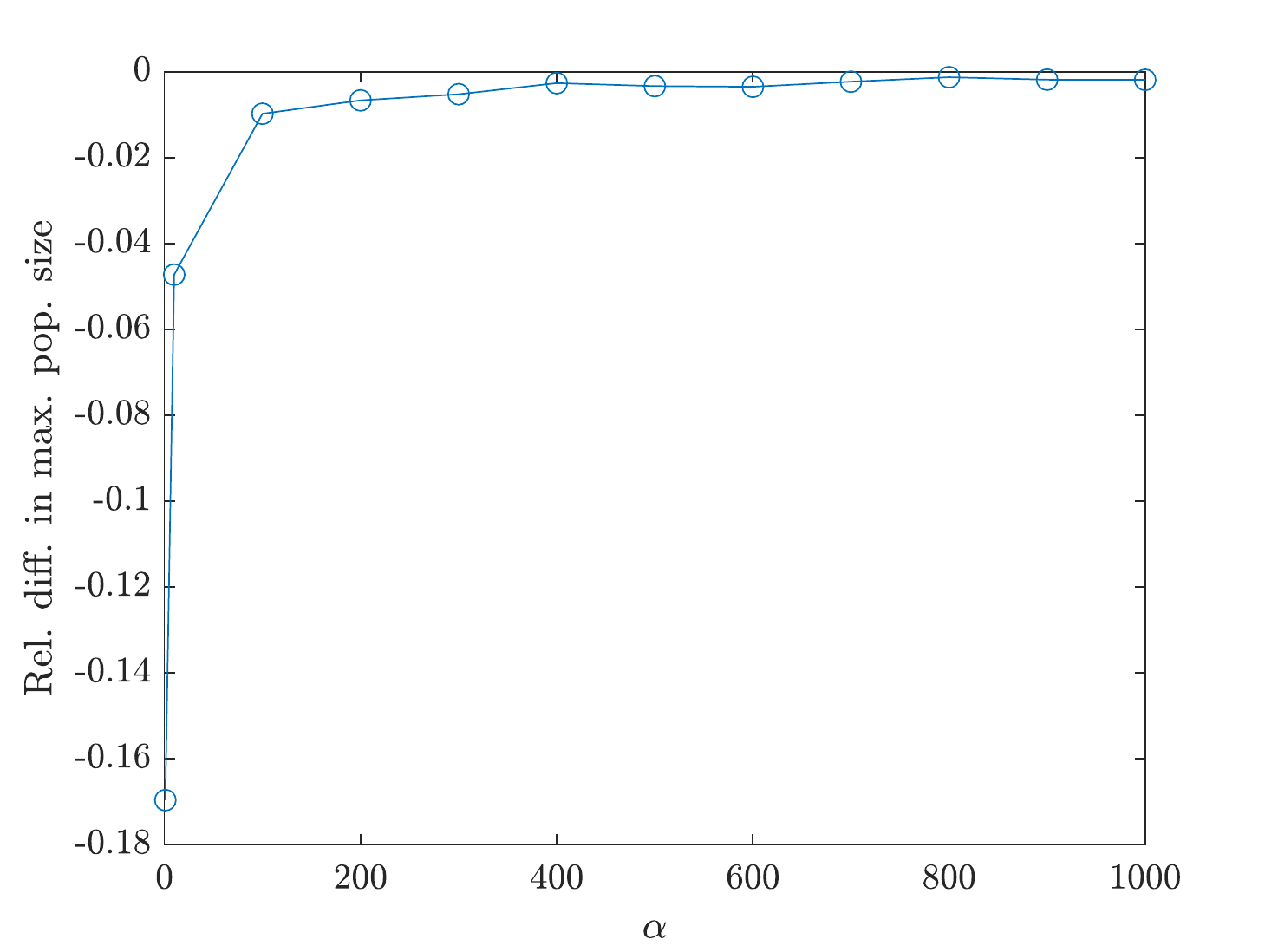}
     \end{subfigure}
        \caption{
          \textbf{Model 2} --- Left: Mean maximum population size as a function of $\alpha$, when $\lambda=1000,\mu = 1$. Right: Relative difference between the average of the simulations and the fluid approximation.}
        \label{fig: model 2: max pop vs alpha}
\end{figure}

From Eqs. (\ref{at max population}) and \eqref{model 2 y1 in y2}, we have:
\begin{align}y_2(t_{max}) &= \alpha \log(\frac{\lambda}{\mu}) \label{model 2: max pop y2}\\
y_1(t_{max}) &= \alpha \log(\frac{\lambda}{\mu}) - \alpha +  \frac{\alpha \mu}{\lambda}e^{\alpha^{-1}}.\label{model 2: max pop y1}
\end{align}
Figures \ref{fig: model 2: max pop vs lambda} and \ref{fig: model 2: max pop vs alpha} show the mean maximum population size as a function of $\lambda$ and $\alpha$, respectively. We see that the fluid approach is accurate, although the fluid curves slightly underestimate the mean maximum population sizes in the stochastic process.

\subsection{Total progeny at extinction}
From Eq. (\ref{model 2 y1 in y2}) and the fact that $y_1(\infty) = 0$, the total progeny at extinction $y_2(\infty)$ satisfies the fixed-point transcendental equation
\begin{equation}
y_2(\infty) = \frac{\alpha \mu}{\lambda}e^{\frac{y_2(\infty)}{\alpha}} - \frac{\alpha \mu}{\lambda}e^{\frac{1}{\alpha}}.
\end{equation}
Although no explicit solution exists, $y_2(\infty)$ can be numerically evaluated (by fixed-point iteration, for example). 

\begin{figure}[H]
     \begin{subfigure}[h]{0.5\linewidth}
         \centering
         \includegraphics[width=\linewidth]{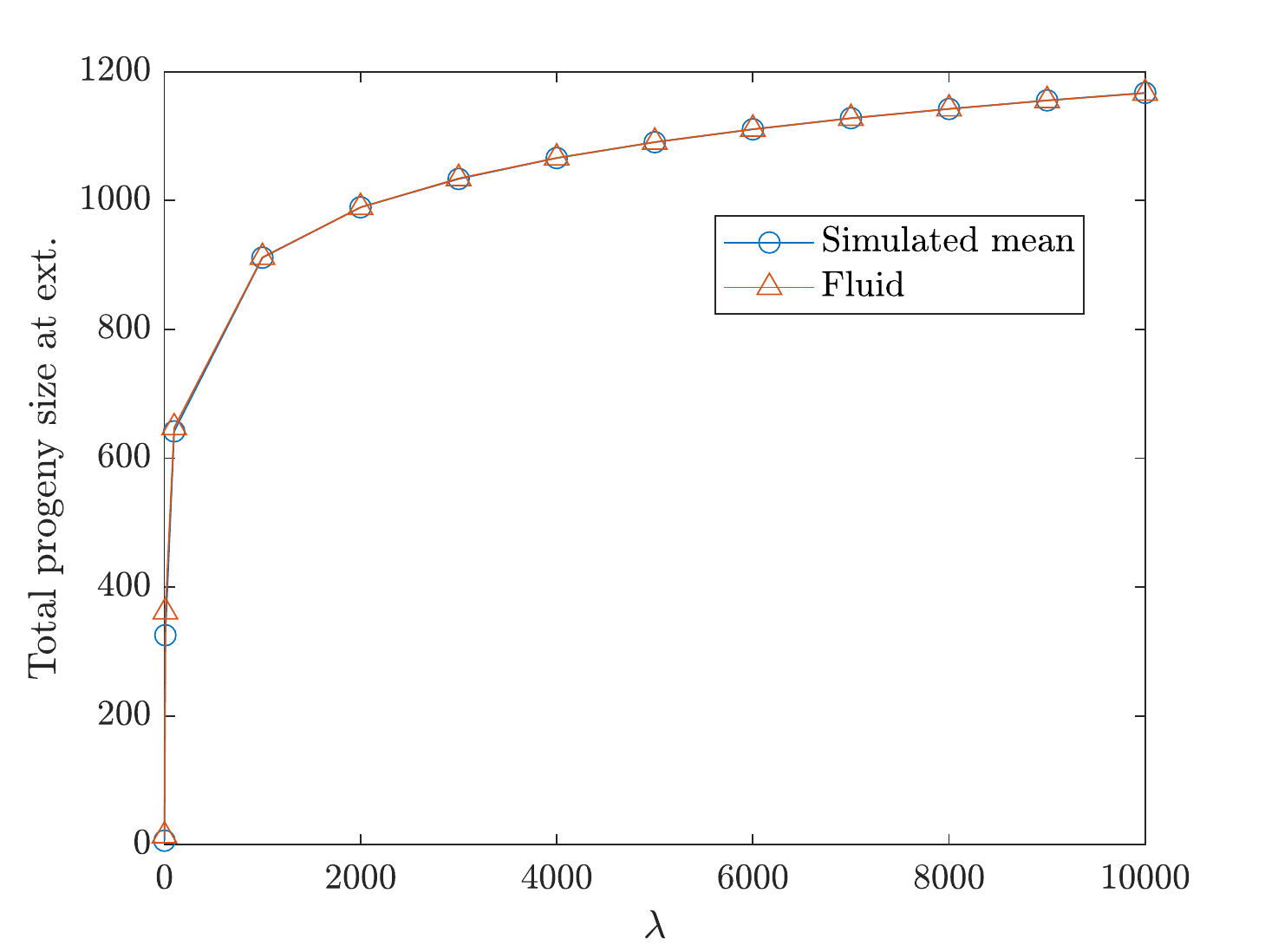}
     \end{subfigure}
     \hfill
     \begin{subfigure}[h]{0.5\linewidth}
         \centering
         \includegraphics[width=\linewidth]{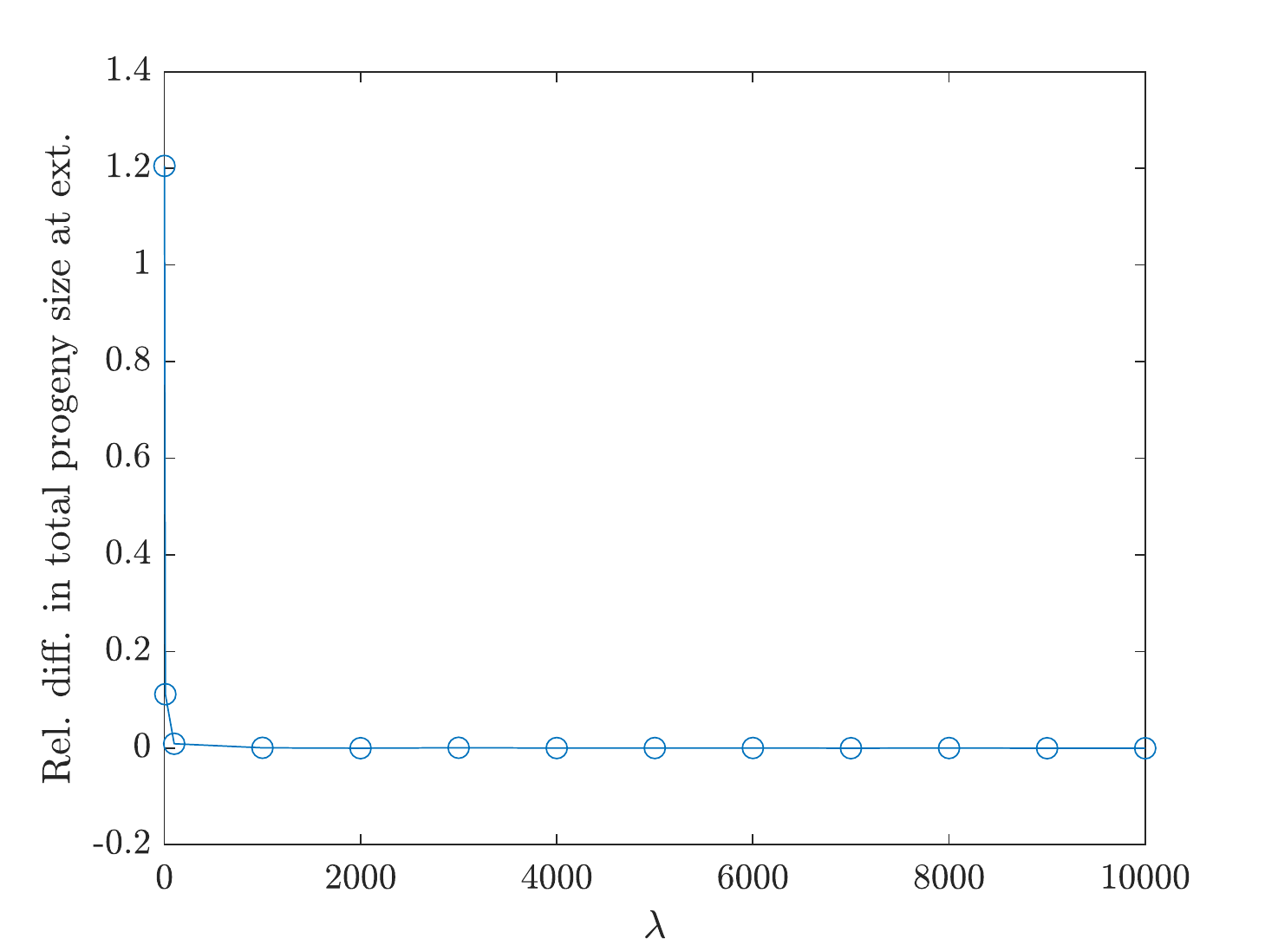}
     \end{subfigure}
        \caption{
          \textbf{Model 2} --- Left: Mean total progeny at extinction as a function of $\lambda$, when $\alpha=100,\mu = 1$. Right: Relative difference between the average of the simulations and the fluid approximation.}
        \label{fig: Model 2: total progeny vs lambda}
\end{figure}

\begin{figure}[H]
     \begin{subfigure}[h]{0.5\linewidth}
         \centering
         \includegraphics[width=\linewidth]{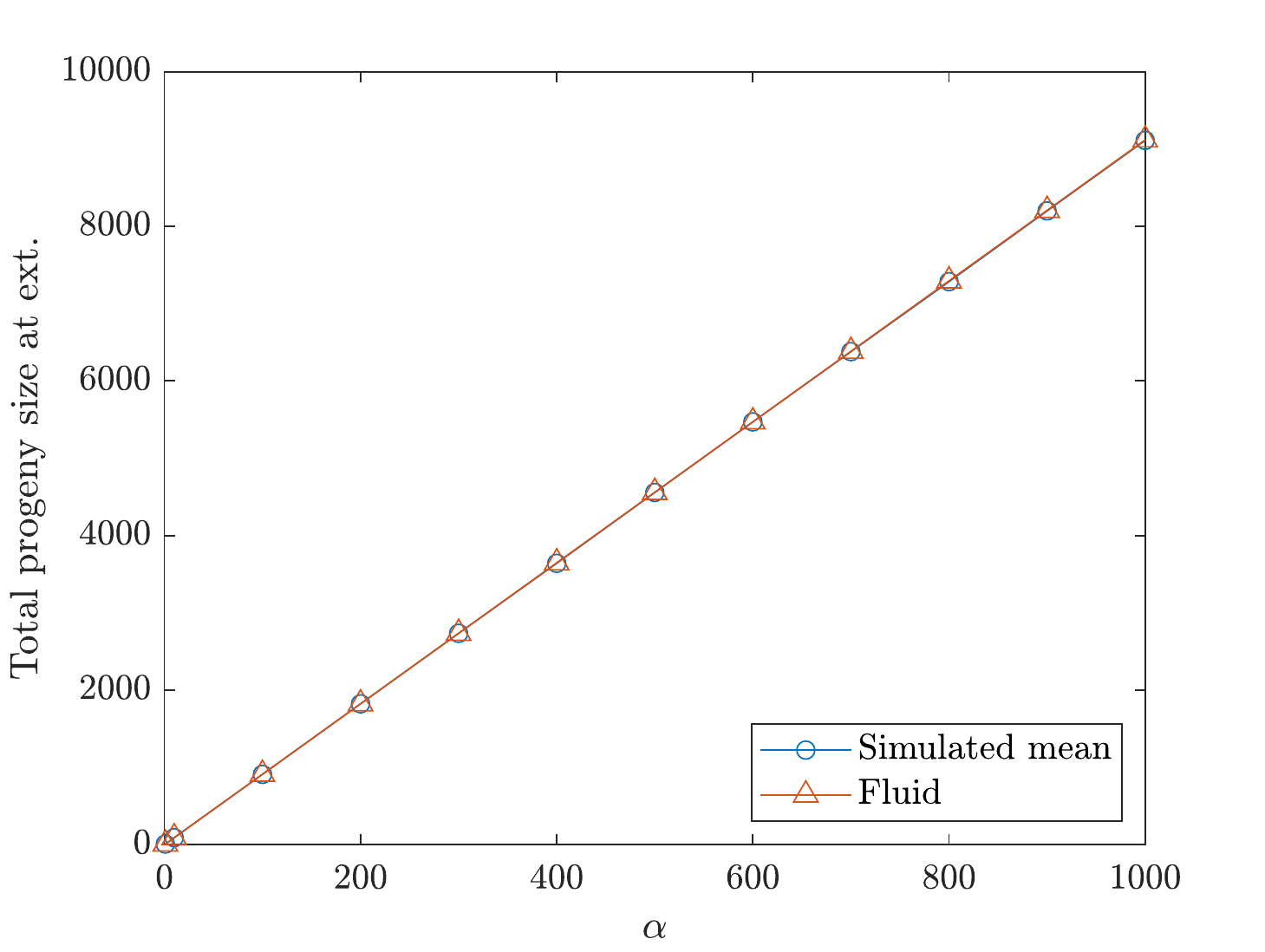}
     \end{subfigure}
     \hfill
     \begin{subfigure}[h]{0.5\linewidth}
         \centering
         \includegraphics[width=\linewidth]{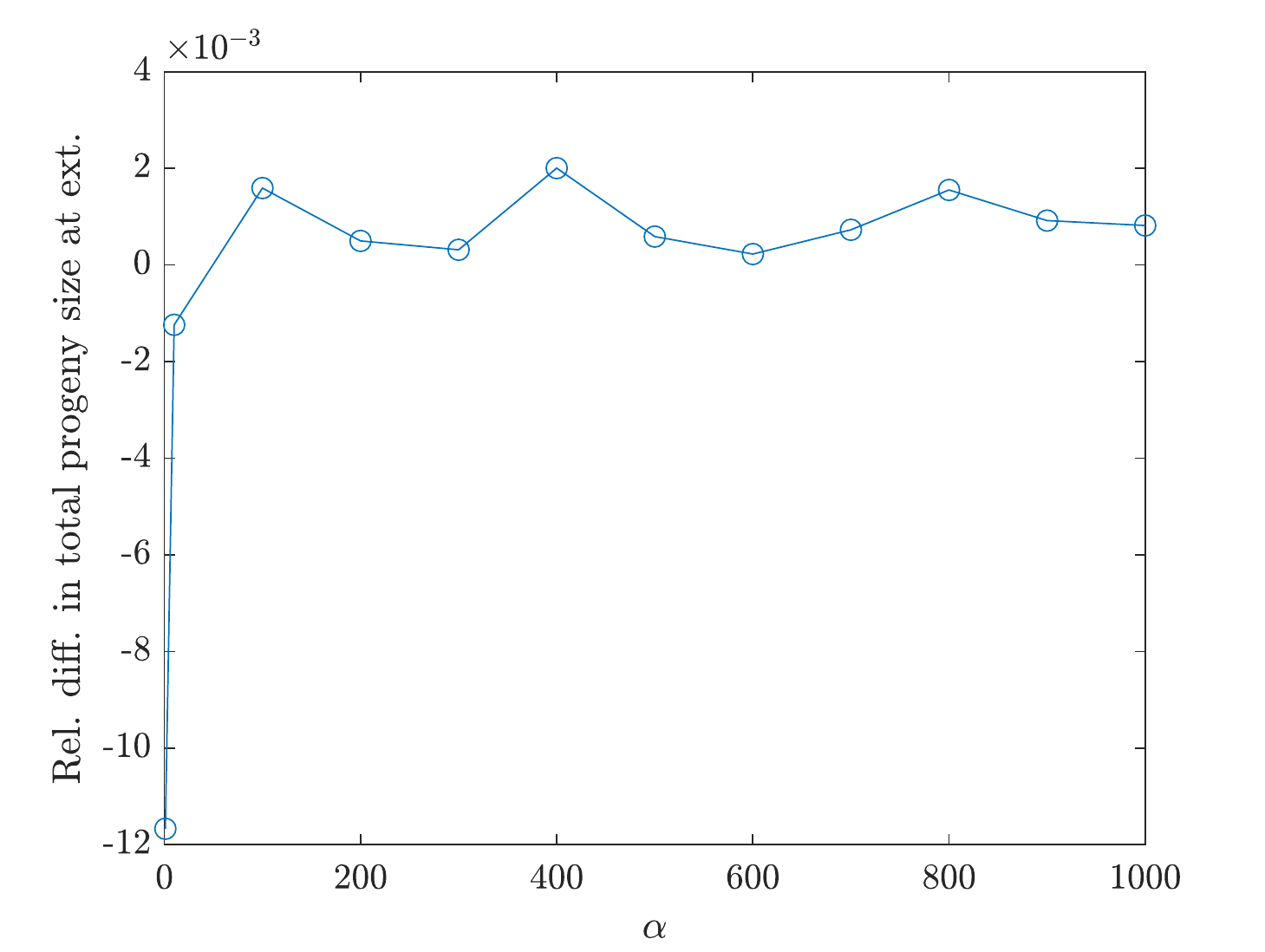}
     \end{subfigure}
        \caption{
          \textbf{Model 2} --- Left: Mean total progeny at extinction as a function of $\alpha$, when $\lambda=1000,\mu = 1$. Right: Relative difference between the average of the simulations and the fluid approximation.}
        \label{fig: Model 2: total progeny vs alpha}
\end{figure}

Figures \ref{fig: Model 2: total progeny vs lambda} and \ref{fig: Model 2: total progeny vs alpha} show the mean total progeny at extinction as a function of $\lambda$ and $\alpha$, respectively. Both figures highlight that the fluid approach is  accurate to approximate the mean total progeny at extinction.

\subsection{Time to reach a given total progeny size}

From Eqs. (\ref{eq2}) and (\ref{model 2 y1 in y2}), we have
\begin{equation}
\frac{dy_2}{dt} = y_1b(y_2) = \left( y_2 - \frac{\alpha \mu}{\lambda}e^{\frac{y_2}{\alpha}} +  \frac{\alpha \mu}{\lambda}e^{\alpha^{-1}} \right) \lambda e^{-\frac{y_2}{\alpha}},
\end{equation}
which implies
\begin{equation} \label{model 2 ode}
 \frac{dt}{dy_2} = \frac{1}{\lambda y_2 e^{-\frac{y_2}{\alpha}} - \alpha \mu + \alpha \mu e^{-\frac{(y_2 -1)}{\alpha}}}.
\end{equation}
Using the initial condition in (\ref{initial condition}), Eq. \eqref{model 2 ode} leads to an integral expression for $t$ as a function of $y_2$,
\begin{equation} \label{model 2: t vs y2}
t(y_2) = \int_1^{y_2} \frac{1}{\lambda u e^{-\frac{u}{\alpha}} - \alpha \mu + \alpha \mu e^{-\frac{(u -1)}{\alpha}}} \,du,
\end{equation}
which, unlike for Model 1, does not have an analytical solution, but can be evaluated numerically.

\subsection{Time to reach the maximum population size}

\begin{figure}[H]
     \begin{subfigure}[h]{0.5\linewidth}
         \centering
         \includegraphics[width=\linewidth]{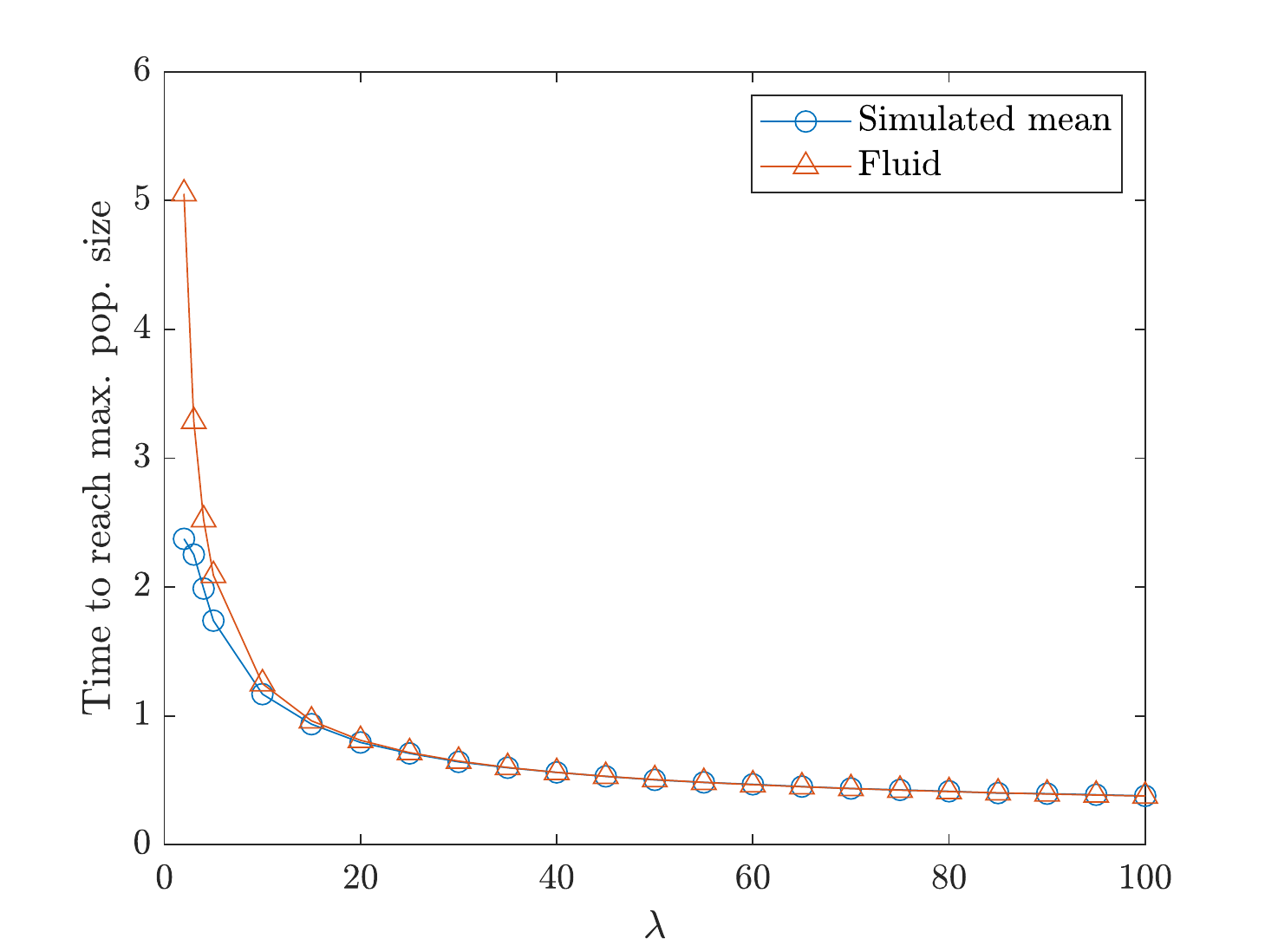}
     \end{subfigure}
     \hfill
     \begin{subfigure}[h]{0.5\linewidth}
         \centering
         \includegraphics[width=\linewidth]{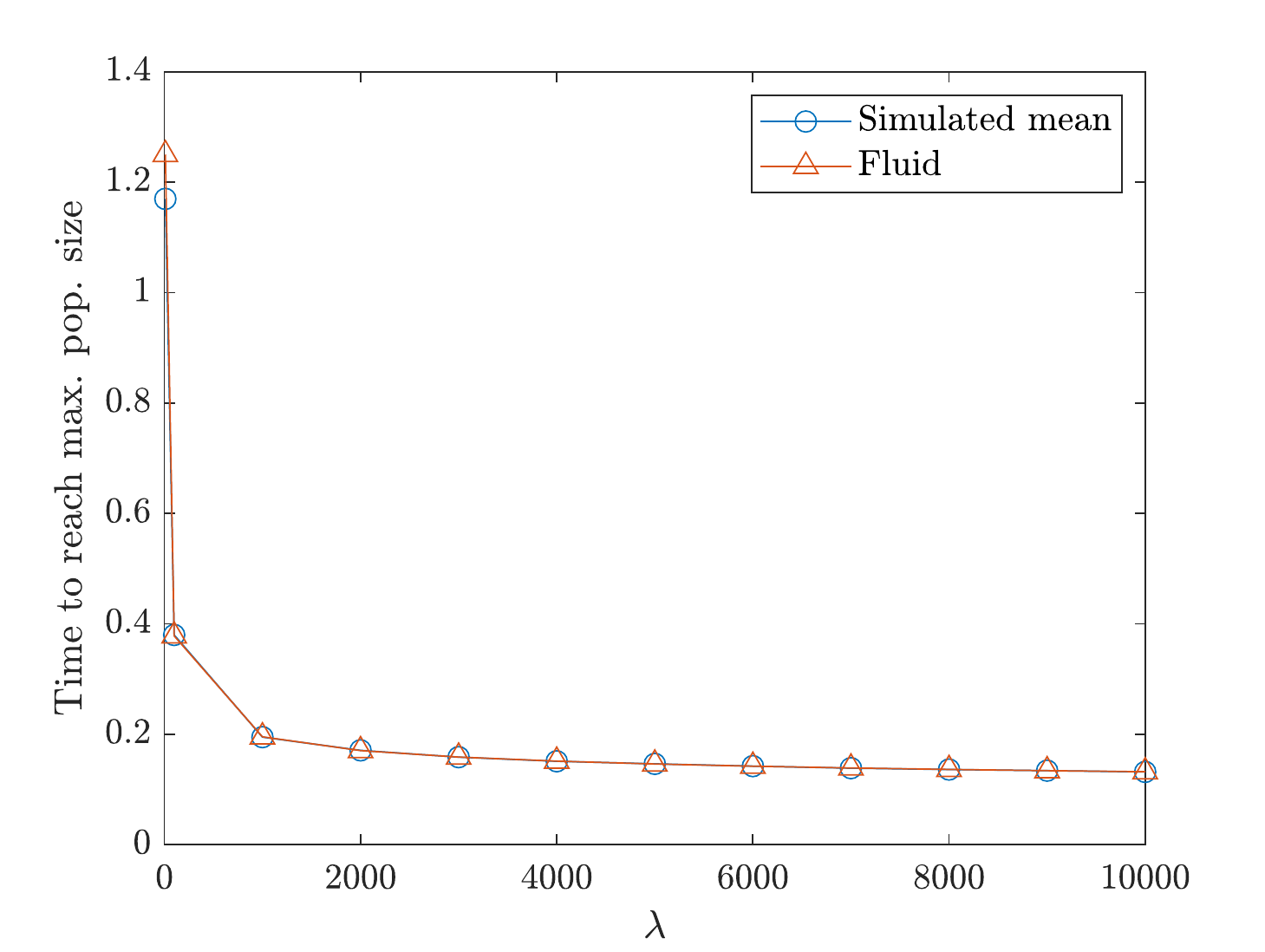}
     \end{subfigure}
        \caption{
          \textbf{Model 2} ---  Mean time to reach the maximum population size as a function of $\lambda$ when $\alpha=100,\mu = 1$. Left: small values of $\lambda$. Right: large values of $\lambda$.}
        \label{fig: model2 tmax lambda}
\end{figure}

By substituting Eq. \eqref{model 2: max pop y2} into Eq. \eqref{model 2: t vs y2}, we obtain the (fluid) time to reach the maximum population size,
\begin{equation}\label{tmax2}
t_{max} = \int_1^{\alpha \log(\lambda/\mu)} \frac{1}{\lambda u e^{-\frac{u}{\alpha}} - \alpha \mu + \alpha \mu e^{-\frac{(u -1)}{\alpha}}} \,du.\end{equation}
Figure \ref{fig: model2 tmax lambda} shows that the fluid approximation is accurate from relatively small values of $\lambda$, and indicates that $t_{max}$ (slowly) decreases to 0 as $\lambda$ increases; this is in contrast with Model 1, where the limit is a strictly positive constant (see Proposition \ref{prop1}). We state this result in the next proposition (whose proof can be found in Section \ref{proofs}):

\begin{prop}\label{prop2}$t_{max} \to 0$ as $\lambda\to\infty$.
\end{prop}
This shows that in this model, the birth rate $b(y_2(t))$ is increasing fast enough with $\lambda$ to ``catch up with'' the logarithmic growth of the maximum population size.

\begin{figure}[H]
\centering
\includegraphics[width=0.8\linewidth]{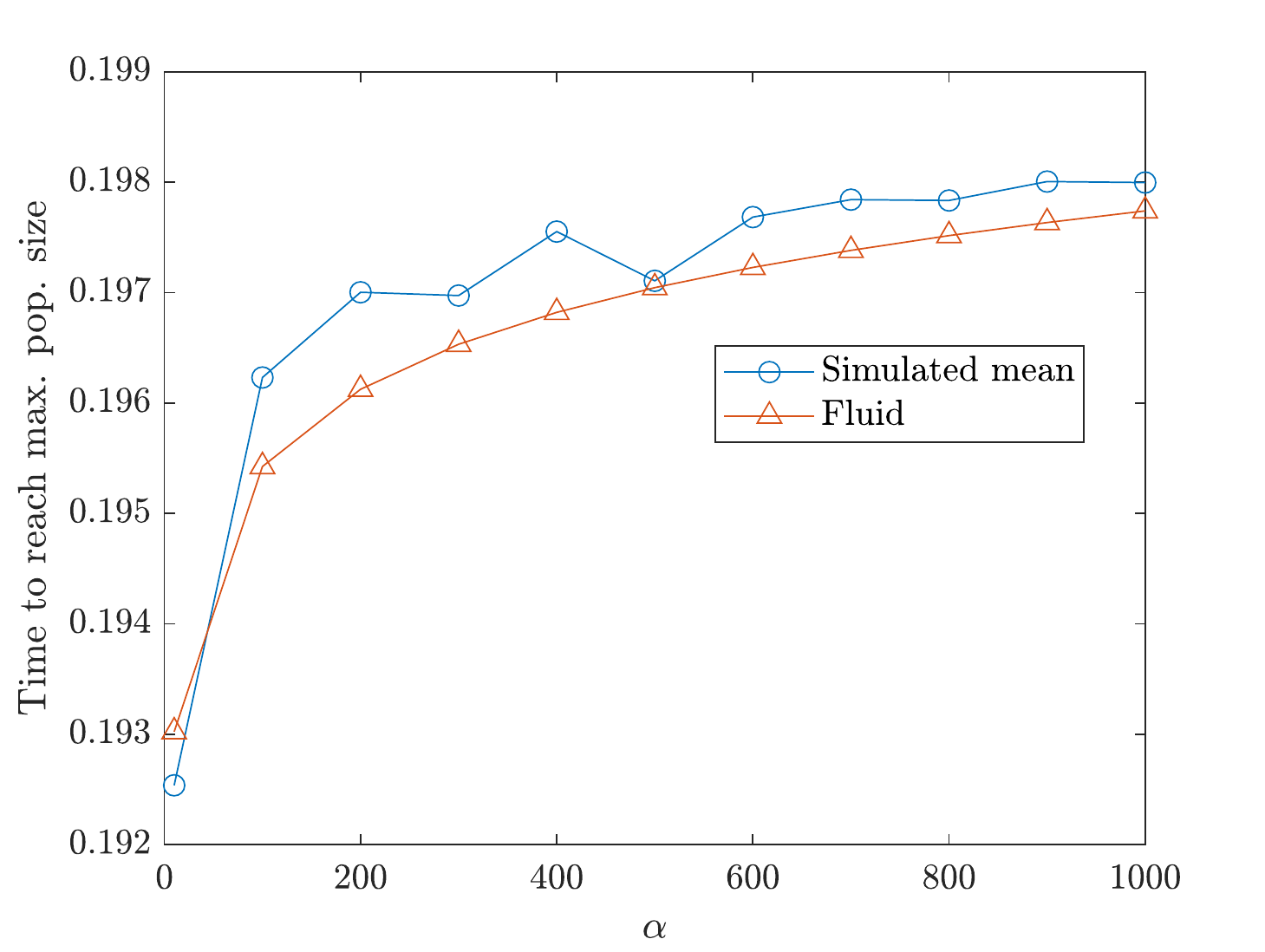}
\caption{\label{fig: model2 tmax alpha}\textbf{Model 2} --- Mean time to reach the maximum population size as a function of $\alpha$ when $\lambda=1000$ and $\mu = 1$. 
}
\end{figure}

Finally, Figure \ref{fig: model2 tmax alpha} shows the mean time to reach the maximum population size as a function of $\alpha$. The figure and further numerical investigations indicate this mean time is increasing with $\alpha$; we did not prove this result analytically.

\subsection{Time until extinction}

The expected time until extinction can  be approximated by $t_{ext}^{(\varepsilon)}$ as in \eqref{text_approx1} and by $t_{ext}^{\star}$ as in \eqref{text_approx2}.
Recall that $t_\varepsilon$ denotes the time at which the (fluid) population reaches size $\varepsilon$ in its descent. 
Let $y_2^*:=y_2(t_\varepsilon)$ be the total progeny at time $t_\varepsilon$; by Eq.~\eqref{model 2 y1 in y2}, $y_2^*$ satisfies
\begin{equation} \label{model 2: y2*, lambda, epsilon}
y_2^* - \frac{\alpha \mu}{\lambda}e^{\frac{y_2^*}{\alpha}} + \frac{\alpha \mu}{\lambda}e^{{\alpha}^{-1}} = \varepsilon,
\end{equation}which can be solved numerically.
We then use Eq. \eqref{model 2: t vs y2} to compute $t_\varepsilon=t(y_2^*)$ and $t(z)$ where $z=y_2(\infty)-1$.

\begin{figure}[h!]
	\begin{subfigure}[h]{0.5\linewidth}
         \centering
         \includegraphics[width=\linewidth]{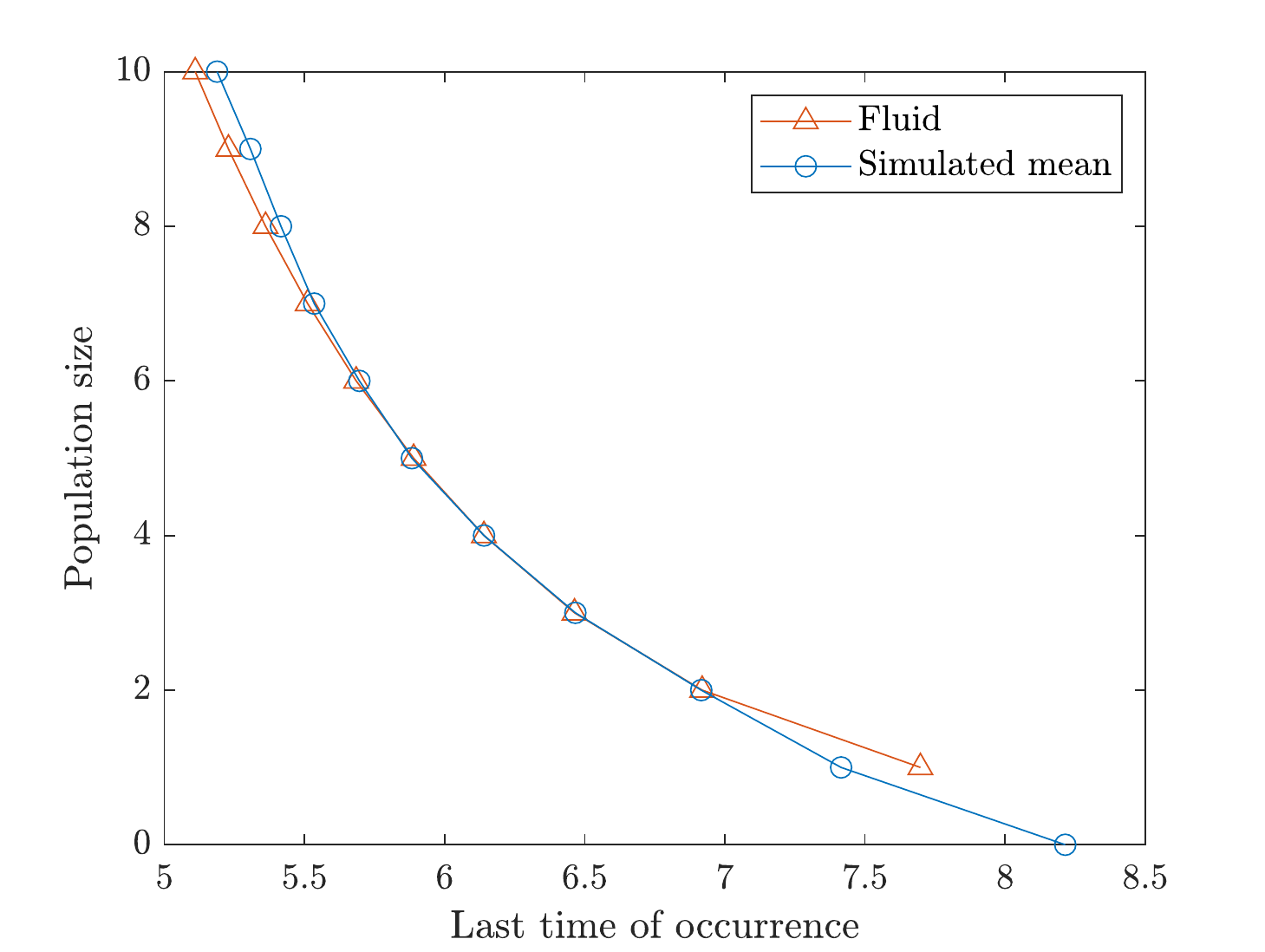}
         \caption{$\lambda = 1000$}
         \label{fig: model 2 small pop size vs t 1000}
     \end{subfigure}
     \hfill
     \begin{subfigure}[h]{0.5\linewidth}
         \centering
         \includegraphics[width=\linewidth]{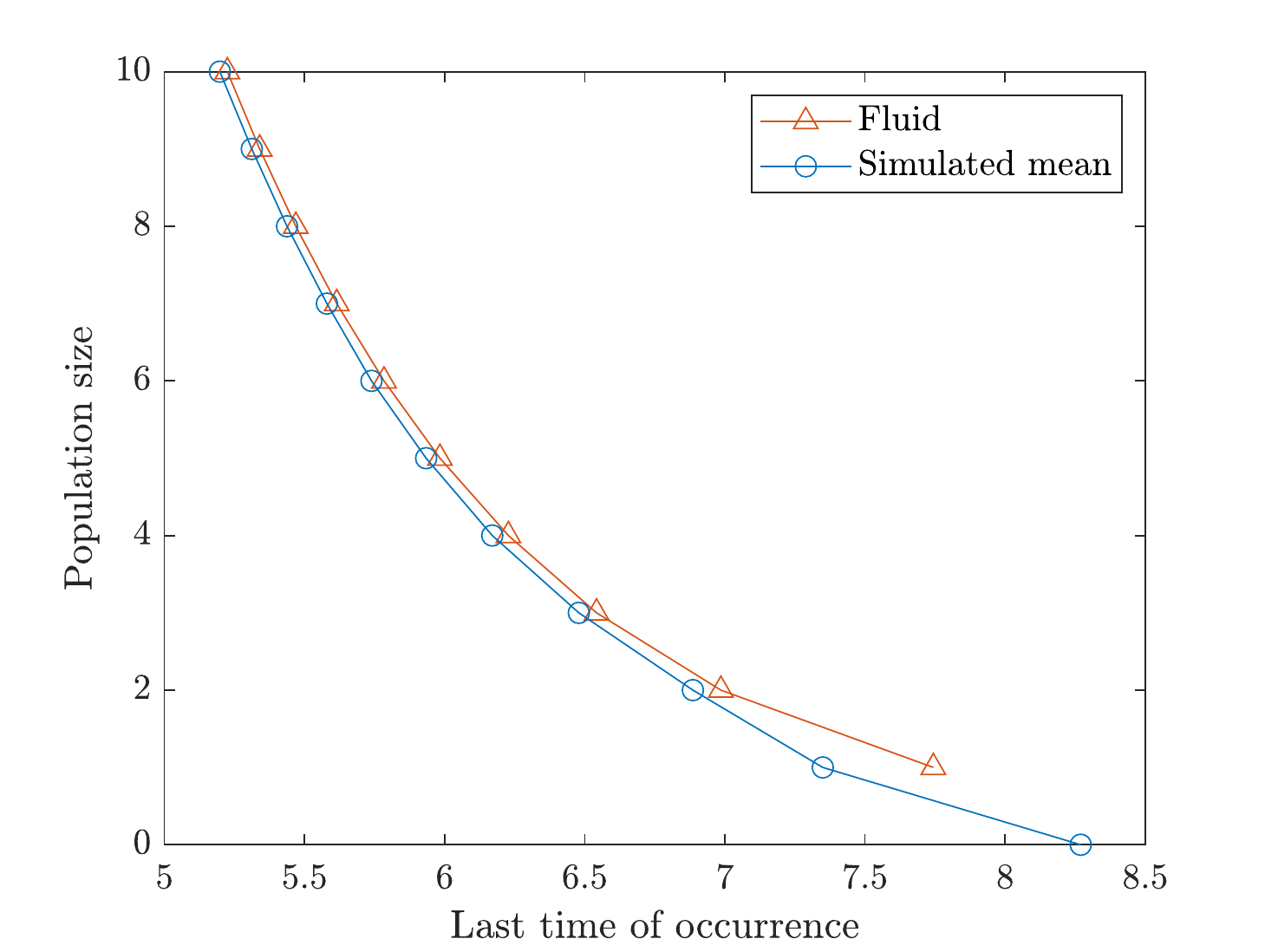}
         \caption{$\lambda = 10000$}
         \label{fig: model 2 small pop size vs t 10000}
     \end{subfigure}

\caption{
\textbf{Model 2} --- Small population sizes and their last occurrence time for $\alpha = 100$.}
\label{fig: model 2 small pop size vs t}
\end{figure}

Similar to Figure \ref{fig: model 1 small pop size vs t}, Figure \ref{fig: model 2 small pop size vs t} plots the average time at which the population reaches some small sizes $\varepsilon$ for the last time before extinction in the stochastic process against $t_{\varepsilon}$. From this figure we expect the fluid approximation $t_{ext}^{(\varepsilon)}$ to be more accurate when $\varepsilon\geq 2$ (like in Model 1).

\begin{figure}[h!]
	\begin{subfigure}[h]{0.5\linewidth}
         \centering
         \includegraphics[width=\linewidth]{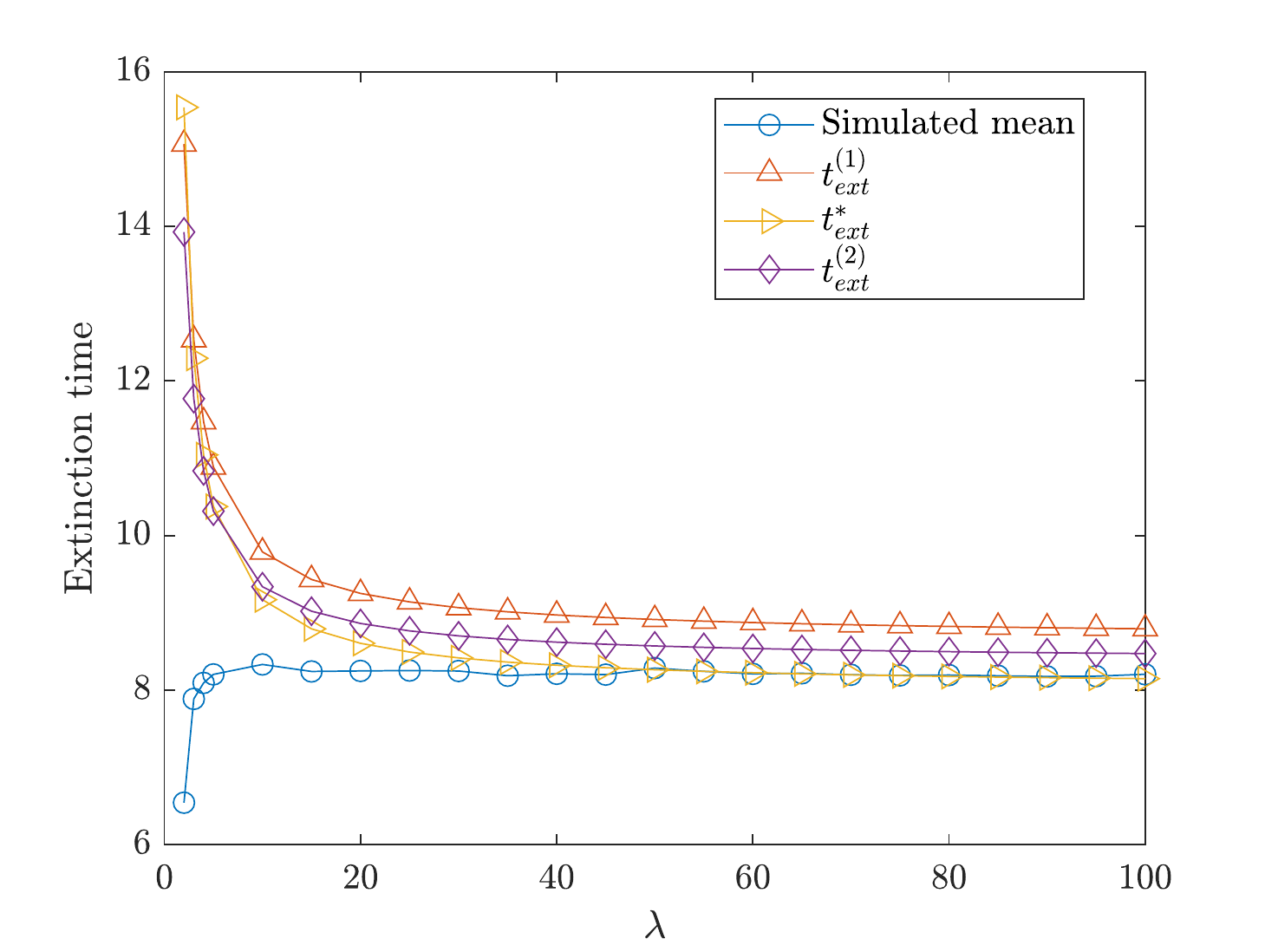}
     \end{subfigure}
     \hfill
     \begin{subfigure}[h]{0.5\linewidth}
         \centering
         \includegraphics[width=\linewidth]{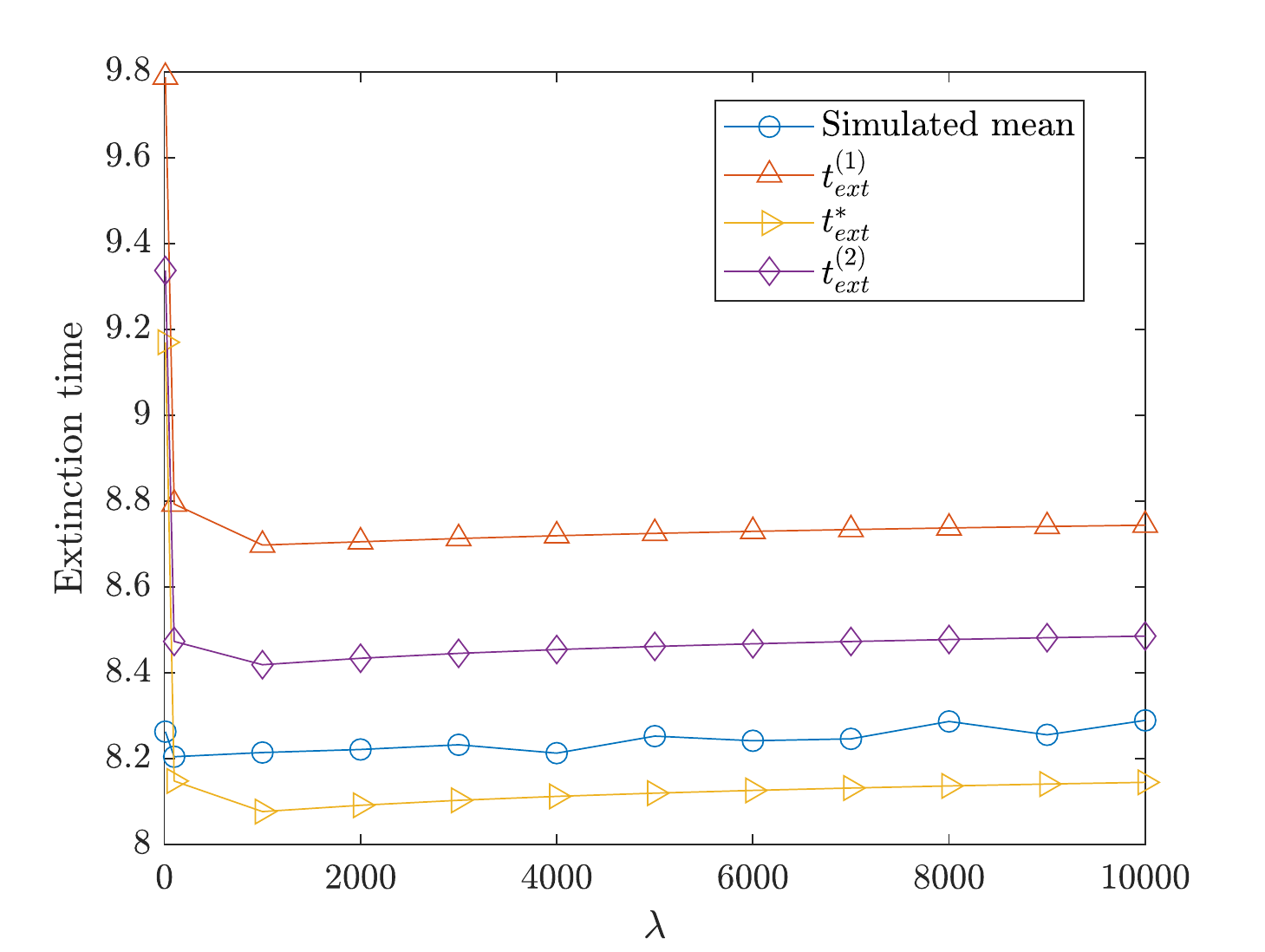}
     \end{subfigure}
\caption{\label{fig: model 2 t vs lambda}
\textbf{Model 2} --- Mean time until extinction obtained with all approximation methods, as a function of $\lambda$ for $\alpha = 100$ and $\mu=1$.  Left: small values of $\lambda$. Right: large values of $\lambda$.}

\end{figure}

\begin{figure}[H]
\centering
\includegraphics[width=0.8\linewidth]{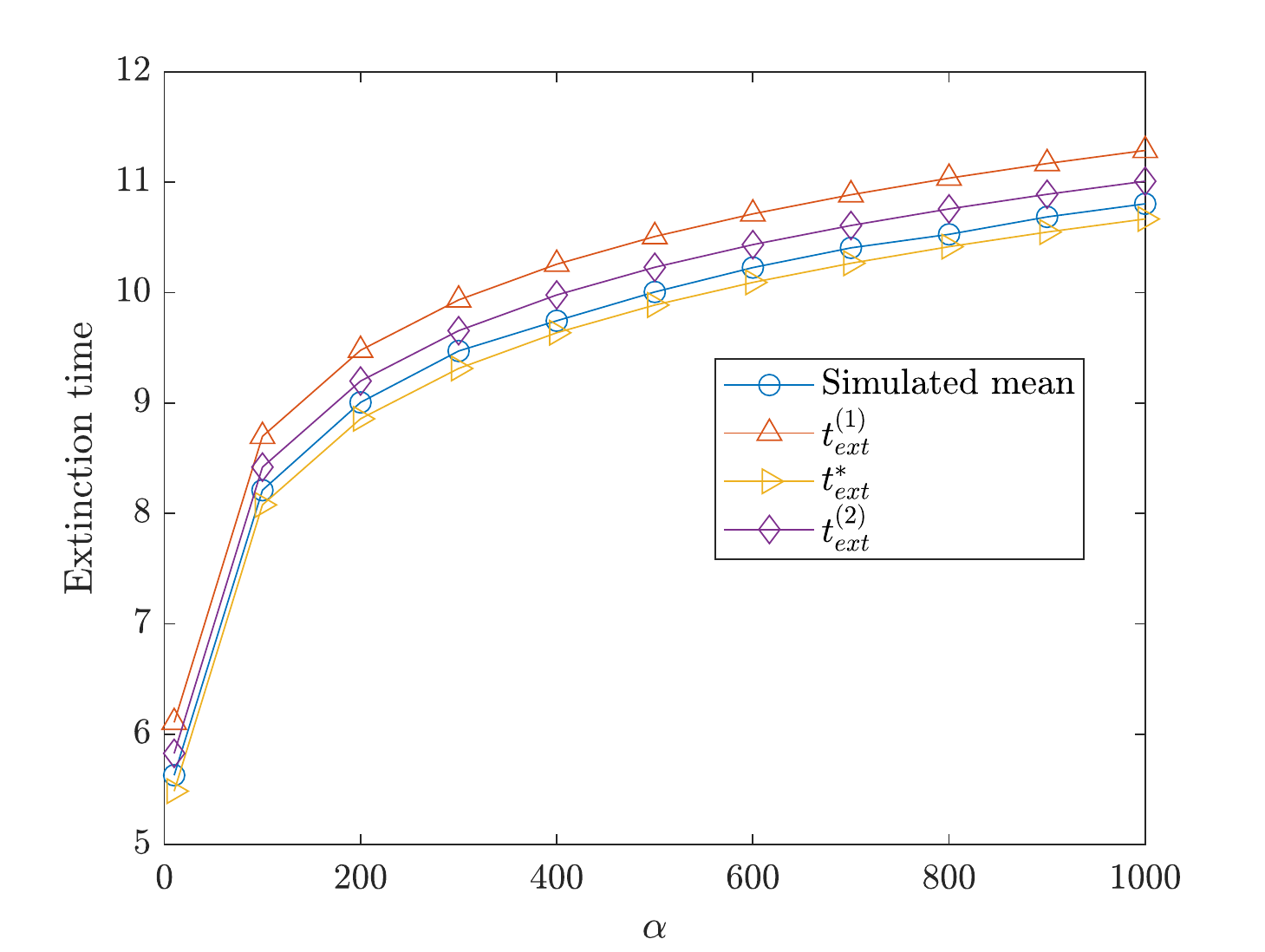}
\caption{\label{fig: model 2 t vs alpha_b}\textbf{Model 2} --- Mean time until extinction obtained with all approximation methods, as a function of $\alpha$ for $\lambda=1000$.
}

\end{figure}

\begin{figure}[h!]
	\begin{subfigure}[h]{0.5\linewidth}
         \centering
         \includegraphics[width=\linewidth]{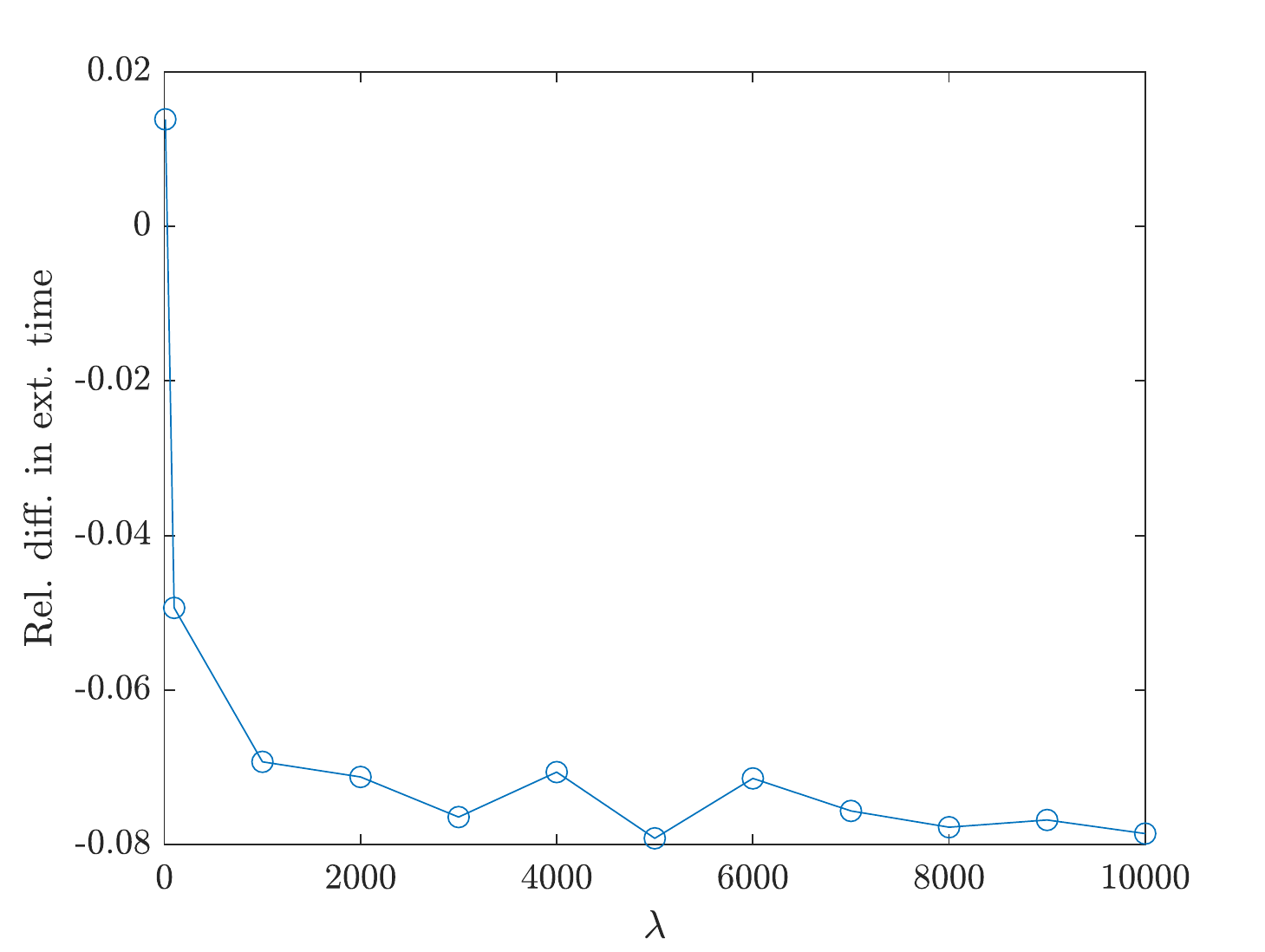}
     \end{subfigure}
     \hfill
     \begin{subfigure}[h]{0.5\linewidth}
         \centering
         \includegraphics[width=\linewidth]{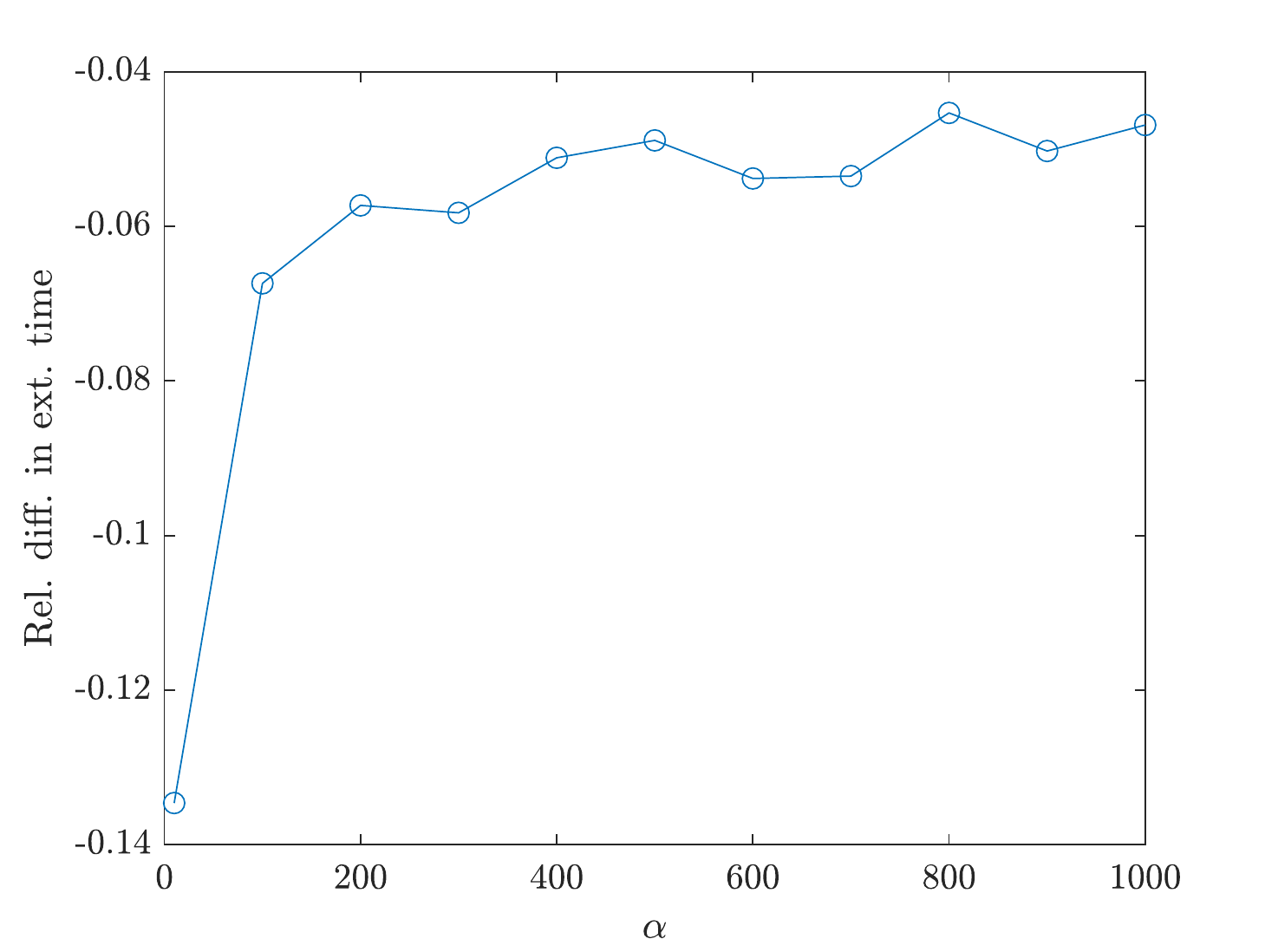}
     \end{subfigure}
\caption{\label{fig: model 2 t re diff}
\textbf{Model 2} --- Relative errors between the average time to extinction based on simulations and the approximation  $t_{ext}^{\star}$ as a function of $\lambda$ (left) and $\alpha$ (right).}

\end{figure}

Figures \ref{fig: model 2 t vs lambda} and \ref{fig: model 2 t vs alpha_b} compare the mean time until extinction obtained by simulations with the approximations $t_{ext}^{(\varepsilon)}$ for $\varepsilon=1,2$, and  $t_{ext}^{\star}$, for values of $\lambda$ and $\alpha$. We observe that the best approximation is given by $t_{ext}^{\star}$ here, followed by $t_{ext}^{(2)}$; the relative error of $t_{ext}^{\star}$ is shown in Figure \ref{fig: model 2 t re diff}. Interestingly, we also see on the left panel of Figure \ref{fig: model 2 t vs lambda} that the time until extinction seems to stabilise to some limiting value as $\lambda$ increases, however the right panel of the figure highlights the potential existence of a minimum in the approximating functions; we further investigate the existence of this minimum at the end of this section.

\begin{figure}[H]
\centering
\includegraphics[width=0.8\linewidth]{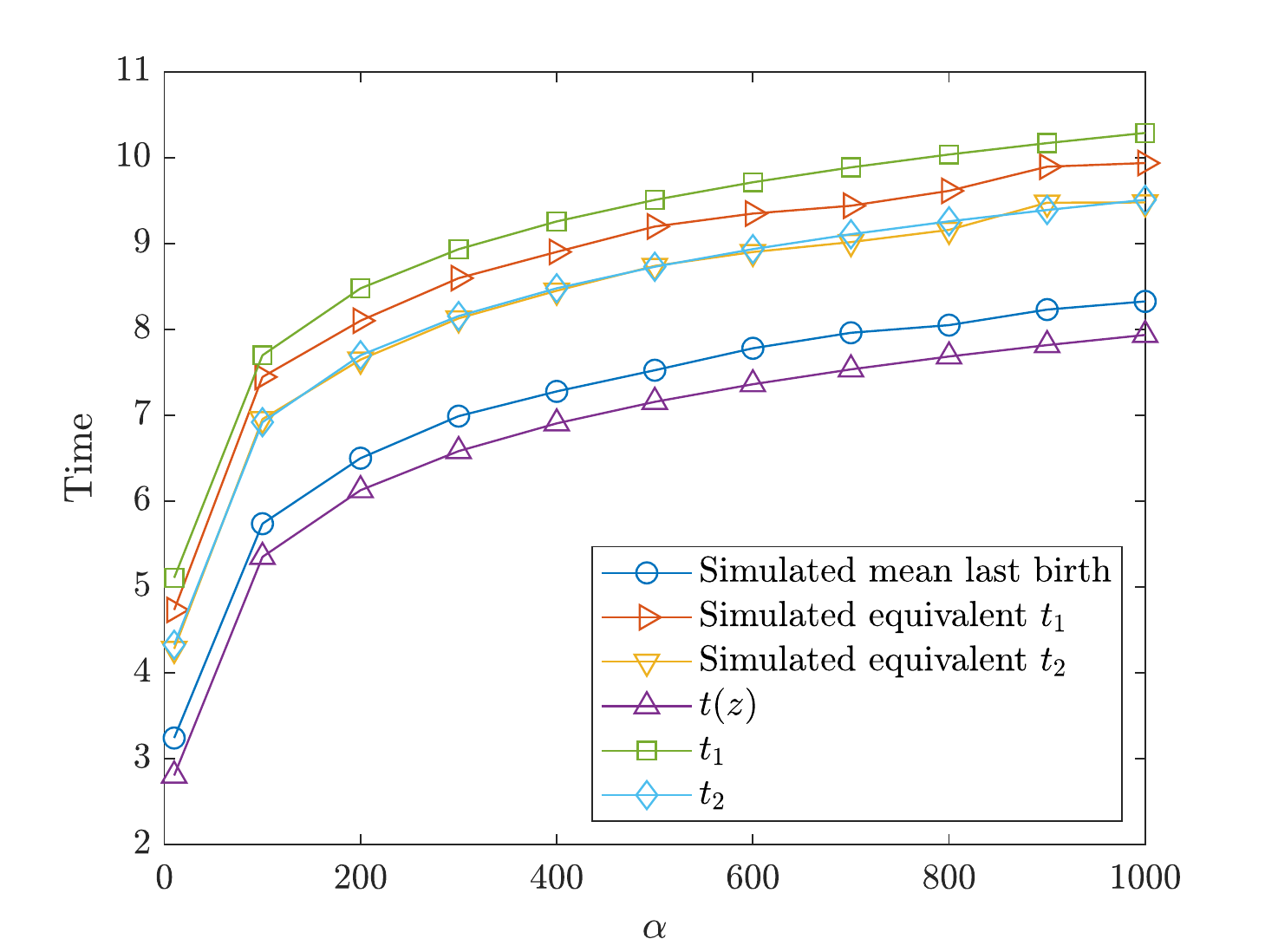}
\caption{\label{fig: model2_time_6_vs_lambda}\textbf{Model 2} --- Time of the last birth, and times to reach population sizes $\varepsilon=1,2$ as a function of $\alpha$ for $\lambda=1000$.
}
\end{figure}

In Figure \ref{fig: model2_time_6_vs_lambda} we compare $t(z)$, the simulated time of last birth, together with $t_{\varepsilon}$ for $\varepsilon=1,2$ and their simulated equivalent, as functions of $\lambda$. This figure shows that $t(z)$ approximates reasonably well the simulated expected time until the last birth, which indicates that our assumption that the process behaves like a pure death process after time $t(z)$ is more appropriate for Model 2 than for Model 1.

In the rest of this section, we show (empirically) the existence of a minimum in $t_\varepsilon$ taken as a function of $\lambda$ (or equivalently, in the approximation $t_{ext}^{(\varepsilon)}$).
Isolating $\lambda$ in Eq. (\ref{model 2: y2*, lambda, epsilon}), we obtain
\begin{equation} \label{model 2: lambda in y2*}
\lambda = \frac{\alpha\mu e^{\frac{y_2^{*}}{\alpha}} - \alpha\mu e^{\frac{1}{\alpha}}}{y_2^{*} - \varepsilon}.
\end{equation}
By substituting Eq. \eqref{model 2: lambda in y2*} into Eq. (\ref{model 2: t vs y2}) with $y_2=y_2^*$, we obtain an expression for $t_{\varepsilon}$ in terms of $y_2^{*}$: 
\begin{equation}
t_{\varepsilon} = \int_1^{y_2^{*}} f(y_2, y_2^{*}) \, dy_2,
\end{equation}
where
\begin{equation}
f(y_2, y_2^{*}) := \dfrac{(y_2^{*} - \varepsilon)}{\alpha\mu}\frac{1}{(y_2^{*} -y_2- \varepsilon)e^{\frac{(1-y_2)}{\alpha}} +y_2 e^{\frac{(y_2^{*}-y_2)}{\alpha}}  - (y_2^{*} - \varepsilon)}.
\end{equation}

As $\lambda \to \infty$, $y_2(t_{max})$ given in Eq. (\ref{model 2: max pop y2}) increases to infinity. If the total progeny at maximum population size increases to infinity, that implies that $ y_2^{*}$ also increases to infinity as $\lambda \to \infty$. More precisely, we show that $ y_2^{*}$ is a strictly monotone function of $\lambda$ (the proof can be found in Section \ref{proofs}):

\begin{lem}\label{lem1}If $\varepsilon\geq 1$ then $d y_2^{*}/d\lambda>0$.
\end{lem}

As a consequence of Lemma \ref{lem1}, if $t_{\varepsilon}$ taken as a function of $y_2^{*}$ exhibits a minimum, then $t_{\varepsilon}$ taken as a function of $\lambda$ also exhibits a minimum. This change of scale allows us to visualise this minimum in a much clearer way: Figure \ref{fig: model 2 t vs y_2^{*}} shows the existence of the minimum in $t_{\varepsilon}$ taken as a function of $y_2^{*}$ for $\varepsilon=1$ and $\alpha=100$. In this case, the minimum can be evaluated numerically and is reached when $y_2^*\approx888.44$,  which according to Eq. \eqref{model 2: lambda in y2*} corresponds to $\lambda\approx 813.3$.

\begin{figure}[h!]
\centering
\includegraphics[width=0.9\linewidth]{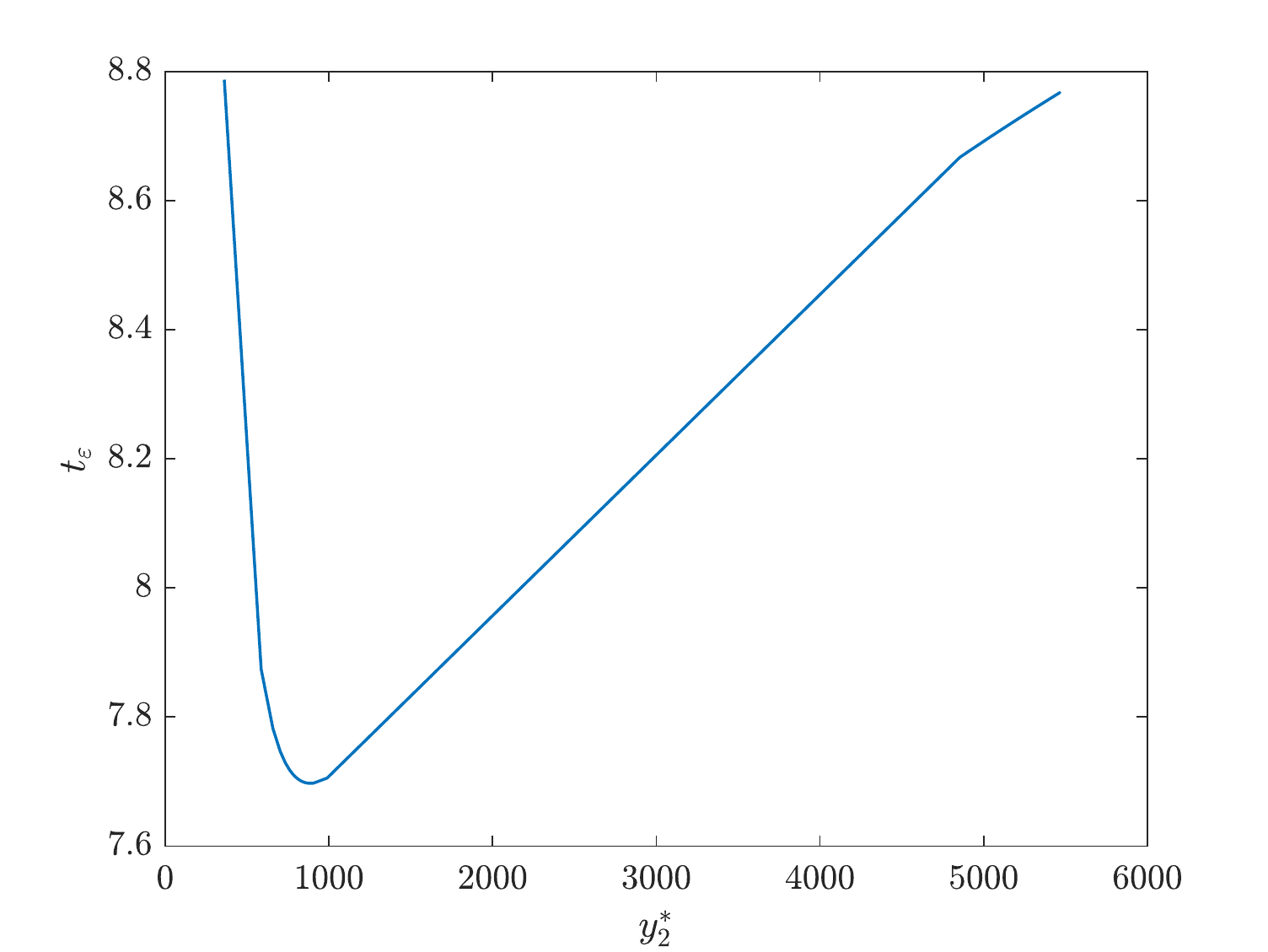}
\caption{\textbf{Model 2} --- $t_{\varepsilon}$ as a function $ y_2^{*}$ for $\varepsilon=1$}
\label{fig: model 2 t vs y_2^{*}}
\end{figure}

\section{Proofs}\label{proofs}

\subsection*{Proof of Proposition \ref{prop1}} 

Using \eqref{model 1 t with y2} with $y_2(t_{max}) = {\lambda}/{\mu}$ and $c^2=1+\frac{\lambda^2}{\mu^2}$, we get
\begin{equation*}
\begin{split}
t({\lambda}/{\mu}) 
&=\frac{1}{\mu}\left[-\log(c^2) + \log(c^2-\frac{\lambda^2}{\mu^2}+\frac{2\lambda}{\mu}-1) + \frac{\lambda}{\mu c}\log\left(\frac{c+\frac{\lambda}{\mu}-1}{c-\frac{\lambda}{\mu}+1}\right)\right]\\
&=\frac{1}{\mu}\left[\log(2)+\log\left(\frac{\frac{\lambda}{\mu}}{1+\frac{\lambda^2}{\mu^2}}\right)  + \frac{\lambda}{\mu \sqrt{1+\frac{\lambda^2}{\mu^2}}}\log\left(\frac{\sqrt{1+\frac{\lambda^2}{\mu^2}}+\frac{\lambda}{\mu}-1}{\sqrt{1+\frac{\lambda^2}{\mu^2}}-\frac{\lambda}{\mu}+1}\right)\right]\\
&=\frac{1}{\mu}\left[\log(2)+(1-\frac{\lambda}{\sqrt{\mu^2+\lambda^2}})\log(\frac{\mu \lambda}{\mu^2+\lambda^2}) \right.\\
& \left.+\frac{\lambda}{\sqrt{\mu^2+\lambda^2}}\log\left(\frac{\frac{1}{\mu} \sqrt{\frac{1}{\lambda^2}+\frac{1}{\mu^2}}+\frac{1}{\mu^2}-\frac{1}{\lambda \mu}}{\frac{1}{\lambda^2}+\frac{1}{\mu^2}}\cdot \frac{1}{\sqrt{1+\frac{\lambda^2}{\mu^2}}-\frac{\lambda}{\mu}+1}\right)\right] \\&\to \frac{2}{\mu}\log(2) \text{ as } \lambda\to\infty.
\end{split}
\end{equation*}
\subsection*{Proof of Proposition \ref{prop2}} 

By the change of variable $u=-\alpha\log x$ in Eq. \eqref{tmax2}, we obtain
\begin{equation}\label{tmax2}
t_{max} = \int_{\mu/\lambda}^{\exp(-1/\alpha)} \frac{1}{-\lambda x^2 \log x +\mu x(x e^{\frac{1}{\alpha}}-1) } \,dx.\end{equation}
For $\lambda$ that is large enough, we have $\mu/\lambda<\exp(-1/\alpha)$, and $-\lambda x^2 \log x +\mu x(x e^{\frac{1}{\alpha}}-1)>0$ for $0<\mu/\lambda\leq x\leq \exp(-1/\alpha)<1$. 

We first show that
\begin{equation} \label{bound}
0\leq t_{max} \leq \int_{\mu/\lambda}^{\exp(-1/\alpha)} g(x) \,dx,\; \text{where}\; g(x)=\frac{M}{-\lambda x^2 \log x},
\end{equation}
and $M$ is a sufficiently large constant.
Let 
$$f(x) = \frac{1}{-\lambda x^2 \log x +\mu x(x e^{\frac{1}{\alpha}}-1) };$$then we have
\begin{equation}\label{g-f}
g(x)-f(x)=\frac{M \mu e^{\frac{1}{\alpha}} x^2 + (M-1) x (-\lambda x \log x - \frac{M}{M-1} \mu ) }{(-\lambda x^2 \log x + \mu x (xe^{\frac{1}{\alpha}} -1))(-\lambda x^2 \log x)}.
\end{equation}
For $\mu/\lambda\leq x\leq \exp(-1/\alpha)$, the denominator in \eqref{g-f} is always positive. For the numerator, let $a(x) = -\lambda x \log x - \frac{M}{M-1} \mu$; we can check that $a(x)$ has one global maximum in the interval of interest, and is positive at the boundaries of the interval, that is, $a({\mu}/{\lambda})>0$ and $a(\exp(-1/\alpha)) > 0$. Therefore, $a(x)>0$, and $g(x)>f(x)>0$ over $[\mu/\lambda, \exp(-1/\alpha)]$, so \eqref{bound} follows.
\hfill \break

Next we show that $\lim_{\lambda\to\infty} \int_{\mu/\lambda}^{\exp(-1/\alpha)} g(x) \,dx=0$.
We have
\begin{equation}
\int_{\mu/\lambda}^{\exp(-1/\alpha)} g(x) \,dx =\frac{M}{\lambda}\int_{\exp(1/\alpha)}^{\lambda /\mu} \frac{1}{\log u} \,du,
\end{equation}
where the right-hand-side involves a logarithmic integral, which has a well-known asymptotic behaviour:
\begin{equation}
\frac{M}{\lambda}\int_{\exp(1/\alpha)}^{\lambda /\mu} \frac{1}{\log u} \,du \sim \frac{M}{\mu \log \frac{\lambda}{\mu} \sum_{k=0}^\infty \frac{k!}{(\log \frac{\lambda}{\mu})^k}} \to 0,
\end{equation}
as $\lambda \to \infty$. Therefore, by \eqref{bound}, the proposition holds.

\subsection*{Proof of Lemma \ref{lem1}}

We consider the cases when $\varepsilon \leq y_1(t_{max})$. For $y_2^* > \max(1, \varepsilon)$, recall that $y_2^*$ satisfies Eq. \eqref{model 2: y2*, lambda, epsilon}. By differentiating \eqref{model 2: y2*, lambda, epsilon}
with respect to $\lambda$, we obtain 
\begin{equation} \label{model 2: y2*, lambda, epsilon, derivative}
\frac{dy_2^*}{d\lambda}=\frac{\frac{\alpha \mu}{\lambda} (e^{\frac{y_2^*}{\alpha}}-e^{\frac{-1}{\alpha}})}{\lambda - \mu e^{\frac{y_2^*}{\alpha}}}.
\end{equation}
The numerator is always negative, and by Eq. \eqref{model 2: lambda in y2*}, we can express the denominator as
\begin{equation}\label{model 2: y2*, lambda, epsilon, derivative, denominator}
\frac{\mu}{y_2^* - \varepsilon}\, f(y_2^*), \quad\text{where}\quad f(y_2^*)=\alpha e^{\frac{y_2^*}{\alpha}} - \alpha e^{\frac{1}{\alpha}} -(y_2^* - \varepsilon) e^{\frac{y_2^*}{\alpha}}.
\end{equation}
The function $f(y_2^*)$ is  decreasing (its derivative is strictly negative). Therefore, for $0<\varepsilon < y_1(t_{max}) = \alpha \log(\frac{\lambda}{\mu}) - \alpha +  \frac{\alpha \mu}{\lambda}e^{\alpha^{-1}}$, by Eq. \eqref{model 2: max pop y2}, $y_2^*\geq y_2(t_{max})= \alpha \log (\frac{\lambda}{\mu})$, thus
\begin{equation}
f(y_2^*) \leq f(\alpha \log (\frac{\lambda}{\mu})) = \frac{\lambda}{\mu} (-\alpha \log(\frac{\lambda}{\mu}) + \alpha -  \frac{\alpha \mu}{\lambda}e^{\alpha^{-1}} + \varepsilon) < 0,
\end{equation}
and from \eqref{model 2: y2*, lambda, epsilon, derivative} and \eqref{model 2: y2*, lambda, epsilon, derivative, denominator}, $d y_2^{*}/d\lambda>0$.

Finally, if $\varepsilon = y_1(t_{max})$, then from Eq. \eqref{model 2: max pop y2}, $y_2^*=\alpha \log (\frac{\lambda}{\mu})$, and the lemma follows as well.

\section*{Acknowledgments}

Sophie Hautphenne would like to thank the Australian Research Council (ARC) for support through her Discovery Project DP200101281.

\FloatBarrier


\begin{thebibliography}{10}

\bibitem{bruss2015resource}
F~Thomas Bruss and Mitia Duerinckx.
\newblock Resource dependent branching processes and the envelope of societies.
\newblock {\em The Annals of Applied Probability}, 25(1):324--372, 2015.

\bibitem{darling2002fluid}
RWR Darling.
\newblock Fluid limits of pure jump {M}arkov processes: a practical guide.
\newblock {\em arXiv preprint math/0210109}, 2002.

\bibitem{hamza2016establishment}
K.~Hamza, P.~Jagers, and F.C. Klebaner.
\newblock On the establishment, persistence, and inevitable extinction of
  populations.
\newblock {\em Journal of Mathematical Biology}, 72(4):797--820, 2016.

\bibitem{jagers1992stabilities}
P.~Jagers.
\newblock Stabilities and instabilities in population dynamics.
\newblock {\em Journal of Applied Probability}, 29(4):770--780, 1992.

\bibitem{jk11}
P.~Jagers and F.C. Klebaner.
\newblock Population-size-dependent, age-structured branching processes linger
  around their carrying capacity.
\newblock {\em Journal of Applied Probability}, 48(A):249--260, 2011.

\bibitem{jagers2020populations}
Peter Jagers and Sergei Zuyev.
\newblock Populations in environments with a soft carrying capacity are
  eventually extinct.
\newblock {\em Journal of Mathematical Biology}, 81(3):845--851, 2020.

\bibitem{kurtz1970solutions}
Thomas~G Kurtz.
\newblock Solutions of ordinary differential equations as limits of pure jump
  {M}arkov processes.
\newblock {\em Journal of Applied Probability}, 7(1):49--58, 1970.

\bibitem{kurtz1971limit}
Thomas~G Kurtz.
\newblock Limit theorems for sequences of jump {M}arkov processes approximating
  ordinary differential processes.
\newblock {\em Journal of Applied Probability}, 8(2):344--356, 1971.

\bibitem{sevast1974controlled}
Boris~Alexandrovich Sevast’yanov and Andrei~Mikhailovich Zubkov.
\newblock Controlled branching processes.
\newblock {\em Theory of Probability \& Its Applications}, 19(1):14--24, 1974.

\bibitem{velasco2017controlled}
Miguel~Gonz{\'a}lez Velasco, In{\'e}s Mar{\'\i}a Del~Puerto Garc{\'\i}a, and
  George~Petrov Yanev.
\newblock {\em Controlled branching processes}.
\newblock John Wiley \& Sons, 2017.

\end{thebibliography}
\end{document}